\documentclass[reqno,oneside, notitlepage]{amsart}

\usepackage{geometry}
\usepackage[dvipsnames]{xcolor}
\numberwithin{equation}{section}

\usepackage{tikz}
\usetikzlibrary{arrows.meta, positioning}
\usetikzlibrary{arrows.meta, automata,  
                positioning,
                quotes}

\usepackage{amsmath}
\usepackage{amssymb, amsthm}
\usepackage{mathtools}
\usepackage{mathrsfs}
\usepackage{esint}
\usepackage{bm}

\usepackage{graphicx}
\usepackage{tikz}
\usepackage{subcaption}
\usetikzlibrary{patterns} 
\usepackage{tikz-3dplot}
\tikzset{tdplot_fillconestyle/.style={fill=gray!10,opacity=0.3}}
\tdplotsetmaincoords{70}{25}

\usepackage[backend=biber,style=numeric,maxnames=99,minnames=1]{biblatex}
\AtEveryBibitem{%
  \clearfield{url}%
  \clearfield{issn}%
  \clearfield{isbn}%
} 
\DeclareFieldFormat[article]{volume}{\mkbibbold{#1}}  
\DeclareFieldFormat[article]{number}{\mkbibparens{#1}}  

\renewbibmacro*{volume+number+eid}{
  \printfield{volume}%
  \setunit{}%
  \printfield{number}%
  \setunit{\addcomma\space}%
  \printfield{eid}}

\usepackage{stackengine}
\usepackage{IEEEtrantools}

\newtheorem{theorem}{Theorem}[section]
\newtheorem{lemma}[theorem]{Lemma}
\newtheorem{proposition}[theorem]{Proposition}
\newtheorem{corollary}[theorem]{Corollary}

\theoremstyle{definition}
 \newtheorem{definition}[theorem]{Definition}

\theoremstyle{remark}
 \newtheorem{remark}[theorem]{Remark}

\newcommand{\Energy}[2]{E_{\epsilon,#1}^{(#2)}(u;\chi)}
\newcommand{\Energyinf}[3]{E_{\epsilon,#1}^{(#2)}(\mathcal{A}_F^{#3})}

\newcommand\subsetsim{\mathrel{%
  \ooalign{\raise0.2ex\hbox{$\subset$}\cr\hidewidth\raise-0.8ex\hbox{\scalebox{0.9}{$\sim$}}\hidewidth\cr}}}

\DeclareMathOperator{\dir}{dir}
\DeclareMathOperator{\per}{per}
\newcommand{\Eelinf}[2]{E_{\rm el, #1}^{(#2)}  (\chi; F)}
\newcommand{\Eel}[2]{E_{\rm el, #1}^{(#2)} (u; \chi)}

\DeclareMathOperator{\dist}{dist}
\DeclareMathOperator{\F}{\mathcal{F}}

\newcommand{\Eelp}[1]{E_{\rm el}^{(#1)} (u)}
\newcommand{\Esurf}[1]{E_{\rm surf, #1} (\chi)}
\DeclareMathOperator{\Per}{Per(\Omega)}
\newcommand{\diag}{{\rm diag}}

\newcommand{\dx}{\di x}
\newcommand{\eps}{\varepsilon}
\newcommand{\N}{\mathbb{N}}
\newcommand{\Z}{\mathbb{Z}}
\newcommand{\R}{\mathbb{R}}

\newcommand{\Id}{\mathbf{Id}}
\newcommand{\defas}{\coloneqq} 
\newcommand*{\di}{\mathop{}\!\mathrm{d}}

\newcommand{\sym}{{\rm sym}}
\newcommand{\loc}{{\rm loc}}
\newcommand{\skw}{{\rm skew}}
\newcommand{\lc}{\mathrm{lc}}

\usepackage{todonotes}
\usepackage{marginnote}

\DeclarePairedDelimiterX\setof[1]\{\}{#1}
\DeclarePairedDelimiterX\abs[1]\lvert\rvert{#1}
\DeclarePairedDelimiterX\norm[1]\lVert\rVert{#1}
\DeclarePairedDelimiterX\sprod[2]\langle\rangle{#1, #2}

\usepackage{enumitem}

\setlist[enumerate,1]{font=\normalfont}
\setlist[itemize,1]{font=\normalfont}
 
\newlist{thmlist}{enumerate}{1}
\setlist[thmlist]{label=(\roman{thmlisti}),
	ref=(\roman{thmlisti}),font=\normalfont,
	noitemsep}

\usepackage[colorlinks=true, allcolors=blue]{hyperref} 
\usepackage[toc]{appendix}
\usepackage{url}
\usepackage{xurl}

	\title[The energy scaling behaviour of singular perturbation models of staircase type]{The energy scaling behaviour of singular perturbation models of staircase type in linearized elasticity for higher order laminates}
	\author{Lennart Machill}
\address{Institute for Applied Mathematics, University of Bonn, Endenicher Allee 60, 53115 Bonn, Germany}
\email{lmachill@uni-bonn.de}
\author{Angkana Rüland}
\address{Institute for Applied Mathematics and Hausdorff Center for Mathematics, University of Bonn, Endenicher Allee 60, 53115 Bonn, Germany}
\email{rueland@uni-bonn.de}

\begin{document}

	\begin{abstract}
We investigate the scaling behaviour of a singular perturbation model within the geometrically linearized theory of elasticity involving data of higher lamination order. 
We study boundary data which are of staircase type and show rather general lower scaling bounds, both in the setting of prescribed Dirichlet data and for periodic configurations with 
 a mean value constraint. 
 In contrast to the setting without gauge invariances, these lower scaling bounds depend on \emph{two} parameters -- the order of lamination of the boundary data as well as the number of involved (non-)degenerate symmetrized rank-one directions. By discussing upper bounds in specific geometries and for a specific constellation of wells, we give evidence of the sharpness of these lower bound estimates. Hence, it is necessary to keep track of the outlined \emph{two} parameters in deducing scaling laws within the geometrically linearized theory of elasticity.
\end{abstract}


	\maketitle

\setcounter{tocdepth}{2}	
\tableofcontents

	\section{Introduction}
Shape-memory alloys display a striking thermodynamic behaviour:
They undergo a solid-solid, diffusionless phase transition in which crystalline symmetry is reduced in the passage from the high temperature phase (austenite) to the low temperature phase (martensite). This gives rise to various variants of martensite which in turn can be concatenated to complex material patterns \cite{B}. 

 These materials have been very successfully studied by (non-quasi-convex) variational models \cite{B3,B,KM2,KM1,M1}, 
 which display rich energy landscapes and complex microstructures. The objective of the present article is to investigate a class of microstructures arising within the geometrically linearized theory of elasticity. More precisely, we analyze the complexity of the microstructures by means of a singularly perturbed problem and prove scaling laws for the minimal energy in terms of the perturbation parameter. 
 Our analysis shows that the scaling laws are fully determined by two ``order parameters'':
\begin{itemize} 
\item the lamination order of the prescribed Dirichlet (or periodic) data,
\item and the number of (non-)degenerate symmetrized rank-one directions which are used to build up the lamination convex hull.
\end{itemize}
The fact that it is necessary to consider \emph{two} order parameters provides a central novelty of the present article and forms a contrast to earlier works (such as \cite{RT23} or \cite{RT24}) in which it was proved that in the absence of gauges only \emph{a single} order parameter determines the complexity of the microstructure.
 
Contrary to 
 recent papers (such as \cite{RT23a,RT23}), in the present article, we consider models within the geometrically linearized theory of elasticity and, hence, include the gauge group $\text{Skew}(d)$ in our model. This additional ingredient does not only have technical implications for our arguments, but is also of physical significance
 as it reflects infinitesimal frame indifference. In addition, mathematically, it further gives rise to important new phenomena and it is at the heart of the necessity to study a second order parameter in analyzing the microstructure's complexity: As it is by now well-known for the two-well problem \cite{CC14, RRT23, RRTT24}, the passage to a model with the gauge group  $\text{Skew}(d)$ leads to novel scaling laws and, hence, to new microstructures compared to the setting without gauges. In the present article, one of our key objectives is to continue the investigation of these novel microstructures and to prove that their presence is also reflected in \emph{higher order laminates}. In particular, compared to the setting from \cite{RT23}, new structures emerge, which are a genuine consequence of the fact that we include the gauge group in our model.

	\subsection{The problem} 
Let us outline the precise setting of our problem. We analyze singularly perturbed variational problems with energies of the form
	\begin{align}
		\Energy{m+1}{p} & \defas \Eel{m+1}{p} + \epsilon \Esurf{m+1} \notag \\
		& \defas \int_\Omega \vert \sym (\nabla u) - \chi \vert^p \dx + \epsilon \sum_{j=1}^{d} \Vert D \chi_j \Vert_{TV(\Omega)}. \label{def:phaseindic}
	\end{align} 
	Here, $\Omega \defas (0,1)^d \subset \mathbb{R}^d$ denotes the reference configuration. The elastic part of the energy depends on the infinitesimal strain tensor $ \sym (\nabla u) \defas \frac{1}{2}(\nabla u + (\nabla u)^T )$ and a parameter $p \in [2,+\infty)$, where the displacement $u\colon \Omega \rightarrow \mathbb{R}^d$ is an element of a suitable class of admissible functions $\mathcal{A}_F$. The choice $p=2$ corresponds to the usual piecewise quadratic Hooke's law. The matrix-valued function
	$\chi \colon \Omega \to K_m$ is interpreted as a phase-indicator with values in a set  $K_m  \defas \{ A_0, \dots, A_{m} \}\subset \R^{d \times d}_{\sym}$, consisting of $m+1$ wells. 
These wells model the variants of martensite at a fixed temperature below the critical temperature.
 Viewing $\chi$ as the pointwise projection of $\sym(\nabla u)$ onto $K_m$, the surface energy penalizes high oscillations of  $\sym(\nabla u)$, where the functions $\chi_j \colon \Omega \to \{0,1\}$ satisfy the relation
\begin{align}
	\chi = \sum_{j=0}^{m} \chi_j A_j \qquad  \text{ and } \qquad \sum_{j=0}^{m} \chi_j = 1.\label{eq:phaseindicprop}
\end{align}
	We seek to prove a scaling law in the singular perturbation parameter for energies of the type
\begin{align}
\label{eq:minprob}
	\Energyinf{m+1}{p}{ }  \defas \inf\limits_{\chi \in BV(\Omega;K_m)} \inf\limits_{u \in \mathcal{A}_F } \Energy{m+1}{p}
\end{align}
  as $\epsilon \to 0$. We view the scaling behaviour in the singular perturbation parameter as a quantitative measure of the complexity of the underlying microstructure -- the slower the scaling tends to zero as $\epsilon$ approaches zero, the more complex the structures become. Hence, while these scaling laws, in general, do not provide detailed information on the microstructure (e.g., on their exact shape), they do extract important information, e.g., on their overall complexity in the form of (ir-)regularity of the solutions and, in particular, on averages of involved length scales.

 In what follows below, we focus on a specific class $K_m$ of wells, 
 which
can be considered as the geometrically linearized analogue of the wells from \cite{RT23} which 
 have been studied there in the absence of a gauge group. A specific mathematical feature of this set of wells is that its symmetric lamination convex hull, denoted by $K_m^{(\lc)}$, is of staircase type, i.e., it is a union of one-dimensional lines such that $K_m^{(\ell)}\setminus K_m^{(\ell-1)} \neq \emptyset$ for all $\ell \in \{ 1, \dots, m\}$ and $K_m^{(m)} = K_m^{(\lc)}$ (we refer to Definition \ref{def:lamconvhull} for the precise definition of laminates of order $\ell$ and the symmetric lamination convex hull).  We analyze the energy scaling behaviour in the singular perturbation parameter for
 these models, 
 both in the setting 
with a prescribed Dirichlet (i.e., displacement) boundary datum and with a prescribed average condition for the periodic setting, respectively. 
 In particular, we show that there is a close relation between the two energy scalings, 
  in that the periodic one is obtained as a ``shifted'' version of the Dirichlet case. Moreover, compared to the setting without gauges, new scaling bounds will be derived which are a genuine consequence of the presence of the gauge group and are 
 \emph{not} present without it.
	
\subsection{Main results for Dirichlet boundary conditions}

We begin by considering the problem \eqref{eq:minprob} with prescribed Dirichlet data.
In this context, our first goal is to prove that the scaling of the energy is determined by \emph{the order of the lamination of the Dirichlet boundary datum} and the \emph{number of (non-)degenerate symmetrized rank-one directions}. 
Given a matrix $F \in \R^{d \times d}_{\sym}$, we set
\begin{align}
	\mathcal{A}_F^{\dir} \defas \{    u \in W^{1,p}_{\loc} (\R^d; \R^d) : u(x) = Fx + b \text{ in } \R^d \setminus \Omega \text{ for some } b \in \R^d \}.   \label{def:boundarycond}       
\end{align}
For such prescribed boundary conditions, we will consider settings with different complexity and dimension which we introduce next.

\subsubsection{Two wells and a first order datum}
In order to illustrate the effect of the presence of the gauge group $\text{Skew}(d)$,  we begin by recalling the setting in which $K_m$ consists of only two symmetrized rank-one connected wells.
 We  
restrict to two dimensions and consider wells $K_1$ to be defined by the matrices  
	\begin{align}\label{eq:second2dwells}
		K_{1,1}&\defas\{A_0, A_1\} \quad \mbox{ with } \quad  A_0 \defas \diag (0,0) \quad \text{ and } \quad A_1   \defas \diag (1,0) \qquad \qquad  \text{ and } \notag \\
		K_{1,2}&\defas\{A_0, A_1\} \quad \mbox{ with }  \quad A_0 \defas \diag (0,0) \quad \text{ and } \quad A_1   \defas \diag (1,-1) .
		\end{align}
Contrary to the setting without gauges, it was first proved in \cite{CC14} and then later systematized in \cite{RRT23, RRTT24} for more general $\mathcal{A}$-free differential inclusions that the scaling laws for the sets $K_{1,1}$ and $K_{1,2}$ display a \emph{different} scaling behaviour.

\begin{theorem}[Scaling of the two-well problem \cite{CC14}]\label{thm:chanandconti}
\label{thm:CC14} Let $d = 2$, $K_1$ as in \eqref{eq:second2dwells}, and $F \in {\rm int}(K_{1}^{(\lc)})$. 
\begin{itemize}
\item[(i)] If $K_1 \defas K_{1,1}$, there exist constants $C,c>0$ such that for every $\epsilon \in (0,1)$ it holds that 
\begin{align*}
 c  \epsilon^{4/5} \leq \Energyinf{2}{2}{\dir}  \leq C \epsilon^{4/5}. 
\end{align*}
\item[(ii)] If $K_1 \defas K_{1,2}$, there exist constants $C,c>0$ such that for every $\epsilon \in (0,1)$ it holds that 
\begin{align*}
 c \epsilon^{2/3} \leq \Energyinf{2}{2}{\dir} \leq C \epsilon^{2/3}. 
\end{align*}
\end{itemize}
\end{theorem}  
In contrast, the analogous two-well problem for the gradient without gauges would always display a scaling behaviour of the form $\epsilon^{2/3}$ \cite{KM2, KM1}.
As observed in \cite{CC14}, the key difference between the two settings of (i) and (ii) consists of the fact that 
in the model \emph{with} $\text{Skew}(2)$ invariance, generically, there are \emph{two} distinct
 compatibility 
 directions between the wells. This is the case for the set $K_{1,2}$ (with the directions $(1,-1)^T$ and $(1,1)^T$),
 see Section~\ref{sec:compat} below. In the setting of case (i) this degenerates and only \emph{one} symmetrized rank-one direction is present (here the direction $(1,0)^T$).

 As worked out in \cite{RRTT24} the dichotomy of one vs two distinct symmetrized rank-one directions also gives rise to different vanishing orders 
 of the compatibility  conditions for being a symmetrized gradient when restricted to the sphere. Indeed, the compatibility conditions for being a gradient (as in the setting without a gauge) is encoded in an annihilator which is a first order operator (the fact that the curl of a gradient field vanishes by the Poincar\'e lemma). The compatibility conditions for being a symmetrized gradient, however, is a second order operator (the Saint-Venant conditions). Considering the symbol of these operators in the direction prescribed by the wells restricted to the sphere in frequency space, the Saint-Venant compatibility conditions may lead to first (in the case of $K_{1,1}$) or second order (in the case of $K_{1,2}$) vanishing conditions towards the non-elliptic set, see also Lemma~\ref{lem:elasticenergybound} below. In contrast, the first order operator (curl) always has a symbol with vanishing order of degree one towards its characteristic set, see e.g.\ \cite[Example 2.1 and Lemma 3.4]{RRTT24}. Mathematically, 
 it is this structure which relates to 
 the different scaling laws in Theorem \ref{thm:CC14}.

This behaviour is by now well-understood for the \emph{two-well} problem which involves laminates up to \emph{first} order. A natural question -- both mathematically and physically -- deals with the effect of this different vanishing behaviour on \emph{higher order laminates} and their scaling laws. Here different combinations of vanishing orders can arise.
It is one of our main purposes in the present article to investigate this situation in the case of a model setting and to outline that the different vanishing orders lead to \emph{novel} scaling laws also for higher order laminates. As one of our main results, for our family of wells, we show that the complexity of the microstructure is determined by the lamination order of the boundary data \emph{and} the number of degenerate vs non-degenerate symmetrized rank-one directions involved in this lamination process.

\subsubsection{Three wells and a second order datum}

\begin{figure}[t]
  \centering
  \begin{minipage}[b]{0.45\textwidth}
    \centering
\begin{tikzpicture}
  \draw (0, 0)   to[out=0, in=180] (4, 0);
  \draw (2, 0) to[out=90, in=-90] (2, 4); 
 
 \node at (-0.3,0) {$A_0$};
\node at (4.3,0) {$A_1$};
\node at (2,4.3) {$A_2$};

\node at (0,0) [circle,fill,inner sep=0.5pt]{};
\node at (4,0) [circle,fill,inner sep=0.5pt]{};
\node at (2,4) [circle,fill,inner sep=0.5pt]{};

\node at (1,0) [circle,fill,inner sep=0.215cm, color=white, opacity= 0.7]{};
   \node (n) at (1,0) {$1$}; 
  \draw[thick] (n) circle[radius=0.3cm];

\node at (2,2) [circle,fill,inner sep=0.215cm, color=white, opacity= 0.7]{};
   \node (n) at (2,2) {$2$}; 
  \draw[thick] (n) circle[radius=0.3cm];

\end{tikzpicture}
  \end{minipage}
  \hfill
  \begin{minipage}[b]{0.45\textwidth}
    \centering
\begin{tikzpicture}
  \draw (0, 0)   to[out=0, in=180] (4, 0);
  \draw (2, 0) to[out=90, in=-90] (2, 4); 
 
 \node at (-0.3,0) {$A_0$};
\node at (4.3,0) {$A_1$};
\node at (2,4.3) {$A_2$};

\node at (0,0) [circle,fill,inner sep=0.5pt]{};
\node at (4,0) [circle,fill,inner sep=0.5pt]{};
\node at (2,4) [circle,fill,inner sep=0.5pt]{};

\node at (1,0) [circle,fill,inner sep=0.215cm, color=white, opacity= 0.7]{};
   \node (n) at (1,0) {$2$}; 
  \draw[thick] (n) circle[radius=0.3cm];

\node at (2,2) [circle,fill,inner sep=0.215cm, color=white, opacity= 0.7]{};
   \node (n) at (2,2) {$1$}; 
  \draw[thick] (n) circle[radius=0.3cm];

\end{tikzpicture}
  \end{minipage}
  \caption{A schematic two-dimensional visualization of the symmetrized lamination convex hull in the cases \eqref{eq:choiceofmatrices} (left) and \eqref{eq:secondcase} (right). The numbers enclosed by the circle represent the number of compatibility directions between the wells. }
    \label{fig:3dwells}
\end{figure}
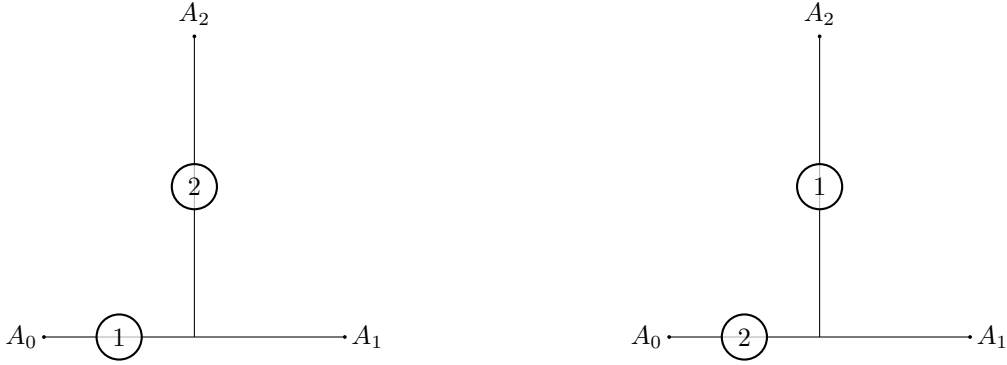

We begin our discussion of higher order laminates by considering  second order laminates. We will here focus on a prototypical choice of the wells:
Let $d = 3$ and assume that $K_{2}$ is given by  
	\begin{align}
		K_{2,12}&\defas \{A_0, A_1, A_2\}   \mbox{ with } A_0 \defas \diag (0,0,0),  \,  A_1   \defas \diag (1,0,0),  \, A_2  \defas \diag \Big(\frac{1}{2},1,-1 \Big) \text{ or}   \tag{W.1}   \label{eq:choiceofmatrices}    \\
		K_{2,21}&\defas\{A_0, A_1, A_2\}   \mbox{ with } A_0 \defas \diag (0,0,0),  \,  A_1   \defas \diag (1,-1,0), \, 	 A_2 \defas \diag \Big(\frac{1}{2},-\frac{1}{2},1 \Big).   \tag{W.2}    \label{eq:secondcase}
	\end{align}
Both sets of wells are chosen such that their symmetrized lamination convex hulls are one-dimensional and that they contain one lamination order determined by a single  symmetrized rank-one connection and one 
lamination order given by two symmetrized rank-one connections. Their ordering however does not coincide, see Figure~\ref{fig:3dwells}. More precisely, the lamination convex hulls of the sets $K_{2,12}$ and $K_{2,21}$ are given by
 $K_{2,21}^{(\lc)} = K_{2,21}^{(2)}$ and $K_{2,12}^{(\lc)} = K_{2,12}^{(2)} $, where 
\begin{align}
\label{eq:3D_wells}
 \hspace*{-0.2cm} K_{2,12}^{(1)} \setminus K_{2,12} &= \left\{ \diag(\alpha, 0,0):  \alpha \in (0,1)\right\} , && \hspace*{-0.2cm} K_{2,12}^{(2)} \setminus K_{2,12}^{(1)} = \left\{ \diag \left(\frac{1}{2}, \alpha, - \alpha\right):  \alpha \in (0,1) \right\}, \notag \\ 
  \hspace*{-0.2cm} K_{2,21}^{(1)} \setminus K_{2,21} 
&= \left\{ \diag(\alpha, -\alpha,0):  \alpha \in (0,1)\right\} , && \hspace*{-0.2cm} K_{2,21}^{(2)}  \setminus K_{2,21}^{(1)} = \left\{ \diag \left(\frac{1}{2}, \frac{1}{2},  \alpha\right):  \alpha \in (0,1) \right\}.
\end{align}
We refer to 
 Definition~\ref{def:lamconvhull} 
for the notion of the symmetrized lamination convex hull.
 Hence, $K_{2,12}^{(1)} \setminus K_{2,12} $ is obtained from a degenerate symmetrized rank-one connection with a single normal while  $K_{2,12}^{(2)} \setminus K_{2,12}^{(1)}$ is associated with 
 two distinct normals. In the case $K_{2,21}$ this behaviour is reversed. We keep track of this ordering in the subindex; the index one indicates the presence of a degenerate symmetrized rank-one connection (with only one normal), the index two that there is a non-degenerate symmetrized rank-one direction (with two distinct normals).

We show that independently of this ordering, the energy scaling of microstructures with boundary data from the second order laminate hull and for $p=2$ feature a lower scaling bound of the size of $\sim \eps^{4/7}$ in both settings while the scaling laws for data from the first order hull depend on the presence of a degenerate/non-degenerate symmetrized rank-one direction.

\begin{theorem}[Lower bounds for the three-well problem]\label{mainth:3d}
Let $K_{2,12}$ and $K_{2,21}$ be as in \eqref{eq:choiceofmatrices} and \eqref{eq:secondcase} and $p \in [2,+\infty)$. 
\begin{itemize}
\item[(i)] \emph{First order laminates.} Let $F \in K_{2,12}^{(1)} \setminus K_{2,12}$. Then, there exists a constant $c>0$ such that for every $\epsilon \in (0,1)$ it holds that
\begin{align*}
	c\epsilon^{2p/(2p+1)} \leq \Energyinf{3}{p}{\dir}.
\end{align*}
 Assume that 
$F \in K_{2,21}^{(1)} \setminus K_{2,21}$. Then, there exists a constant $c>0$ such that for every $\epsilon \in (0,1)$ it holds that
\begin{align*}
	c\epsilon^{2p/(2p+2)} \leq \Energyinf{3}{p}{\dir}.
\end{align*}
\item[(ii)] \emph{Second order laminates.} Let $K_2$ be either $K_{2,12}$ or $K_{2,21}$,
 and assume that $F \in K^{(2)} \setminus K^{(1)}$. Then, there exists a constant $c>0$ such that for every $\epsilon \in (0,1)$ it holds that
\begin{align*}
	c\epsilon^{2p/(2p+3)}\leq \Energyinf{3}{p}{\dir}.
\end{align*}
\end{itemize}
\end{theorem}

We highlight that the lower scaling bound for second order laminates is new and is independent of the ordering of the degenerate/non-degenerate symmetrized rank-one connections. 
In contrast, in the analogous setting in \cite{RT23} for models without gauges for $p=2$, one obtains a scaling behaviour of the order $\epsilon^{\frac{1}{2}}$ (an analogous bound is proved in Theorem \ref{mainth:md} below when considering two non-degenerate rank-one directions). We highlight that compared to the setting without the presence of gauges, our new scaling bound (which for $p=2$ takes the form $\epsilon^{\frac{4}{7}}$) requires less energy. Similarly, as in the case of Theorem \ref{thm:CC14}(i) this is due to the fact that in our new setting it is possible to leverage the genuinely vectorial structure of the problem.  

 Moreover, we emphasize that, to the best of our knowledge, Fourier localization methods have been successfully employed for quadratic elastic energies, relying on the isometric $L^2$-isomorphism of real and Fourier space. By the simple observation that $L^p$-based elastic energies control the $L^2$-based counterpart for $p \geq 2$, our results also show that the Fourier analytical approach is robust for $L^p$-based elastic energies in this case. For the sake of completeness, we mention that real-space localization methods have already been successfully employed for a general $p \in [1,+\infty)$ in the context of two-well problems, see e.g.\ \cite{CDZ17} or \cite[Section~3.2]{RT23}.   

 Since three-dimensional branching constructions can become complex if 
boundary conditions are prescribed on the whole boundary, 
we restrict ourselves to the case \eqref{eq:choiceofmatrices} when justifying the optimality of the lower bound.
 In this setting, we consider both the case of \emph{first} and \emph{second} order lamination data as boundary conditions.
For the first order lamination data in two-dimensional settings as outlined in Theorem \ref{thm:CC14}(i) above, the known upper bound microstructures that generate the optimal scaling behaviour are based on branching constructions. We highlight that the upper bound construction in Theorem \ref{thm:CC14}(i) heavily exploits the vectorial character of the problem. As a consequence, the transfer of these constructions to three dimensions, does not simply follow by an affine extension to the additional dimension if boundary conditions are prescribed on the \emph{whole} boundary.
 In this case, one would expect that such constructions also need to branch out in the additional direction.
 While this has been successfully analyzed in the absence of gauges by a suitable rotation of the two-dimensional construction \cite{RT23}, the situation becomes more complex in the presence of gauges if only one degenerate compatibility direction exists.
 Here, both components of the competitor for the two-dimensional upper bound constructions need to be taken into account and complicate the latter approach in three dimensions,
 see also Section~\ref{sec:challenges} below. We provide an upper bound construction in a cylindrical domain.

\begin{proposition}[Upper bound for a cylinder in three dimensions]\label{prop:Uppercylinder}
	Consider a three-dimensional cylinder $\Omega = (0,1) \times B_1(0)\subset \R^3$, the wells
	$K_1 \defas \{ \diag(0,0,0), \ \diag(1,0,0) \}$, $p \in [1,+\infty)$, and $F = \diag(\frac{1}{2},0,0)$. 
	 Then, there exists a constant $C>0$ such that for every $\epsilon \in (0,1)$ it holds that 
\begin{align*}
  \Energyinf{2}{ p }{\dir}  \leq C \epsilon^{ 2p/(2p+1) }.
\end{align*}
\end{proposition}

 For the case of \emph{second order boundary conditions}, we also turn to a specific geometry which is particularly suited to our problem. 
 More precisely, we provide an upper bound for a rotated cube $(0,1) \times \Phi^{-1}((0,1)^2)$,
 which is designed such that the associated compatibility directions match the basis vectors of the cube.
For our choice of wells, the associated rotation is given as $\Phi\colon \R^2 \to \R^2$,  $\Phi(z) \defas \frac{\sqrt{2}}{2} (z_1 + z_2, z_2 - z_1)^T$.
\begin{theorem}[Upper bound for a second order laminate in the case \eqref{eq:choiceofmatrices}]\label{thm:mainresultupper}
Let $p \in [1,3]$, $\Omega \defas (0,1) \times \Phi^{-1}((0,1)^2)$, and
\begin{align*}
	 F = \diag\bigg(\frac{1}{2}, \frac{1}{2}, -\frac{1}{2}\bigg) \in K_2^{(2)} \setminus K_2^{(1)},
\end{align*}
where $K_2 = K_{2,12}$.
Then, there exists a constant $C>0$ such that for every $\epsilon \in (0,1)$ it holds that
\begin{equation*}
\Energyinf{3}{p}{\dir} \leq C\epsilon^{2p/(2p+3)} .
\end{equation*}
\end{theorem}
The upper bound is constructed via two nested (three-dimensional) branching constructions on two different length scales.
More precisely, the innermost branching pattern is based on a suitable version of Proposition~\ref{prop:Uppercylinder} on cubes instead of cylinders.
Here, the geometry of the domain enters which is why we resort to a specific cut-off argument:
As the branching pattern associated to a single compatibility direction occurs on the smaller length scale (see the discussion below Proposition~\ref{lem:2well3d}), we recover matching scaling bounds for $p\in [2,3]$ while our cut-off is energetically too expensive for $p>3$. 
Nevertheless, heuristic calculations suggest that the lower bound in the case \eqref{eq:secondcase} is sharp as well, which is discussed in Subsection~\ref{sec:heuristics}. 
The combination of the two different vanishing orders hence yields a novel scaling law.  
 Both Proposition~\ref{prop:Uppercylinder} and Theorem~\ref{thm:mainresultupper} are proved in Section~\ref{sec:upper}, following from Proposition~\ref{prop:optimalfirst-order} and Proposition~\ref{prop:upperboundcons}, respectively. 

One may wonder whether an even lower scaling law arises for the setting in which both of the involved symmetrized rank-one connections are degenerate. As is shown below, this case however is rather different from the previous ones as the lamination convex hull here becomes a higher dimensional set and no longer consists of single line segments (see Proposition \ref{prop:staircaseconstraint}). Due to this structural change and the expectation that this leads to a substantial amount of additional flexibility, we do not discuss this case at this point but postpone it to future work.

\subsubsection{Four wells and a third order datum}
In a similar spirit as in the case of second order laminates, one can analyze third order laminate data, provided that the lamination convex hull does not contain a structure where two consecutive laminates with only one compatibility direction exist. As an example of such a setting, we deduce the following lower scaling bound. 
	\begin{theorem}[Lower bounds for a four-well problem]\label{mainth:4d}
	Let $d = 4$, $p \in [2,+\infty)$, and assume that $K_{3}$ is defined via the matrices
	\begin{align}
		A_0 &\defas \diag (0,0,0,0),   && A_1   \defas \diag (1,0,0,0), \notag \\
		 A_2 &\defas \diag \left(\frac{1}{2},1,-1,0 \right),  && A_3  \defas \diag \left(\frac{1}{2},\frac{1}{2},-\frac{1}{2},1 \right).  \label{eq:Kin4d}  
	\end{align}
Let $F \in K_{3}^{(3)} \setminus K_{3}^{(2)}$. Then, there exists a constant $c>0$ such that for every $\epsilon \in (0,1)$ it holds that
\begin{align*}
	c\epsilon^{2p/(2p+4)} \leq \Energyinf{4}{p}{\dir}.
\end{align*}
\end{theorem}  
 We note that the wells are such that $A_0,A_1$ are symmetrized rank-one connected via exactly one degenerate symmetrized rank-one direction and there are no direct further rank-one connections between the other wells. However, $A_2$ is symmetrized rank-one connected with one element of the line $[A_0, A_1]$ via a non-degenerate symmetrized rank-one connection, and $A_3$ is symmetrized rank-one connected with the line $[\frac{1}{2} (A_0 + A_1), A_2]$ via a single degenerate symmetrized rank-one connection.
As above, the scaling bound in Theorem \ref{mainth:4d} is \emph{new} and part of a larger family of new scaling laws outlined in the next section. For $p=2$ its counterpart for the gradient from \cite{RT23} displays a scaling bound of the order $\epsilon^{\frac{2}{5}}$. Again, the presence of a degenerate symmetrized rank-one connection hence reduces the energetic cost, which for $p=2$ is of the order $\epsilon^{\frac{1}{2}}$ in the above setting. In our discussion below, we treat the case of four wells separately as it provides important intuition for the general $d$-dimensional setting.

\subsubsection{$m+1$ wells and higher order data} 
With the previous examples in mind, we now turn to the general construction of our set of wells and the discussion of the associated scaling laws.
More precisely, for $m \in \N$ we seek to construct a set $K_m$ such that $K^{(\lc)}_m= K^{(m)}_m$ and $ \emptyset \neq K^{(j)}_m \setminus K^{(j-1)}_m$ are one-dimensional lines for all $j \in \{1,\dots,m \}$.
We denote by 
\begin{align*}
 k_m \defas \# \big\{ j \in \{1,\dots,m \} :  & \text{ the only compatibility direction in }  \\
 & K^{(j)}_m \setminus K^{(j-1)}_m \text{ is a degenerate connection} \big\} 
\end{align*}
the number of distinct one-dimensional lines in $K^{(\lc)}_m$, which are formed by degenerate symmetrized rank-one directions.
 Labeling the order of the laminate in $K^{(\lc)}_m$ by $\{1, \dots, m \}$   ($j\in\{1, \dots, m \}$ corresponds to the line $K^{(j)}_m \setminus K^{(j-1)}_m$), 
we associate a function $f \colon \{1, \dots, m \} \to \{1,2\}$ which indicates the number of compatibility directions of the corresponding laminate, i.e., it encodes the (non-)degeneracy of these compatibility directions. 
 With this notation in hand we define for each $r\in \{  0  ,\dots,m\}$ the number of degenerate symmetrized rank-one directions up to the $r$th order by 
\begin{align}
\label{eq:kr}
 k_{r} \defas \# \{ z \in \{1,\dots, r\}: f(z) = 1 \},
\end{align}
where we note that $k_0 \defas 0$ and $r-k_r = \# \{ z \in \{1,\dots, r\}: f(z) = 2 \} $.

Due to the structural result in Proposition~\ref{prop:staircaseconstraint} below, 
 we do not consider settings with two consecutive degenerate symmetrized rank-one connections, and define the class of admissible laminates by
 \begin{align}\label{def:admissiblelaminates}
	\mathcal{S}_m \defas \big\{ f \colon  \{1, \dots, m  \}  \to \{1,2\} : & \text{ for no index } i \in \{1, \dots, m -1 \}  \notag \\ 
	& \qquad    \ \ \text{ it holds that } f(i) = f(i+1) = 1 \big\}.
 \end{align}
  Indeed, by Proposition~\ref{prop:staircaseconstraint}, these would give rise to lamination convex hulls, which no longer consist of one-dimensional segments but include higher dimensional parts. We thus expect these sets to be substantially more flexible and postpone their analysis to future work.

\begin{remark}[Cardinality of $\mathcal{S}_m$]
	There exist $2^m$ choices for the functions $f \colon \{1, \dots, m  \} \to \{1,2\}$.
	As we exclude functions $f$ with $f(i) = f(i+1) = 1$, the sequence $(\# \mathcal{S}_m)_m$ is a Fibonacci sequence.
	Indeed, for $m = 1$ we have $\# \mathcal{S}_1 = 2$. Given an admissible element, we can always attach a non-degenerate connection to construct a valid choice for $m =2$.
	If $f(1) = 2$, a further degenerate connection is admissible which can be added for one element in $\mathcal{S}_1$, implying that $\# \mathcal{S}_2 = 3$.
	Admissible functions for $m = 3$ can then again be formed by either adding non-degenerate connections for every element in $\mathcal{S}_2 $ or by attaching degenerate connections which are only admissible if the previous connection is non-degenerate.
	By construction the number of possible ways for the latter formation coincides with $\# \mathcal{S}_1$. By an inductive iteration, we obtain $\# \mathcal{S}_m = \# \mathcal{S}_{m-1} + \# \mathcal{S}_{m-2}$.
\end{remark} 

 \medskip

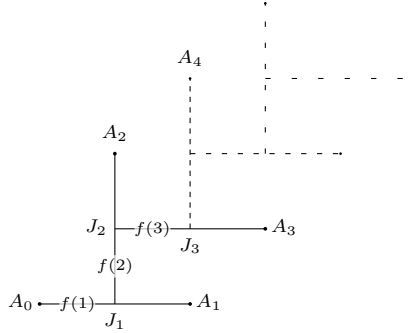
\begin{figure}
	\centering
	\begin{tikzpicture}
  \draw (0, 0)   to[out=0, in=180] (2, 0);
  \draw (1, 0) to[out=90, in=-90] (1, 2);
  \draw (1, 1) to[out=0, in=180] (3, 1) ;
  \draw  [dash pattern={on 2pt off 2pt on 2pt off 2pt}] (2, 1) to[out=90, in=-90] (2, 3);
	\draw [dash pattern={on 2pt off 4pt on 2pt off 4pt}](2, 2) to[out=0, in=180] (4, 2) ;
  \draw [dash pattern={on 2pt off 6pt on 2pt off 6pt}](3, 2) to[out=90, in=-90] (3, 4);
  \draw [dash pattern={on 2pt off 8pt on 2pt off 8pt}] (3, 3) to[out=0, in=180] (5, 3) ;
 \node[below, font=\scriptsize] at (-0.25,0.25) {$A_0$};
\node[below,font=\scriptsize] at (2.25,0.25) {$A_1$};
\node[below,font=\scriptsize] at (1,0) {$J_1$};
\node[below,font=\scriptsize] at (1,2.5) {$A_2$};
\node[below,font=\scriptsize] at (0.75,1.25) {$J_2$};
\node[below,font=\scriptsize] at (2,1) {$J_3$};
\node[below,font=\scriptsize] at (3.25,1.25) {$A_3$};
\node[below,font=\scriptsize] at (2,3.5) {$A_4$};

\node at (0,0) [circle,fill,inner sep=0.5pt]{};
\node at (2,0) [circle,fill,inner sep=0.5pt]{};
\node at (1,2) [circle,fill,inner sep=0.5pt]{};
\node at (3,1) [circle,fill,inner sep=0.5pt]{};

\node at (2,3) [circle,fill,inner sep=0.2pt]{};
\node at (3,4) [circle,fill,inner sep=0.2pt]{};
\node at (4,2) [circle,fill,inner sep=0.2pt]{};
\node at (5,3) [circle,fill,inner sep=0.2pt]{};

\node at (0.5,0) [circle,fill,inner sep=4.8pt, color=white, opacity= 0.7]{};
 \node[font=\tiny] at (0.5,0) {$f(1)$};
 \node at (1,0.5) [circle,fill,inner sep=3pt, color=white, opacity= 0.7]{};
\node[font=\tiny] at (1,0.5) {$f(2)$};
\node at (1.5,1) [circle,fill,inner sep=4.8pt, color=white, opacity= 0.7]{};
\node[font=\tiny] at (1.5,1) {$f(3)$};

\end{tikzpicture}
\caption{ Schematic visualization of the staircase in \eqref{eq:constructiondim}. 
The nodes represent the set $K_m$ while the one-dimensional lines 
visualize the symmetrized lamination convex hull $K_m^{(\lc)}$. The function $f$ represents the number of compatibility directions between matrices on these lines.
Notice that in two dimensions,
 $K_m^{(\lc)}$ cannot be one-dimensional, see Proposition~\ref{prop:topoprop} below,
 i.e., this figure should be interpreted schematically. }
\label{fig:staircase}
\end{figure}

\noindent\textbf{Explicit construction of the wells.} 
 For a given $m \in \N$, a given function $f:\{1,\dots,m\} \rightarrow \{1,2\}$ and parameter $k_m$ as above, the set $K_m := \{A_0, \dots, A_m\} \subset \R^{d\times d}_{\rm diag}$ is constructed in dimension $d = 2(m- k_m) + k_m $ for $m \geq 2$. 
We begin by setting
  $A_0 = J_0 = \diag(0,\dots, 0) \in \R^{d \times d}_{\rm diag}$. Then, we inductively define $A_i := J_{i-1} + M_{i}$, where the matrices are defined as
\begin{align}\label{eq:constructiondim}
	J_i = \frac{1}{2} \sum_{k=1}^i M_k \quad \text{and} \quad (M_i)_{l \tilde l} := \begin{cases}
		1 & \text{if }    l = 1 + \sum_{j=1}^{i-1} f(j)  \ \text{and } l = \tilde l , \\
		-\delta_{f(i)2} & \text{if }  l = 2 + \sum_{j=1}^{i-1} f(j)  \ \text{and } l = \tilde l,  \\
		0 & \text{else}
	\end{cases}
\end{align}
for any $i = 1, \dots, m$, where $\delta_{f(i)2}$ indicates the Kronecker-delta function with mass concentrated in  $(f(i),2) \in \R^2$.  
 We note that these wells satisfy the following properties:
\begin{itemize}
\item The wells $A_0$ and $A_1$ are symmetrized rank-one connected. This rank-one connection is degenerate if $f(1) = 1$ and non-degenerate if $f(1)=2$.
\item There are no further symmetrized rank-one directions among the other wells in $K_m$.
\item For $i > 1$, each well $A_i$ is symmetrized rank-one connected with the auxiliary well $J_{i-1}$. This symmetrized rank-one connection builds up the lamination convex hull which consists of one-dimensional segments.
\end{itemize} 
A schematic visualization of this construction is provided in Figure~\ref{fig:staircase}.
In this setting, we deduce the following lower scaling bound which includes the ones from the previous sections as special cases: 
	\begin{theorem}[Lower bounds for the $m+1$-well problem]\label{mainth:md}
		 Let $p \in [2,+\infty)$,
		 $m \in \N$, $\ell \in \{1,\dots,m\}$ and $k_{\ell}$ as in \eqref{eq:kr} above. 
  Assume that $K_m$ is constructed as in \eqref{eq:kr}--\eqref{eq:constructiondim} for $f \in \mathcal{S}_m$ and
suppose that $F \in K^{(\ell)}_m \setminus K^{(\ell-1)}_m$.
 Then, there exists a constant $c>0$ such that for every $\epsilon \in (0,1)$ it holds that
\begin{align*}
	c\epsilon^{2p / (2p+   2(\ell-k_{\ell}) + k_{\ell} )} \leq  \Energyinf{m+1}{p}{\dir} .
\end{align*}
 Moreover, as $\dist(F , K^{(\ell-1)}_m) \rightarrow 0$ the constant $c$ tends to zero. 
\end{theorem}  

\begin{remark}

It would be possible to make the constant $c>0$ explicit in its dependence on $\dist(F, K^{(\ell-1)}_m)$ which would yield a power type behaviour. As we expect that these estimates are not optimal in the exponent, we do not pursue the discussion of this dependence further in this article but postpone it to a future investigation.

\end{remark}

Let us comment on the result of Theorem \ref{mainth:md}. Firstly, we note that if $k_{\ell} = 0$, i.e., if only non-degenerate symmetrized rank-one connections are 
 used, our 
scaling behaviour turns into 
\begin{align*}
c\epsilon^{ p / (p+  \ell ) } \leq  \Energyinf{m+1}{p}{\dir} .
\end{align*}
This bound recovers the estimate from \cite{RT23} in which such a bound for $p = 2$ is deduced for a gradient differential inclusion with laminates of the order $\ell$ \emph{without} gauge invariances.
 This is in accordance with the expectation (see \cite{RRTT24} for the two-well problem) that if only non-degenerate symmetrized rank-one connections are present, then the scaling laws for the gradient and symmetrized gradients coincide.
Secondly, if $k_{\ell} \neq 0$, we obtain \emph{new} lower scaling bounds. If $\ell = 1$ and $k_{\ell} = 1$ this coincides with the estimate from \cite{CC15} (see Theorem \ref{thm:CC14}(i)) while for $\ell \neq 1$ the result provides new estimates.
Theorem \ref{mainth:md}, thus, shows that for Dirichlet data the lower scaling behaviour is determined by the two parameters $\ell$ and $ k_{\ell}$ which determine the ``complexity'' of the microstructure and the associated scaling law. Finally, we note that while we do not carry out the upper bound constructions in the present article, based on the constructions from \cite{RT23}, we believe that the here obtained upper bounds are sharp. 
 More specifically, by assuming that one can construct upper bounds via nested branching constructions on $\ell$ different length scales,
which is made rigorous in the model without gauges \cite[Theorem~1.4]{RT23},
we present a heuristic argument for the sharpness of the lower scaling bounds, see Lemma~\ref{lem:heuristic} in Section \ref{sec:heuristics}. We expect that a rigorous argument for this requires overcoming substantial technical difficulties and novel constructions (e.g.~iterative combinations of curved and linear branching structures) and hence postpone a rigorous analysis to future projects.

\subsection{Main results for periodic boundary conditions}\label{sec:periodicboundary}
As a second scenario, we consider the case in which periodic boundary conditions and a mean value constraint are imposed instead of Dirichlet boundary conditions.
More precisely, given a matrix $F \in \R^{d \times d}_{\sym}$, we set  
\begin{align}
	\mathcal{A}_F^{\per} \defas \{    u \in W^{1,p}_{\loc} (\R^d; \R^d) : \langle \nabla u \rangle = F , \ \nabla   u \text{ is } [0,1]^d\text{-periodic} \},   \label{def:periodicboundarycond}     
\end{align}
 where $\langle \nabla u \rangle \defas \int_{[0,1]^d} \nabla u \di x $. 
In this context, we will provide lower bounds for 
\eqref{eq:minprob} and, for $m \in \N$, $f\colon\{1,\dots,m\} \rightarrow \{1,2\}$ and $k_m$ given, we consider the sets $K_m$ from the construction in \eqref{eq:constructiondim}.

As in the setting with the Dirichlet data, we will also in this context prove that the lower bounds are completely determined by the two parameters $\ell$ and $k_{\ell}$, i.e., the lamination order of the data and the number of the degenerate symmetrized rank-one connections associated with these data. In comparison to the Dirichlet setting, one however expects the periodic setting to be energetically stricly less costly 
 (Dirichlet boundary conditions are more restrictive and can be viewed as a smaller subset of the periodic setting in the minimization problem \eqref{eq:minprob}, i.e., $\mathcal{A}_F^{\dir} \subsetsim \mathcal{A}_F^{\per} $ via the divergence theorem). 
 In fact this is reflected in our scaling bounds 
 in 
  that they are ``shifted'' by one order compared to the Dirichlet case.

	\begin{theorem}[Lower bounds for the $m+1$-well problem in the periodic case]\label{mainth:mdperiodic}
	 Let $p \in [2,+\infty)$,
	 $m \in \N$, 
		and $\ell \in \{1,\dots,m\}$.
  Assume that $K_m$ is constructed as in \eqref{eq:kr}--\eqref{eq:constructiondim} for $f \in \mathcal{S}_m$ and
suppose that $F \in K^{(\ell)}_m \setminus K^{(\ell-1)}_m$. 
Then, there exists a constant $c>0$ such that for every $\epsilon \in (0,1)$ it holds that
\begin{align*}
	  c\epsilon^{2p / (2p+  2((\ell-1)- k_{\ell-1}) + k_{\ell-1} ) } \leq  \Energyinf{m+1}{p}{\per} .
\end{align*}
Moreover, $c \rightarrow 0$ as $\dist(F, K_m^{(\ell-1)}) \rightarrow 0$.
\end{theorem}
Compared to Theorem \ref{mainth:md}, the exponents in Theorem \ref{mainth:mdperiodic} differ due to the larger class of competitors in the periodic setting. For mean value data in the first order lamination hull, the scaling is always of the order $\epsilon$. The corresponding upper bound construction simply consists of a simple laminate. The more complex structures only enter in higher order laminates. Here the outermost laminate should be given by a simple laminate which does not need to satisfy Dirichlet data. 
Only in the next orders of lamination 
the regularity of the deformation requires the finer branching constructions to have matching Dirichlet boundary conditions with the constructed first order laminate. 
More specifically, one expects that $\ell-1$ suitably nested branching constructions on $\ell-1$ different length scales can generate optimal upper bounds. 
This accounts for the shift in the scaling law compared to the Dirichlet setting. 
 As in the Dirichlet setting, the optimality of the lower bounds is discussed in Subsection~\ref{sec:heuristics}. Here, we restrict ourselves to heuristic arguments
once again due to the technicalities involved in these constructions. 

\subsection{Some remarks on the challenges associated with the presence of gauges}\label{sec:challenges}

As outlined above, the presence of gauges distinguishes the present work from earlier results on scaling laws for higher order laminates such as, for instance, \cite{RTTZ25,RT23a,RT23}. Contrary to the setting without gauges, it gives rise to the presence of \emph{two} relevant parameters in our scaling laws -- the \emph{order of lamination} of the Dirichlet boundary data and the \emph{number of (non-)degenerate symmetrized rank-one connections}. Also from a more technical point of view, the presence of gauges and -- related to this -- the second order differential equations for the compatibility conditions for symmetrized gradients lead to novel challenges and phenomena. 
Let us only comment on two of these aspects which play a fundamental role in our arguments below:
\begin{itemize}
\item \emph{The coercivity conditions.} One central ingredient in our Fourier-based argument consists in the quantification of the compatibility conditions for symmetrized gradients. The compatibility conditions are reflected in highly anisotropic, nonlocal elastic energies which are coercive outside of certain regions in Fourier space. In the setting \emph{without} gauges, in similar situations involving wells of staircase type (see e.g.\ \cite{RT23}), the non-coercive regions corresponded to cones along the rank-one directions between the wells. In the context of our present setting \emph{with} gauges, the situation changes (see Lemma \ref{lem:elasticenergybound}): While the elastic energy still provides strong localization around the cones associated with the non-degenerate symmetrized rank-one directions, we only obtain much weaker coercivity around the degenerate symmetrized rank-one directions. More precisely, using purely linear estimates stemming from the elastic energy (without the additional nonlinear structure of the wells), we cannot separate these degenerate symmetrized rank-one directions. As a consequence, we only obtain control outside of a cone around a vector space spanned by \emph{all} degenerate directions.
 This leads to major additional effort in the overall bootstrap argument in deducing the lower bound estimates. In addition, and as indicated already in the two-well case in \cite{RRT23,RRTT24}, the two types of symmetrized rank-one directions also differ in their rate of vanishing towards the non-elliptic set. These different vanishing orders imply that in contrast to the setting without gauges, in our proofs of the lower bound estimates we need to work with two types of cones with different opening angles. 

\item \emph{The upper bound constructions and their building blocks.} Also in the upper bound estimates the presence of gauges manifests itself in novel challenges. This is particularly the case in settings involving degenerate symmetrized rank-one conditions. Due to the work of Chan-Conti \cite{CC14} it is known that in these settings upper bound constructions which saturate the lower bound estimates exploit the full vectorial nature of the problem.
 In the setting of higher order laminate boundary conditions and higher dimensions, this thus requires substantial amount of additional care in constructing \emph{and} concatenating the structures on the different lamination levels. In particular,
 one cannot simply resort to three-dimensional constructions from settings without gauges (e.g.\ from \cite{RT23}) as 
 curvilinear interfaces and more careful cut-off arguments have to be used. This is already the case in our upper bound construction for the specific setting from \eqref{eq:choiceofmatrices}, see the discussion of Lemma~\ref{lem:2well3d}.  We expect that it will play an even stronger role in the general upper bound construction for which we only give heuristic arguments (see Section \ref{sec:heuristics}) and which we postpone to future work.
\end{itemize}

 Hence, although we follow a Fourier-based strategy for the derivation of lower bounds which had been introduced in earlier works such as \cite{RTTZ25,RT23a,RT22}, substantial novelties arise due to the presence of a gauge invariance. Besides general technical difficulties that are already generated in real space, challenges in Fourier space are particularly reflected in combining the natural coercivity constraints from the energy with the nonlinear structure of the wells.
 Additionally, branching constructions for higher order data in three dimensions require a careful treatment due to the vectorial character of the problem. 
 In conclusion, the challenges induce both new phenomena and technical difficulties which have to be overcome.

\subsection{Relation to the literature}
Let us connect our article to the literature on scaling laws for models inspired by shape-memory alloys. As in \cite{RT23, RT23a}, it is our aim to rigorously determine the complexity of microstructure involving \emph{higher order laminates} by ``simple'' order parameters. In contrast to the setting without gauges from \cite{RT23}, the results from above strongly suggest that for the geometrically linearized setting a single order parameter does not encode the full behaviour. Additionally to the order of lamination which was identified as the only necessary order parameter in the setting of \cite{RT23}, in the framework \emph{with gauges} we further identify the number of degenerate symmetrized rank-one connections as a second order parameter. This is directly related to earlier work due to Chan-Conti 
\cite{CC14}, who showed that in the two-well setting the scaling law in the degenerate and non-degenerate cases differs. We also refer to \cite{RRT23, RRTT24} for a systematic study of this for quantifications of a family of $\mathcal{A}$-free differential inclusions, and to \cite{PF25} for scaling laws associated with incompatible two-well problems. 

Our results are closely related to and build on ideas from the large literature on scaling behaviour of singular perturbation models inspired by shape-memory alloys. Let us particularly highlight the seminal works \cite{KM1, KM2} which for the first time rigorously proved that 
 branching structures can be energetically preferable to simple laminate structures in terms of scaling. Moreover, our geometrically linearized model and our Fourier analytic approach towards it are strongly inspired by the articles \cite{CO,CO1,KKO13,KO19}, see also \cite[Chapter 11]{B} for the modelling and \cite{KW14, KW16, PW22} for closely related models from compliance minimization. We also mention \cite{RT22,RT23,RT23a,RT24,IKRTZ24} for recent works set and analyzed in a Fourier-based framework. 
 For the analysis of quasiconvex hulls in geometrically linearized models which we use to illustrate the difference between the settings in which there are two consecutive degenerate symmetrized rank-one directions as well as the two-dimensional case, we particularly also mention \cite{CS13,CS15,CM22, CM22a} and point to \cite{BKreisbeck19,Merabet25} for recent developments on symmetric polyconvexity. We further highlight the works \cite{GZ23,GRTZ24} on scaling laws for low energy microstructures, \cite{CC14,CC15} for scaling laws in
 the geometrically nonlinear framework, the analysis of complex microstructures through regularity results \cite{RTZ19,RZZ18, RZZ19}, the investigation of nucleation results \cite{KK11,CDMZ20,AKKR24,RT23a,TZ24} as well as models inspired by the study of hysteresis \cite{Z14,CDZ17, CZ16}. Fine properties of minimizers have been studied in \cite{C1, CDKZ21, TS21, TS21a}. The effect of different forms of surface energy regularizations are the topic of \cite{RTTZ25}.  
 The latter result also rigorously relates \eqref{eq:minprob} with elastic energies perturbed by $\Vert D^2 u \Vert_{TV(\Omega)}$. 
We emphasize that this only represents a non-exhaustive list on the large literature on singular perturbation models inspired by shape-memory alloys and refer to the above articles and their bibliography for further references.

 \subsection{Outline of the remainder of the article}

The remainder of the article is structured as follows. In Section \ref{sec:prelim} we begin by introducing important notation and structural results. In particular, we explain why the case in which there are two subsequent degenerate symmetrized rank-one directions is different from the setting outlined above. We also recall the central Fourier analytical ideas to separate the different energetic contributions in phase space and adapt these to our setting with gauges. Next, in Section \ref{sec:proofs} we then present the arguments for the lower scaling bounds. In order to introduce the central ideas, we here first discuss the specific cases outlined in the introduction before proving the general statement from Theorems \ref{mainth:md} and \ref{mainth:mdperiodic} in Section \ref{sec:mwell}. In Section \ref{sec:upper} we then complement our discussion
 on lower bounds by a matching upper bound construction in the case \eqref{eq:choiceofmatrices} and prove Theorem~\ref{thm:mainresultupper}. 
Here a combination of different branching constructions becomes necessary. At the end of this section, we argue heuristically that also in the more general case one may hope that our lower scaling bounds are sharp -- we however expect that a rigorous justification of this requires overcoming substantial additional technical challenges. Finally, in the Appendix, for the convenience of the reader, we provide proofs for auxiliary results which have appeared in similar forms in the literature.

\section{Preliminaries}
\label{sec:prelim}

We begin by collecting various auxiliary results. In Section \ref{sec:compat}, we first recall the notion of compatibility
 while structure results on lamination convex hulls with consecutive degenerate symmetrized rank-one directions are deduced in Section \ref{sec:struc}. Next, in Section \ref{sec:boundary} we outline how the Dirichlet boundary conditions are connected with the elastic energy. Finally, we collect various Fourier-based perspectives on the phase-space distribution of elastic and surface energies in Section \ref{sec:Fourier}.

\subsection{Compatibility of wells}
\label{sec:compat}

A for our discussion central notion is that of the compatibility of two wells, which we hence recall.

\begin{definition}\label{def:Compatibility}
	We say that two distinct matrices $A,B \in \R^{d \times d}_{\sym}$ are \emph{compatible} if there exist two vectors $a,b \in \R^d \setminus \{ 0 \}$ such that
	\begin{align}
		A-B = a \odot b \defas \frac{1}{2} (a \otimes b + b \otimes a).\label{eq:def:compatible}
	\end{align}
	We say that
	
	\begin{itemize}
	\item there is \emph{one (degenerate) compatibility direction} if $a$ is  parallel or antiparallel to $b$, 
	\item and there are \emph{two (non-degenerate) compatibility directions} if $a$ and $b$ are neither parallel nor antiparallel.
	\end{itemize}
	
\end{definition}

As it is well-known in the literature (see, for instance, \cite{CS13,CS15}), compatibility can be characterized rather concisely.

\begin{lemma}\label{lem:compatibility}
	Let $A,B \in \R^{d \times d}_\sym$   with $A \neq B$. $A$ and $B$ are compatible if and only if $A-B$ has exactly two non-zero eigenvalues $\lambda_1$ and $\lambda_2$ with $\lambda_1 \lambda_2 < 0$ or exactly one non-zero eigenvalue.
	In the former case there are two non-degenerate compatibility directions, and in the latter one there is only one degenerate compatibility direction.
\end{lemma}
 We refer to \cite[Lemma~4.1]{Kohn91} for a precise discussion of the previous statement,
but, for the sake of completeness, we present a short proof for diagonal strains which we however postpone to the Appendix~\ref{sec:Appendix}. 
For the convenience of the reader, we further recall the definition of the symmetrized lamination convex hull of a set $K \subset \R^{d\times d}_{\sym}$.

\begin{definition}
\label{def:lamconvhull}
Let $K \subset \R^{d \times d}_{\sym}$. Then, the \emph{symmetrized lamination convex hull $K^{(\lc)}$} of the set $K$ is given by
\begin{align*}
K^{(\lc)}\defas \bigcup\limits_{j=0}^{\infty} K^{(j)},
\end{align*}
with $K^{(0)} = K$ and for $j\geq 1$ we define 
\begin{align*}
K^{(j)} &\defas \{ A \in \R^{d\times d}_{\sym}: \ A = \lambda B_1 + (1-\lambda) B_2 \mbox{ for } \lambda \in (0,1), \ B_1, B_2 \in K^{(j-1)} \\
& \qquad \mbox{ and } B_1 \text{ and } B_2  \mbox{ are compatible} \}.
\end{align*}
Elements of $K^{(j)}\setminus K^{(j-1)}$ will be referred to as \emph{laminates of order $j$}.
\end{definition}

\subsection{Structural results of the symmetric lamination convex hull}
\label{sec:struc}
 
The following result provides information on the topological structure of the lamination convex hull in two dimensions.
In particular, we show that a one-dimensional ``staircase type laminate'' does not exist in two dimensions and that thus in our discussion we necessarily have to consider three and higher dimensions. We point to \cite{CM22, CM22a} for an earlier discussion of semi-convex hulls in two-dimensional symmetrized elasticity.

\begin{proposition}[Topological property of the lamination convex hull I]\label{prop:topoprop}
	Let $K \subset \R^{2\times 2}_{\sym}$ be a set of wells such that 
			$K^{(2)} \setminus K^{(1)} \neq \emptyset$. Then,
			 $K^{(\lc)}$ contains a (two-dimensional) surface. 
\end{proposition}
 The proposition could be deduced via
\cite[Proposition~1]{CM22a}, where the symmetrized lamination convex hull for three wells has been characterized explicitly.
For convenience of the reader, we provide a short self-contained proof. 
\begin{proof}[Proof of Proposition~\ref{prop:topoprop}]
	If $K$ contains less than three matrices, the statement is trivial. In what follows, we thus assume that $\# K \geq 3$. We argue by contradiction and suppose that 
	$ K^{(2)} \setminus K^{(1)} \neq \emptyset $ while $K^{(\lc)}$ does \emph{not} contain a two-dimensional surface. 
Hence, there exist distinct matrices $A,B,C \in \R^{2 \times 2}_{\sym}$ such that $\lambda A + (1-\lambda) B \in K^{(1)}$ for all $\lambda \in (0,1)$, and some $\lambda_0 \in (0,1)$ such that for $Z \coloneq \lambda_0 A + (1-\lambda_0) B$ we have  $ \mu C + (1-\mu) Z \in K^{(2)} \setminus K^{(1)} $ for all $\mu \in (0,1)$.
 Notice that in two dimensions Lemma~\ref{lem:compatibility} characterizes the compatibility between two matrices by means of the determinant. 
This gives $
		\det (C - Z ) = \det (C - (\lambda_0 A + (1-\lambda_0) B   ) ) \leq 0 $. 
	This inequality cannot be strict due to the continuity of the determinant as $K^{(\lc)}$ is one-dimensional.
	Further, $\det(B-C) \neq 0 $ as otherwise Lemma~\ref{lem:compatibility} implies that $B_s \defas s Z + (1-s) B$ and $C_s = s Z + (1-s) C$
	are compatible for $s>0$ as
	$	\det (s Z + (1-s) B      -      s Z + (1-s) C) = 0$, contradicting the one-dimensionality of $K^{(\lc)}$.
	Since $A$ and $B$ are compatible, we have $\det (A-B) \leq 0$.
	If $\det(A-B) <0$, the definition of $Z$ implies that $\det (A-Z) < 0$, and therefore
	$\det (A - (\mu C + (1-\mu)Z )   ) <0$ for $\mu$ small, contradicting that $K^{(\lc)}$ does not contain a two-dimensional surface. 
	If $\det (A-B) = 0$, we can write $A-B = v \otimes v$ for a vector $v \neq 0$.
	Hence, we find that   
	\begin{align*}
		0 & = \det(C-Z) =  \det (C  - \lambda_0 A - (1-\lambda_0) B ) = \det (C-B) \det (\Id -  \lambda_0      ( ( C-B)^{-1} v) \otimes v  )  \\
		& =\det (C-B)  (1  -  \lambda_0       ( C-B)^{-1} v \cdot v ) . 
	\end{align*}
	The right-hand side is an affine function in $\lambda_0$ with $\det (C-B) \neq 0$. In particular, for $\epsilon>0$ small and $\lambda \in (\lambda_0 - \epsilon, \lambda_0)$ or $\lambda \in (\lambda_0, \lambda_0 + \epsilon)$, the expression $(1  -  \lambda_0       ( C-B)^{-1} v \cdot v ) $ will have a fixed sign. Thus, there exist further compatible matrices by choosing $\lambda \in (\lambda_0 - \eps , \lambda_0 + \eps) $ for some small $\eps$ for which it holds that $\det (C-B)  (1  -  \lambda_0       ( C-B)^{-1} v \cdot v ) \leq 0$, 
	contradicting again the fact that $K^{(\lc)}$ is one-dimensional.
\end{proof}

The following proposition implies that in any dimension ``staircase type" structures do not exist, provided that
there are consecutive degenerate compatibility directions in the symmetrized lamination convex hull. 
\begin{proposition}[Topological property of the lamination convex hull II]\label{prop:staircaseconstraint}
	Let $K \subset \R^{d \times d}_{\sym}$ be a set of wells such that there exist matrices $A,B\in K$ and a vector $n \in \R^d$ with
	$A-B = (\pm n) \odot n$.
	Setting $C_\lambda := \lambda A + (1-\lambda) B  $ for  $\lambda \in [0,1]$,
	we further assume that there exists $C \in K $, $\lambda_0 \in (0,1)$, and a vector
	$\nu \in \R^d$ with $\nu \neq n$ such that $C- C_{\lambda_0} = ( \pm \nu )   \odot \nu $.
	Then, $K^{(\lc)}$ contains a (two-dimensional) surface.
\end{proposition}

\begin{proof}
	We show that $C$ and $C_{\lambda_0 + s}$ are compatible for $s>0$ \emph{or} $s<0$ with $\vert s \vert$ sufficiently small.
	Indeed, according to Definition~\ref{def:Compatibility} we have
\begin{align*}
	C - C_{\lambda_0 + s} & = C - C_{\lambda_0} - s (A-B)  \\
	&=\begin{cases}
		( \pm \nu )   \odot \nu - s  (\pm n) \odot n = \pm (   \nu  + \sqrt{s}  n  ) \odot ( \nu - \sqrt{s} n  )  &\text{for } s \in [0, 1- \lambda_0], \\
		( \pm \nu )   \odot \nu - s  (\mp n) \odot n =  \pm  (  \nu  + \sqrt{ - s} n ) \odot ( \nu - \sqrt{ - s}   n   )   & \text{for } s \in [-\lambda_0, 0]. \\
	\end{cases}
	\end{align*}  
\end{proof}

We highlight that this more complex lamination convex hull is one of the central reasons why we do not discuss the case of consecutive degenerate symmetrized rank-one connections in this article.

\subsection{Exploiting the Dirichlet boundary conditions}
\label{sec:boundary}

We next turn to relating the Dirichlet boundary conditions and the elastic energy.
Consider the elastic energy of the form 
\begin{align}
	\Eelp{p} \defas	\int_\Omega \dist^p( \sym(\nabla u) , K_m) \di x.
\end{align}
 For convenience of notation, we also introduce the following two sets: 
\begin{align}
	S_t \defas  \{ l \in \{1,\dots ,d \} : l = 1 + \sum_{j=1}^{i-1} f(j) \text{ and }f(i) = t \text{ for some } i =1, \dots, m    \}    ,\quad t = 1,2. \label{eq:components}
\end{align}
 $S_1$ collects the ``effective indices'' associated with a degenerate symmetrized rank-one direction.
Since there are two ``effective indices'' for every non-degenerate symmetrized rank-one direction,
$S_2$ is the collection of the respective smaller indices.

With this notation fixed, we obtain the following relation between the deformation and the elastic energy which depends on the presence of a single degenerate or two non-degenerate rank-one directions.

\begin{lemma} \label{lem:boundarycondition}
	   Assume that $K_m$ is constructed as in \eqref{eq:kr}--\eqref{eq:constructiondim} for $f \in \mathcal{S}_m$ and $p \in [2,+\infty)$. 
Then,	there exists a constant $C>0$ such that for any
	  $u \in \mathcal{A}_F^{\dir} $ it holds that
	\begin{align}\label{eq:lemboundary}
		\int_{\Omega} \vert v_{ l } \vert^2 \di x \leq \begin{cases}
		C  \Big( \Eelp{2} + \sqrt{\Eelp{2}} \Big) & \quad \text{if }{ l } \in S_1, \\
		C  \Eelp{2}  & \quad \text{if } { l } \in S_2 ,\\
		\end{cases} 
	\end{align}
	where $v(x) = u(x) - Fx - b \in W_0^{1,p}(\Omega;\R^d) \subset H_0^1(\Omega;\R^d)$, and $S_j$, $j \in \{1,2\}$, as in \eqref{eq:components}. 
 \end{lemma}
 Notice that we have $\Eelp{2} \leq \Eel{m+1}{2}$ for any $\chi$ as in \eqref{eq:phaseindicprop} with equality if $\chi\colon \Omega \to \R^{d \times d}_{\sym}$ is the pointwise projection of $\sym(\nabla u)$ onto $K_m$.

\begin{proof}
	The definition of $\mathcal{A}_F^{\dir}$ in \eqref{def:boundarycond} implies that $v \in W_0^{1,p}(\Omega;\R^d) $. Let $\chi$ be the pointwise projection of $\sym(\nabla u)$ onto $K_m$ and let $\tilde \chi \defas \chi - F$. 

\emph{Step 1 (${ l } \in S_2$):} Assume that ${ l } = 1$. 
By the choice of $K_m$, we have $\tilde \chi_{11} = - \tilde \chi_{22}$. Using  the triangle inequality, we find that
\begin{align}\label{eq:trianglechoicewells}
	    \int_{\Omega} \vert \partial_1 v_1 + \partial_2 v_2 \vert^2 + \vert \partial_1 v_2 + \partial_2 v_1 \vert^2 \di x   &	   \leq 2 \int_{\Omega}  \sum_{j=1}^2 \vert \partial_j v_j - \tilde \chi_{jj} \vert^2 + \vert \partial_1 v_2 + \partial_2 v_1 \vert^2 \di x \notag \\
& \leq 2 \Eel{m+1}{2} . 
\end{align}
Recall that $v \in W_0^{1,p}(\Omega;\R^d) $ and that we can extend the function by $0$ to $(-1,2)^d$.
Setting $\nu_1 = ( 1, -1, 0, \dots, 0)^T$ and $\nu_2 = (1,1, 0, \dots, 0)^T$, we can use  a change of variables, the fundamental theorem of calculus, and derive for $i, j \in \{1,2\}$ with $i \neq j$ that
\begin{align}\label{useboundarycond}
	  \int_{(0,1)^2} \vert \nu_i \cdot v (x) \vert^2   \di x_1 \di x_2 & \leq   \int_{(-1,2)^2 } \vert \nu_i \cdot v ( \lambda \nu_j + r \nu_i ) \vert^2   \di r \di \lambda  \notag \\
	   &\leq  \int_{(-1,2)^2 } \left\vert \int_{-1}^r \nu_i \cdot \nabla v ( \lambda \nu_j + s \nu_i ) \nu_i \di s \right\vert^2   \di r \di \lambda \notag \\
	   & \leq C \int_{(-1,2)^2} \vert  \partial_2 v_2 + \partial_1 v_1     \vert^2  + \vert \partial_2 v_1 + \partial_1 v_2 \vert^2   \di r \di \lambda   \notag \\   
	   	   & \leq C \int_{(0,1)^2} \vert  \partial_2 v_2 + \partial_1 v_1     \vert^2  + \vert \partial_2 v_1 + \partial_1 v_2 \vert^2   \di x_1 \di x_2.
\end{align}
As $\nu_2 + \nu_1 = (2,0,\dots, 0)^T$, \eqref{eq:trianglechoicewells} and \eqref{useboundarycond} yield \eqref{eq:lemboundary}.
 The statement for a general ${ l }$ can be deduced by symmetry. 
 
	\emph{Step 2 (${ l } \in S_1$):} Assume that ${ l } =3$. Then, \eqref{def:admissiblelaminates} implies that the neighbouring effective indices lie in $S_2$. In particular, $ 1 \in S_2 $. 
Then, the fundamental theorem of calculus, the vanishing Dirichlet conditions and Jensen's inequality yield
 \begin{align}\label{eq:fundaclevertrick}
	  \int_{\Omega} \vert v_3 \vert^2 \di x 
	\leq \sum_{j=1}^2   \int_{\Omega} \vert \partial_{j} v_3 \vert^2 \di x \leq  \sum_{j=1}^2  \int_{\Omega} \vert \partial_{j} v_3 \vert^2+ \vert \partial_{3} v_j \vert^2   \di x.
 \end{align}
 For $v \in C_c^\infty(\Omega;\R^d)$ we can use integration by parts and the symmetry of second derivatives, which gives
 \begin{align}\label{eq:estimateforsymgrad}
	\sum_{j=1}^2 \int_{\Omega} \vert \partial_j v_3 \vert^2 + \vert \partial_3 v_j \vert^2 \di x &= \sum_{j=1}^2 \int_{\Omega} \vert \partial_j v_3 + \partial_3 v_j \vert^2 - 2 \partial_j v_3  \partial_3 v_j \di x \notag \\
	& = \sum_{j=1}^2 \int_{\Omega} \vert \partial_j v_3 + \partial_3 v_j \vert^2 \di x - 2 \int_\Omega (\partial_1 v_1 + \partial_2 v_2  ) \partial_3 v_3 \di x.
 \end{align}
 By density the previous identity also holds for $v \in W_0^{1,p}(\Omega;\R^d) \subset H^1_0 (\Omega;\R^d)$. 
Using once again that $\tilde \chi_{11} = - \tilde \chi_{22}$, we can proceed similarly to \eqref{eq:trianglechoicewells} and find together with Hölder's inequality that
 \begin{align}
& \quad \	\left\vert \int_\Omega (\partial_1 v_1 + \partial_2 v_2  ) \partial_3 v_3 \di x \right\vert \leq \left\vert \int_\Omega (\partial_1 v_1 - \tilde \chi_{11} + \partial_2 v_2 - \tilde \chi_{22} ) \partial_3 v_3 \di x \right\vert \notag \\
& \leq \left\vert \int_\Omega ( \partial_1 v_1 - \tilde \chi_{11} + \partial_2 v_2 - \tilde \chi_{22} ) ( \partial_3 v_3 - \tilde \chi_{33} ) \di x + \int_\Omega ( \partial_1 v_1 - \tilde \chi_{11} + \partial_2 v_2 - \tilde \chi_{22} )  \tilde \chi_{33} \di x \right\vert \notag \\
& \leq 2 \Eel{m+1}{2} + C \sqrt{\Eel{m+1}{2}}. \label{eq:cleverHölder}
 \end{align}
 By combining \eqref{eq:fundaclevertrick}, \eqref{eq:estimateforsymgrad}, and \eqref{eq:cleverHölder}, we obtain \eqref{eq:lemboundary} for ${ l } =3$. 
 The case of a general ${ l }$ can be deduced similarly by choosing a neighbouring index in $S_2 \neq \emptyset$.
In the case $K_{1,1}$, we have $S_2 = \emptyset$ and notice that the diagonal entry $ \tilde  \chi_{22} $ is zero. Then, one has to replace the sum $\sum_{j=1}^2$ in \eqref{eq:fundaclevertrick} with the single index $2$ and the index $3$ with $1$.
By proceeding along similar lines, this concludes the proof. 
\end{proof}

\subsection{The Fourier-based approach}
\label{sec:Fourier}
   Our goal is to prove the lower bounds in Theorems~\ref{mainth:3d}--\ref{mainth:mdperiodic}
via a Fourier-based framework. 
 At this point, we highlight once again that the construction of $K_m$ in \eqref{eq:kr}--\eqref{eq:constructiondim} is consistent with the wells $K_1, K_2$, and $K_3$ addressed in Theorems~\ref{thm:chanandconti}, \ref{mainth:3d} and \ref{mainth:4d}.
 Hence, the results in the remaining part of this section will be applicable in all of these cases. 
The general strategy is to control the $L^2$-norm of the phase indicator $\chi$, see
\eqref{eq:phaseindicprop}, in its representation as a Fourier series.
More specifically, one seeks to control the
Fourier coefficients via the elastic and surface energies.
As these bounds will hold uniformly among the class of phase indicators,
one obtains a lower bound for energies as defined in \eqref{eq:minprob}. 

\noindent For a $\mathbb{T}^d$-periodic function $f \in L^2\left(\mathbb{T}^d\right)$ the Fourier transform at $\xi \in \mathbb{Z}^d$ is given by  
$$
\mathcal{F}(f)(\xi)=\hat{f}(\xi)= \int_{\mathbb{T}^d} e^{- 2 \pi i \xi \cdot x} f(x) \, {\rm d}x ,
$$
 where the evaluations at $\xi$ are also called Fourier coefficients. Here $\mathbb{T}^d$ denotes the $d$-dimensional torus of length one.
In Section~\ref{sec:Fourier1} and Section~\ref{sec:Fourier2}, we reformulate the elastic energy as a Fourier series.
 Some preliminary estimates for low and high frequencies are provided in Section~\ref{sec:Fourier3} and Section~\ref{sec:Fourier4}.
 As these estimates do not allow for the optimal control in the case of higher order boundary data, a refined technique is presented in Section~\ref{sec:Fourier5}. 

\subsubsection{The symmetric gradient as a differential operator in Fourier space} \label{sec:Fourier1}
 Our next goal is to reformulate the elastic energy as a Fourier series. 
With a slight abuse of notation, we do not distinguish between the Dirichlet and periodic energies and, in this section, simply set  
\begin{align}\label{eq:Dirichletsetandperiodic}
	\Eelinf{m+1}{p} \defas \inf\limits_{u \in \mathcal{A}_F^{\dir} } \Eel{m+1}{p} \quad \text{ and } \quad \Eelinf{m+1}{p}{ } \defas \inf\limits_{u \in \mathcal{A}_F^{\per} } \Eel{m+1}{p} .
\end{align}
\begin{lemma}[Fourier representation of the elastic energy]\label{lem:Fourierchar}   
 Let $p \in [2,+\infty)$ and let $\chi \in L^\infty(\Omega;\R^{d \times d}_{\diag})$ be as in \eqref{eq:phaseindicprop}.  
  Then, it holds that
  \begin{align}
      &\quad \  \Eelinf{m+1}{p}^{2/p}   \notag \\
	   & \geq \sum\limits_{ \xi \in \mathbb{Z}^d \setminus \{0 \}}    \left( \sum_{i = 1}^d \sum_{j = i+1}^d \vert  \hat \xi_i^2 \F \tilde \chi_{jj}(\xi) +  \hat \xi_j^2 \F \tilde \chi_{ii}(\xi) \vert^2 +  \sum_{i=1}^d  \sum_{j \neq i }^d \sum_{k \neq j ,\ k \neq i}^d \hat \xi_k^2  \hat \xi_j^2 \vert \F \tilde \chi_{ii}(\xi) \vert^2 \right)  \notag \\
	   & \qquad \qquad +  \vert \F \tilde \chi (0) \vert^2,  \label{ineq:foursym} 
\end{align} 
where $\hat \xi = \vert \xi \vert^{-1} \xi $ denote the normalized frequencies and $\F \tilde \chi$ denotes the Fourier transform of \begin{align*}
	\tilde \chi \defas  
		\chi - F   \text{(viewed as a } \mathbb{T}^d \text{-periodic function)} .
\end{align*}
\end{lemma}

 We note that $\F \tilde{\chi} (\xi) = \F \chi(\xi)$ for all $\xi \neq 0$ and $\F \tilde{\chi}(0) = \F \chi(0) -F$.

The lemma is essentially proved in \cite[Lemma~3.1]{CO1} and \cite[Lemma~4.1]{KKO13} in dimension three and for $p = 2$.  We also refer to \cite[Lemma 3.1]{RRT23}, where the elastic energy in Fourier space has been characterized for more general differential operators, and to the simplification of the Saint-Venant compatibility conditions on diagonal matrices (see \eqref{eq:compatibilitystrain} below). 
We postpone a self-contained proof in a $d$-dimensional setting (which also deals with the case $p \in [2,\infty)$) to the Appendix~\ref{sec:Appendix}.

\subsubsection{ Incorporating the structure of the wells in Fourier space}\label{sec:Fourier2}

Heuristically, the elastic energy should vanish for e.g.\ laminate constructions where $\sym(\nabla u)$ only attains values in $K_m$. 
Hence, one should also observe this effect by viewing the elastic energy in its Fourier characterization.
 In this section, we take the structure of the wells into account and indeed show (quantitatively) that
the energy vanishes if the projection of $\sym(\nabla u)$ onto $K_m$ oscillates in the directions induced by the compatibility of the wells.

Recalling Definition~\ref{def:Compatibility}, we observe that the matrices $M_i$, $i = 1, \dots, m$, in \eqref{eq:constructiondim}
define some of the compatibility directions in $K^{(\lc)}$. We complement them by a second set of compatibility directions (which coincide with the previous ones in the degenerate case) by defining
\begin{align*}   (\tilde M_i)_{l \tilde l} := \begin{cases}
		1, &    l = 1 + \sum_{j=1}^{i-1} f(j)  \text{ and } l = \tilde l, \\
		\delta_{f(i)2}, & l = 2 + \sum_{j=1}^{i-1} f(j) \text{ and } l = \tilde l,  \\
		0, & \text{else}.
	\end{cases}
\end{align*}
With this definition in hand, we consider the sets $V_i \subset \R^d$ given by  
\begin{align}\label{def:vectorspaces}
	V_i \defas \begin{cases}
		\langle \bigcup_{j: f(j) = 1} \{ M_j \} \rangle  & \text{if }f(i) = 1, \\
		\langle   \{ M_i \} \rangle \cup \langle \{ \tilde  M_i \} \rangle  & \text{if }f(i) = 2, \\
	\end{cases}
\end{align}
for $i = 1,\dots, m$,
where we identify the diagonal entries of the matrices $M_i, \tilde M_i \subset \R^{d \times d}_\sym$ as a vector in $\R^d$, and 
$\langle S \rangle \subset \R^d $ denotes the smallest vector space generated by the set $S \subset \R^d $.

Eventually, the following lemma provides the lower bound for the elastic energy using the structure of the wells.
\begin{lemma}[Lower bound on the elastic energy in its Fourier characterization]\label{lem:elasticenergybound}
	 Let $p \in [2,+\infty)$.	   Assume that $K_m$ is constructed as in \eqref{eq:kr}--\eqref{eq:constructiondim}  for $f \in \mathcal{S}_m$. 
	Then, there exists a constant $c>0$ such that for all $i \in \{1,\dots,m\}$ with $l = l(i) = 1 + \sum_{j=1}^{i-1} f(j)$ it holds that 
	\begin{align}
	\Eelinf{m+1}{p}^{2/p}   & \geq  c \sum\limits_{ \xi \in \mathbb{Z}^d \setminus \{0 \}}   
		 \dist (\hat \xi, V_i)^{ 4  \delta_{f(i)1}   +   2\delta_{f(i)2}  }   \vert \hat \chi_{ll}(\xi) \vert^2 +  \vert \hat \chi (0)  - F\vert^2 .\label{eq:conclusionFourierchar} 
  \end{align} 
 \end{lemma}
	 Observe that the number of compatibility directions is reflected in the order of the multiplier.  
	 \begin{proof}
 We divide the proof into three steps. In the first one, we recall and introduce relevant notation and formulas
for a convenient presentation. In Step 2, we provide first preliminary estimates which are used for the conclusion in Step 3.
		
 \emph{Step 1 (Notation and formulas):} 
	 Let us first recall the definition of $S_1$ and $S_2$ in \eqref{eq:components}, 
	and observe that the (disjoint) union of $S_1$ and $S_2$ does not cover the set $\{1, \dots, d\}$ if $S_2 \neq \emptyset$ since $S_2$ only collects the smaller effective indices. 
	 It however holds that
\begin{align}
\label{eq:splitting}
\{1,\dots,d\} = S_1 \cup S_2 \cup \{j+1 \in \{1,\dots, d\}: \ j \in S_2\} =: S_1 \cup S_2 \cup \tilde{S}_2,
\end{align}	
which we will use frequently in the following argument.

 Throughout the proof the upper limit in the respective sums is not displayed and implicitly given by $d$ if an equality symbol appears in the lower limit. Moreover, as our proof will be given on each fixed frequency $\xi \in \Z^d$, with slight abuse of notation, we do not spell out the argument of the functions $\hat{\chi}_{ij}$ in what follows below but will understand it as evaluation on the frequency $\chi \in \Z^d$.

	We will rely on the lower bound in Lemma~\ref{lem:Fourierchar}, where we notice that Lemma \ref{lem:Fourierchar} already contains the coercivity estimate for the zero frequency.
	Hence, it suffices to consider the case $\xi \neq 0$, for which we have $\F \tilde{\chi}(\xi) = \hat{\chi}(\xi)$. 	
	 By the choice of wells it holds that $\chi_{j_0j_0} = - \chi_{( j_0 +1) (j_0+1) }$, provided that $j_0 \in S_2$, and hence $\hat \chi_{j_0j_0} = -\hat \chi_{( j_0 +1) (j_0+1) } $.   
Using this relation and the structure of $S_1$ and $S_2$ (see \eqref{eq:splitting}), the first term on the right-hand side of \eqref{ineq:foursym} can be rewritten as    

\begin{align}\label{eq:indiceschange1}
	&\quad \ \sum_{i = 1} \sum_{j = i+1} \vert  \hat \xi_i^2 \hat \chi_{jj} +  \hat \xi_j^2 \hat \chi_{ii} \vert^2  \notag \\
		&= \sum_{i \in S_1 } \sum_{j = i+1} \vert  \hat \xi_i^2 \hat \chi_{jj} +  \hat \xi_j^2 \hat \chi_{ii} \vert^2 
		+ \sum_{i \in S_2} \sum_{j = i+1} \vert  \hat \xi_i^2 \hat \chi_{jj} +  \hat \xi_j^2 \hat \chi_{ii} \vert^2 + \sum_{i \in \tilde{S}_2 } \sum_{j = i+1} \vert  \hat \xi_i^2 \hat \chi_{jj} +  \hat \xi_j^2 \hat \chi_{ii} \vert^2 \notag \\
			&= \sum_{i \in S_1 } \sum_{j = i+1} \vert  \hat \xi_i^2 \hat \chi_{jj} +  \hat \xi_j^2 \hat \chi_{ii} \vert^2 
		+  \sum_{i \in S_2} \vert   \hat \xi_{i+1}^2  - \hat \xi_i^2 \vert^2 \vert \hat \chi_{ii} \vert^2
 + \sum_{i \in S_2} \sum_{j > i+1} \vert  \hat \xi_i^2 \hat \chi_{jj} +  \hat \xi_j^2 \hat \chi_{ii} \vert^2 \notag	\\	
& \quad		+ \sum_{i \in S_2 } \sum_{j = i+2} \vert  \hat \xi_{i+1}^2 \hat \chi_{jj} +  \hat \xi_j^2 \hat \chi_{(i+1)(i+1)} \vert^2 \notag  \\
	&= \sum_{i \in S_1 } \sum_{j = i+1} \vert  \hat \xi_i^2 \hat \chi_{jj} +  \hat \xi_j^2 \hat \chi_{ii} \vert^2 \notag \\
	&\qquad + \sum_{i \in S_2} \vert   \hat \xi_{i+1}^2  - \hat \xi_i^2 \vert^2 \vert \hat \chi_{ii} \vert^2 + \sum_{i \in S_2} \sum_{j = i+2}   \left(   \vert  \hat \xi_i^2 \hat \chi_{jj} +  \hat \xi_j^2 \hat \chi_{ii} \vert^2   +  \vert  \hat \xi_{i+1}^2 \hat \chi_{jj} -  \hat \xi_j^2 \hat \chi_{ ii } \vert^2   \right).
\end{align}
 In order to rewrite all of the previous sums only in terms of indices in $S_1$ and $S_2$, we observe that
\begin{align}\label{eq:indiceschange2}
	\sum_{i \in S_1 } \sum_{j = i+1} \vert  \hat \xi_i^2 \hat \chi_{jj} +  \hat \xi_j^2 \hat \chi_{ii} \vert^2   &  =   \sum_{i \in S_1 } \sum_{j \in S_1 , j > i} \vert  \hat \xi_i^2 \hat \chi_{jj} +  \hat \xi_j^2 \hat \chi_{ii} \vert^2   +      \sum_{j \in S_1 } \sum_{i \in S_2 , i > j } \vert  \hat \xi_j^2 \hat \chi_{ii} +  \hat \xi_i^2 \hat \chi_{jj} \vert^2    \notag \\
	& \qquad  +     \sum_{j \in S_1 } \sum_{i \in S_2, i > j } \vert  \hat \xi_j^2 \hat \chi_{ii} -  \hat \xi_{i+1}^2 \hat \chi_{jj} \vert^2  ,
\end{align}
where we used the relation $\chi_{j_0j_0} = - \chi_{( j_0 +1) (j_0+1) }$ for $j_0 \in S_2$ once again and interchanged indices $j$ and $i$ in the last two terms.
Similarly, we find that
\begin{align}\label{eq:indiceschange3}
	& \quad \sum_{i \in S_2} \sum_{j = i+2}   \left(   \vert  \hat \xi_i^2 \hat \chi_{jj} +  \hat \xi_j^2 \hat \chi_{ii} \vert^2   +  \vert  \hat \xi_{i+1}^2 \hat \chi_{jj} -  \hat \xi_j^2 \hat \chi_{ ii } \vert^2   \right) \notag  \\
	& =  \sum_{i \in S_2} \sum_{j \in S_2 , j > i}   \left(   \vert  \hat \xi_i^2 \hat \chi_{jj} +  \hat \xi_j^2 \hat \chi_{ii} \vert^2   +  \vert  \hat \xi_{i+1}^2 \hat \chi_{jj} -  \hat \xi_j^2 \hat \chi_{ ii } \vert^2   \right) \notag  \\ 
		& \qquad + \sum_{i \in S_2} \sum_{j \in S_1, j> i }   \left(   \vert  \hat \xi_i^2 \hat \chi_{jj} +  \hat \xi_j^2 \hat \chi_{ii} \vert^2   +  \vert  \hat \xi_{i+1}^2 \hat \chi_{jj} -  \hat \xi_j^2 \hat \chi_{ ii } \vert^2   \right) \notag  \\
			& \qquad +  \sum_{i \in S_2} \sum_{j \in S_2, j > i}   \left(   \vert  \hat \xi_i^2 \hat \chi_{jj} -  \hat \xi_{j+1}^2 \hat \chi_{ii} \vert^2   +  \vert  \hat \xi_{i+1}^2 \hat \chi_{jj} +  \hat \xi_{j+1}^2 \hat \chi_{ ii } \vert^2   \right) .
\end{align}
By combining \eqref{eq:indiceschange1}, \eqref{eq:indiceschange2}, and \eqref{eq:indiceschange3} it therefore holds that 
\begin{align}	
	\label{eq:rewriting1}	
	&\quad \ \sum_{i = 1} \sum_{j = i+1} \vert  \hat \xi_i^2 \hat \chi_{jj} +  \hat \xi_j^2 \hat \chi_{ii} \vert^2  \notag \\
	 &= \sum_{i \in S_1 } \sum_{j > i, j \in S_1} \vert  \hat \xi_i^2 \hat \chi_{jj} +  \hat \xi_j^2 \hat \chi_{ii} \vert^2  + \sum_{i \in S_2} \vert   \hat \xi_{i+1}^2  - \hat \xi_i^2 \vert^2 \vert \hat \chi_{ii} \vert^2   \notag \\
	& \qquad    + \sum_{i \in S_1 } \sum_{ j \in S_2 } \left( \vert  \hat \xi_i^2 \hat \chi_{jj} +  \hat \xi_j^2 \hat \chi_{ii} \vert^2   + \vert \hat \xi_{j+1}^2 \hat \chi_{ii} -    \hat \xi_i^2 \hat \chi_{jj}\vert^2   \right) \notag \\
	& \qquad + \sum_{i \in S_2} \sum_{j >i, j \in S_2 }   \Big(   \vert  \hat \xi_i^2 \hat \chi_{jj} +  \hat \xi_j^2 \hat \chi_{ii} \vert^2   +  \vert  \hat \xi_{i+1}^2 \hat \chi_{jj} -  \hat \xi_j^2 \hat \chi_{ii} \vert^2  \notag \\
	& \qquad \qquad  \qquad \qquad   +    \vert   \hat \xi_{j+1}^2 \hat \chi_{ii}- \hat \xi_i^2 \hat \chi_{jj}  \vert^2   +  \vert  \hat \xi_{i+1}^2 \hat \chi_{jj} +  \hat \xi_{j+1}^2 \hat \chi_{ii} \vert^2   \Big)  .  
	 	 \end{align}

		  \emph{Step 2 (Preliminary estimates):} 
We keep the first two terms on the right-hand side of \eqref{eq:rewriting1} and start manipulating the third one.
Notice that Young's inequality and a quadratic expansion imply that $2  \hat \xi_i^2  \vert  \hat \xi_j^2 -\hat \xi_{j+1}^2 \vert  \vert \hat \chi_{ii} \vert \vert \hat \chi_{jj} \vert \leq \frac{3}{2}    \hat \xi_i^4 \vert \hat \chi_{jj} \vert^2 + \frac{2}{3} \vert  \hat \xi_j^2 -\hat \xi_{j+1}^2 \vert^2 \vert \hat \chi_{ii} \vert^2 $ and $\vert \hat \xi_j^2 -\hat \xi_{j+1}^2 \vert^2 \leq \hat \xi_j^4 + \hat \xi_{j+1}^4 $.
Thus, we find that
\begin{align}
	& \quad \ \vert  \hat \xi_i^2 \hat \chi_{jj} +  \hat \xi_j^2 \hat \chi_{ii} \vert^2   + \vert \hat \xi_{j+1}^2 \hat \chi_{ii} -    \hat \xi_i^2 \hat \chi_{jj}\vert^2 \notag \\
	&=  2 \hat \xi_i^4 \vert \hat \chi_{jj}  \vert^2 + 2  \hat \xi_i^2  ( \hat \xi_j^2 -\hat \xi_{j+1}^2 )  \text{Re}[\overline{ \hat \chi_{ii}} \hat \chi_{jj}] + (\hat \xi_j^4 +\hat \xi_{j+1}^4 )\vert \hat \chi_{ii} \vert^2 \notag \\
	&\geq  \frac{1}{2} \hat \xi_i^4 \vert \hat \chi_{jj}  \vert^2  + \frac{1}{3} (\hat \xi_j^4 +\hat \xi_{j+1}^4 )\vert \hat \chi_{ii} \vert^2 . \label{eq:rewriting2}
\end{align}
In a similar spirit, one can estimate the fourth term on the right hand side of \eqref{eq:rewriting1}. Thus, \eqref{eq:rewriting1} and \eqref{eq:rewriting2} yield
\begin{align}
 \sum_{i = 1} \sum_{j = i+1} \vert  \hat \xi_i^2 \hat \chi_{jj} +  \hat \xi_j^2 \hat \chi_{ii} \vert^2 &\geq 
   \sum_{i \in S_2} \vert   \hat \xi_{i+1}^2  - \hat \xi_i^2 \vert^2 \vert \hat \chi_{ii} \vert^2   \notag \\
   & \qquad  + \sum_{i \in S_1 } \sum_{ j \in S_2 } \left(  \frac{1}{2} \hat \xi_i^4 \vert \hat \chi_{jj}  \vert^2  + \frac{1}{3} (\hat \xi_j^4 +\hat \xi_{j+1}^4 )\vert \hat \chi_{ii} \vert^2   \right) \notag \\
   &  \qquad   + \sum_{i \in S_2} \sum_{j \neq i, j \in S_2 }   \left(        \frac{1}{3}  ( \hat \xi_{j}^4 + \hat \xi_{j+1}^4) \vert \hat \chi_{ii}  \vert^2    \right)   . \label{eq:rewriting3}
\end{align} 
 Next, we deal with the second term on the right-hand side of \eqref{ineq:foursym}. 
 Using the fact that $\hat \chi_{j_0 j_0} = - \hat \chi_{(j_0 +1) (j_0 +1) }$ for $j_0 \in S_2$ once again, we have
\begin{align}
	   & \quad \ \sum_{i=1}  \sum_{j \neq i } \sum_{k \neq j ,\ k \neq i} \hat \xi_k^2  \hat \xi_j^2 \vert \hat \chi_{ii} \vert^2 \label{eq:rewriting4} \\
& =  \sum_{i \in S_1}  \sum_{j \neq i } \sum_{k \neq j ,\ k \neq i} \hat \xi_k^2  \hat \xi_j^2 \vert \hat \chi_{ii} \vert^2 + \sum_{i \in S_2}  \sum_{j \neq i } \sum_{k \neq j ,\ k \neq i} \hat \xi_k^2  \hat \xi_j^2 \vert \hat \chi_{ii} \vert^2 + \sum_{i \in S_2}  \sum_{j \neq {i+1} } \sum_{k \neq j ,\ k \neq {i+1}  } \hat \xi_k^2  \hat \xi_j^2 \vert \hat \chi_{ii} \vert^2. \notag
\end{align}
Our next goal is to show that we additionally control $\hat \xi_j^2 \vert \hat \chi_{ii} \vert^2 $ for $i \in S_2$ and $j \neq i$, while \eqref{eq:rewriting3} currently only controls $\hat \xi_j^4 |\hat \chi_{ii}|^2$.
To this end, we observe that for any $\hat \xi_j^2 \in [\frac{1}{2},1]$ it holds that $ 2 \hat \xi_j^4 \vert \hat \chi_{ii} \vert^2 \geq \hat \xi_j^2 \vert \hat \chi_{ii} \vert^2$.
 If $\hat \xi_j^2 < \frac{1}{2}$, we have $\sum_{k \neq j}\hat \xi_k^2 > \frac{1}{2}$, implying that
 $ \hat \xi_j^2 \vert \hat \chi_{ii} \vert^2 <  2 \hat \xi_j^2  (\sum_{k \neq j}\hat \xi_k^2 ) \vert \hat \chi_{ii} \vert^2 $.
Thus, a combination of \eqref{eq:rewriting3} and \eqref{eq:rewriting4} implies that there exists a constant $c>0$ such that
\begin{align}\label{eq:rewriting5}
	&\quad \ \sum_{i = 1} \sum_{j = i+1} \vert  \hat \xi_i^2 \hat \chi_{jj} +  \hat \xi_j^2 \hat \chi_{ii} \vert^2 + \sum_{i=1}  \sum_{j \neq i } \sum_{k \neq j ,\ k \neq i} \hat \xi_k^2  \hat \xi_j^2 \vert \hat \chi_{ii} \vert^2  \notag \\
	&\geq c \sum_{i \in S_1 } \sum_{ j \in S_2 } \left(   \hat \xi_i^2 \vert \hat \chi_{jj}  \vert^2  +  ( \hat \xi_j^4 +\hat \xi_{j+1}^4 ) \vert \hat \chi_{ii} \vert^2   \right) + c \sum_{i \in S_2} \sum_{j \neq i, j \in S_2 }   \left(          ( \hat \xi_{j}^2 + \hat \xi_{j+1}^2) \vert \hat \chi_{ii}  \vert^2    \right)
\notag  \\
	& \qquad    +	  \sum_{i \in S_2} \vert   \hat \xi_{i+1}^2  - \hat \xi_i^2 \vert^2 \vert \hat \chi_{ii} \vert^2  .  
   \end{align}

  \emph{Step 3 (Conclusion):} 
 Recall that any vector $v \in \R^d$ can be rewritten as $v = \sum_{i = 1}^d \langle v, \zeta_i \rangle \zeta_i$ for
an orthonormal basis $\{ \zeta_i \}_{i = 1,\dots, d}$ and that the Pythagorean theorem implies that the distance of $v$ to a vector space $\langle \zeta_{i_1}, \dots, \zeta_{i_k} \rangle$ for $\{ i_1, \dots , i_k \} \subset \{1,\dots , d\}$ with complement $P_c =   \{1,\dots , d\} \setminus \{ i_1, \dots , i_k \}  $ is given by 
$\dist^2 (v, \langle \zeta_{i_1}, \dots, \zeta_{i_k} \rangle) = \sum_{j \in P_c} \vert \langle v, \zeta_j \rangle  \vert^2$. 

Let $l \in S_1$. In view of \eqref{def:vectorspaces} and Young's inequality we obtain a constant $C>0$ such that 
\begin{align}
	\dist^4 (\hat \xi, V_l)   \vert \hat \chi_{ll} \vert^2 =  \Bigg(\sum_{j \in S_2} (\vert \hat \xi_j \vert^2 + \vert \hat \xi_{j+1} \vert^2) \Bigg)^2 \vert \hat \chi_{ll} \vert^2 \leq C \sum_{j \in S_2} ( \hat \xi_j^4 + \hat \xi_{j+1}^4)  \vert \hat \chi_{ll} \vert^2 . \label{eq:rewriting6}
\end{align}
 The case $l \in S_2$ is slightly more involved to derive a suitable formula for $\dist^2 (\hat \xi, V_l)   \vert \hat \chi_{ll} \vert^2$ and we distinguish between two cases.   
Let us assume that 
$ \vert \hat \xi_l + \hat \xi_{l +1}  \vert = \vert \langle \hat \xi , \tilde M_l \rangle \vert \leq \vert \langle \hat \xi , M_l \rangle \vert = \vert \hat \xi_l - \hat \xi_{l +1} \vert $, recalling that we identify the matrices $\tilde M_l$ and $M_l$ as vectors in $\R^d$, see \eqref{def:vectorspaces}. 
In this case, it holds that
\begin{align}\label{ineq:distcase2}
	\dist^2 (\hat \xi, V_l)   \vert \hat \chi_{ll} \vert^2 = \dist^2 (\hat \xi, \langle  \{ M_l \}  \rangle )   \vert \hat \chi_{ll} \vert^2 =  \bigg( \vert \langle \hat \xi, \tilde  M_l  \rangle \vert^2 + \sum_{j \neq l, \ j \neq l +1 } \vert \langle \hat \xi, e_j \rangle \vert^2   \bigg) \vert \hat \chi_{ll} \vert^2.
\end{align}
If $ \vert  \hat \xi_l - \hat \xi_{l +1}  \vert > 1/2$, then 
\begin{align}\label{ineq:distcase3}
	\vert \langle \hat \xi, \tilde  M_l  \rangle\vert^2  \vert \hat \chi_{ll} \vert^2  =  \frac{( \hat \xi_l + \hat \xi_{l +1}  )^2 ( \hat \xi_l - \hat \xi_{l +1}  )^2 }{( \hat \xi_l - \hat \xi_{l +1}  )^2 }  \vert \hat \chi_{ll} \vert^2  <  4 ( \hat \xi_l^2 - \hat \xi_{l +1}^2   )^2  \vert \hat \chi_{ll} \vert^2.
\end{align}
If $ \vert  \hat \xi_l - \hat \xi_{l +1}  \vert \leq 1/2$, we particularly have $  \vert  \hat \xi_l + \hat \xi_{l +1}  \vert \leq 1/2$. 
These inequalities also imply that $\vert \hat \xi_l \vert +  \vert \hat \xi_{l +1} \vert \leq 1/2$, which in turn yields $ \vert \hat \xi_l \vert^2 +  \vert \hat \xi_{l +1} \vert^2  \leq 1/2$. 
Using that $\vert \hat \xi \vert = 1$, we find that  $\sum_{i \notin \{ l, l+1\} } \hat \xi_i^2 > 1/2$. Hence, we can estimate
\begin{align}\label{ineq:distcase4}
	\vert \langle \hat \xi, \tilde M_l \rangle\vert^2  \vert \hat \chi_{ll} \vert^2 \leq 2 ( \vert \hat \xi_l\vert^2  + \vert \hat \xi_{l +1} \vert^2  )  \vert \hat \chi_{ll} \vert^2  \leq   2 \sum_{i \notin \{ l, l+1\} } \hat \xi_i^2   \vert \hat \chi_{ll} \vert^2.
\end{align}
A combination of \eqref{ineq:distcase2}, \eqref{ineq:distcase3} and \eqref{ineq:distcase4} leads to
\begin{align}
	  \dist (\hat \xi, V_l)^{ 2}   \vert \hat \chi_{ll} \vert^2 \leq   3 \sum_{j \in S_1}  \hat \xi_j^2  \vert \hat \chi_{ll} \vert^2 + 3 \sum_{j \in S_2 , j\neq l} (  \hat \xi_j^2 +  \hat \xi_{j+1}^2)  \vert \hat \chi_{ll} \vert^2    + 4 ( \hat \xi_l^2 - \hat \xi_{l +1}^2   )^2  \vert \hat \chi_{ll} \vert^2  . \label{eq:rewriting7}
\end{align}
One can proceed similarly if $\vert \langle \hat \xi , M_l \rangle \vert > \vert \langle \hat \xi , \tilde M_l \rangle \vert$, i.e., \eqref{eq:rewriting7} remains true in this case.
Thus, by \eqref{ineq:foursym}, \eqref{eq:rewriting5}, \eqref{eq:rewriting6}, and \eqref{eq:rewriting7}, we conclude the proof of 
\eqref{eq:conclusionFourierchar}. 
 \end{proof}

\subsubsection{Low frequency estimates}\label{sec:Fourier3}
As a corollary of Lemma~\ref{lem:elasticenergybound}, we can control frequencies of $\hat \chi_{ll}$ outside of a neighbourhood of the sets $V_i$ by the elastic energy.

\begin{corollary}[Coercivity bound]\label{cor:lowfrequencyelastic}
		 	Let $p \in [2,+\infty)$.   Assume that $K_m$ is constructed as in \eqref{eq:kr}--\eqref{eq:constructiondim}  for $f \in \mathcal{S}_m$.  There exists a constant $C>0$ such that for all $i \in \{1,\dots,m\}$ with $l = l(i) = 1 + \sum_{j=1}^{i-1} f(j)$ and $\kappa \in (0,1)$ it holds that 
	\begin{align*}
		\sum_{\substack{\xi \in \Z^d\\ \dist(\hat \xi; V_i) > \kappa }}\vert \hat \chi_{ll}(\xi) \vert^2 \leq C \kappa^{-4\delta_{1f(i)} - 2 \delta_{2f(i)}}   \Eelinf{m+1}{p}^{2/p} .
	\end{align*}
\end{corollary}
\begin{proof}
	The proof follows immediately from \eqref{eq:conclusionFourierchar} as
	\begin{align*}
		 \sum_{\substack{\xi \in \Z^d\\ \dist(\hat \xi; V_i) > \kappa }}\vert \hat \chi_{ll} (\xi) \vert^2 \leq  \kappa^{-4  \delta_{f(i)1}   -  2\delta_{f(i)2}}  \sum_{\substack{\xi \in \Z^d\\ \dist(\hat \xi; V_i) > \kappa }} \dist(\hat \xi; V_i)^{ 4  \delta_{f(i)1}   +   2\delta_{f(i)2}  }  \vert \hat \chi_{ll}(\xi) \vert^2.
	\end{align*}
\end{proof}
 We also need to control Fourier coefficients at low frequencies, where the elastic energy is not coercive. To this end, we will use the following lemma in the setting of prescribed Dirichlet data.
\begin{lemma}[Estimates for low frequencies using boundary conditions]\label{lem:lowfrequencies1}
		Let $p \in [2,+\infty)$ and let $S_j$, $j\in \{1,2\}$, be as in \eqref{eq:components}. 	 	   Assume that $K_m$ is constructed as in \eqref{eq:kr}--\eqref{eq:constructiondim}  for $f \in \mathcal{S}_m$. Moreover, consider the Dirichlet setting, i.e., the class of admissible displacements is given by $\mathcal{A}_F^{\dir} $ in \eqref{eq:Dirichletsetandperiodic}. 
	Then, there exists a constant $C>0$ such that for every $\mu_l \in (1,+\infty)$ it holds that
	\begin{align}
				\sum_{ \{ \xi \in \Z^d : \, \vert \xi_l \vert \leq \mu_l \} } \vert \F \tilde \chi_{ll} (\xi) \vert^2  \leq \begin{cases}
		C \mu_l^2 \Big( \Eelinf{m+1}{p}^{2/p}  + \Eelinf{m+1}{p}^{1/p} \Big) & \quad \text{if }l \in S_1, \\
		C  \mu_l^2 \Eelinf{m+1}{p}^{2/p}   & \quad \text{if } l \in S_2 . 
		\end{cases}  
	\end{align}
\end{lemma}
We emphasize that the proof in \cite[Lemma 4.2]{RT23} is different as we need to prove a control via the symmetric gradient instead of the full gradient.
More precisely, the boundary conditions cannot be exploited simply by using the fundamental theorem of calculus, which is why we need to employ the quantitative estimate from Lemma~\ref{lem:boundarycondition}.    
Moreover, this lemma refines the estimate from \cite{RT23} by proving a control on a whole tube instead of a localized one-dimensional line. 
 Finally, we mention that we resort to a real-space argument based on Lemma~\ref{lem:boundarycondition}.
We refer e.g.\ to
\cite[Lemma~3.7(a)]{RRT23} for a Fourier-based argument in a setting without the gauge invariance $\text{Skew}(d)$. 
 
\begin{proof}  
	By the direct method in the calculus of variations the minimization problem 
	$\Eelinf{m+1}{p} = \inf_{u \in \mathcal{A}_F^{\dir} } \Eel{m+1}{p}$ admits a minimizer $u^*\in \mathcal{A}_F^{\dir}$. We recall that 
	$\tilde \chi =  
		\chi - F  $ and 
	set $v^*(x) = u^*(x) - Fx - b$ such that $v^* \in W^{1,p}_0(\Omega;\R^d) \subset H^1_0 (\Omega;\R^d)  $ can be viewed as a periodic function on $\mathbb{T}^d$. 
	 Plancherel's formula and Hölder's inequality imply that 
	\begin{align}\label{eq:newchoiceofv}
  \sum_{\xi \in \Z^d} \vert \mathcal{F} (\partial_l v^*_l)( \xi) - \mathcal{F} ( \tilde \chi_{ll})( \xi)  \vert^2	\leq 	\Vert \sym (\nabla v^*) - \tilde \chi \Vert^2_{L^2( \mathbb{T}^d )} \leq \Eelinf{m+1}{p}^{2/p}  .
	\end{align}
  Using  Lemma~\ref{lem:boundarycondition}, Plancherel's formula, and Hölder's inequality once again, there exists a constant $C>0$ such that we have
\begin{align}\label{v_ellestimatenew}
	 \sum_{\xi \in \Z^d}  \vert \mathcal{F} ( v^*_l )( \xi)   \vert^2 = \int_{\Omega}  \vert v^*_l \vert^2 \di x \leq \begin{cases}
		C  ( \Eelinf{m+1}{p}^{2/p}  + \Eelinf{m+1}{p}^{1/p} ) & \  \text{if }l \in S_1 ,\\
		C  \Eelinf{m+1}{p}^{2/p}  & \ \text{if } l \in S_2 . 
		\end{cases}           
\end{align}
Using  \eqref{eq:newchoiceofv}, \eqref{v_ellestimatenew}, the triangle inequality, and the identitiy $ \mathcal{F} (\partial_l v^*_l ) (\xi)  =  2\pi i \xi_l \mathcal{F} (v^*_l ) (\xi) $, we find that    
\begin{align*}
&\quad \sum_{ \{ \xi \in \Z^d : \, \vert \xi_l \vert \leq \mu_l \} } \vert \mathcal{F} ( \tilde \chi_{ll} )( \xi)   \vert^2 \\
  & \leq 2	\sum_{ \{ \xi \in \Z^d : \, \vert \xi_l \vert \leq \mu_l \} }\vert \mathcal{F} ( \partial_l v^*_l )( \xi)    - \mathcal{F} ( \tilde \chi_{ll} ) ( \xi) \vert^2 + 8 \pi^2 \sum_{ \{ \xi \in \Z^d : \, \vert \xi_l \vert \leq \mu_l \} }\vert \xi_l \vert^2 \vert \mathcal{F} ( v^*_l  )( \xi)  \vert^2  \\
	& \leq \begin{cases}
		C \mu_l^2 \left( \Eelinf{m+1}{p}^{2/p}  + \Eelinf{m+1}{p}^{1/p} \right) & \quad \text{if }l \in S_1, \\
		C  \mu_l^2 \Eelinf{m+1}{p}^{2/p}   & \quad \text{if } l \in S_2 . 
		\end{cases}     
\end{align*}
Thus, we conclude the lemma.
\end{proof}

\subsubsection{High frequency estimates}\label{sec:Fourier4}
For high frequencies, we will resort to standard estimates in terms of  the surface energy, see e.g.\ \cite[Lemma~4.5]{RT23} and \cite[Lemma~4.3]{KKO13}.
To this end, we recall the definition of the surface energy from \eqref{def:phaseindic}.  
\begin{lemma}\label{lem:highfrequencyoutsidecones}
 	Let $\chi \in L^\infty(\Omega;\R^{d \times d}_{\rm diag})$ be as in \eqref{eq:phaseindicprop}. 
	Then, there exists a constant $C>0$ such that for all $\mu_1 >0$ and $j \in \{1, \dots , d \}$ it holds that 
	\begin{align}
		 \sum_{ \{ \xi \in \Z^d : \, \vert \xi \vert \geq \mu_1 \} } \vert \hat \chi_{jj} (\xi) \vert^2 \leq C \mu_1^{-1} ( \Esurf{m+1} + \Per ).
	\end{align} 
\end{lemma}
\begin{proof}
	The lemma is proved in \cite[Lemma~4.5]{RT23}, which itself is based on \cite[Lemma~4.3]{KKO13}.
\end{proof}

\subsubsection{Quantitative interaction of Fourier coefficients}\label{sec:Fourier5}
While it suffices to apply Lemma~\ref{lem:highfrequencyoutsidecones} in the case of a first order datum,
 one needs to relate the oscillations of the different components with each other in the case of data of higher lamination order. 
More precisely,
in Proposition~\ref{prop:conereduction} below we show \emph{quantitatively} that
high frequencies of the component that resembles the finest oscillation, already control 
the periodicity of the other components.
 This concept is not new and has already been employed for $T_N$-structures as well as in the context of higher order data.
We refer e.g.\ to \cite{RT22} for its first application in the context of the singularly perturbed Tartar square and to the recent result \cite[Corollary 2.7]{RTTZ25} that neglects gauge invariances.  As outlined in the introduction, the presence of gauges however gives rise to substantial additional difficulties which are connected with the presence of two different types of symmetrized rank-one directions. 

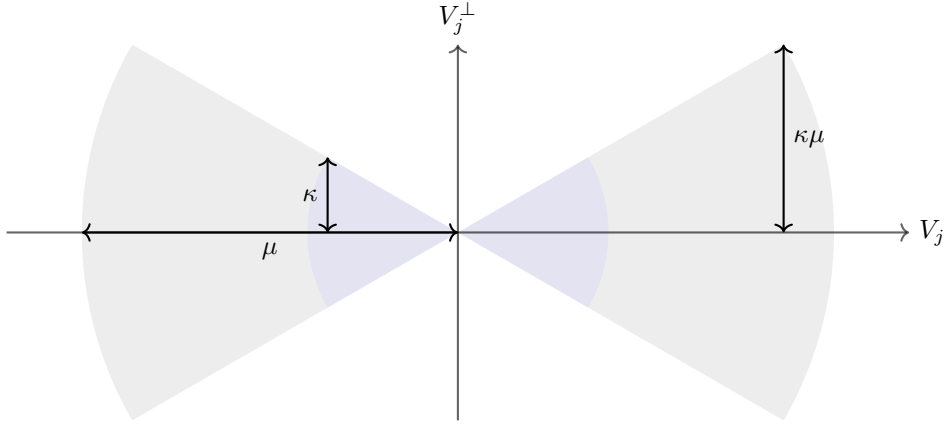
\begin{figure}
	\centering
	\begin{tikzpicture}[thick, scale= 1]

\def\mutikz{5}     
\def\kappatikz{30}

\fill[gray!20, opacity=0.7]
  (0,0) --
  ({\mutikz*cos(\kappatikz)},{\mutikz*sin(\kappatikz)}) arc[start angle=\kappatikz, end angle=-\kappatikz, radius=\mutikz] --
  (0,0) -- cycle;

  \fill[gray!20, opacity=0.7]
  (0,0) --
  ({\mutikz*cos(\kappatikz+180)},{\mutikz*sin(\kappatikz+180)}) arc[start angle=\kappatikz+180, end angle=-\kappatikz+180, radius=\mutikz] --
  (0,0) -- cycle;

\fill[blue!20, opacity=0.3]
  (0,0) --
  ({2*cos(\kappatikz)},{2*sin(\kappatikz)}) arc[start angle=\kappatikz, end angle=-\kappatikz, radius=2] --
  (0,0) -- cycle;

  \fill[blue!20, opacity=0.3]
  (0,0) --
  ({2*cos(\kappatikz+180)},{2*sin(\kappatikz+180)}) arc[start angle=\kappatikz+180, end angle=-\kappatikz+180, radius=2] --
  (0,0) -- cycle;

\draw[<->] ({5*cos(\kappatikz)},0) -- ({5*cos(\kappatikz)},{5*sin(\kappatikz)})  node[midway, right]{ $ \kappa \mu $}; 

\draw[<->] (-5,0) -- (0,0)  node[midway, below]{ $ \mu $}; 

\draw[<->] ({2*cos(\kappatikz+180)},0) -- ({2*cos(\kappatikz+180)},{-2*sin(\kappatikz+180)}) node[midway, left]{ $\kappa $};

\draw[->,opacity=0.6] (-6,0) -- (6,0) node[anchor=west] { $  $ };
\draw[->,opacity=0.6] (0,-2.5) -- (0,2.5) node[anchor=south] { $ $ };
\draw (6,0) node[right]  {$V_j$};
\draw (0,2.5) node[above]  {$V_j^\perp$};

\end{tikzpicture}
\caption{Visualization of the cone $C_{j, \kappa, \mu}$ in \eqref{eq:defcones}.
The horizontal axis indicates the vector space $V_j$ whereas the vertical axis represents the orthogonal vector space $V_j^\perp$.
The surface in blue depicts the intersection of the unit disc with the cone.}
\label{fig:defcone}
\end{figure}

Following the notation of \cite{RT23},
we introduce specific cones around the sets $V_i$ introduced in \eqref{def:vectorspaces}, allowing us to formulate the \emph{cone reduction} result.
For arbitrary but fixed $\kappa, \mu>0$, $j \in \{1,\dots ,m \}$, we set
\begin{align}\label{eq:defcones}
C_{j, \kappa, \mu} \defas \{\xi \in \Z^d: \dist(\hat \xi, V_j)  \leq \kappa, \ \vert \xi \vert \leq \mu\}.
\end{align}
Here, $\kappa$ should be interpreted as the angle width, while $\mu$ prescribes the length of the cones.
We refer to Figure~\ref{fig:defcone} for a visualization of these quantities.
Additionally, we associate Fourier multipliers $m_{j, \kappa,\mu}(D)$ satisfying $m_{j, \kappa,\mu}(\xi)=1$ for $\xi \in C_{j, \kappa,\mu}$ and $m_{j, \kappa,\mu}(\xi)=0$ for $\xi \in\Z^d \setminus C_{j,2 \kappa,2\mu}$.
Moreover, they can be chosen such that for any $p \in (1,+\infty)$ there exists a constant $C >0$ independent of $\kappa$ and $\mu$ such that for every $u\in L^p(\mathbb{T}^d)$ it holds that
\begin{align}\label{eq:decaymultiplay}
	 \Vert m_{j,\kappa,\mu}(D) u \Vert_{L^p} \leq C \Vert u \Vert_{L^p}.
\end{align}
Indeed, an explicit construction for such a multiplier is e.g.\ given in \cite[Section~4.3]{RT23}.

The main goal of this section is to prove the following result.

\begin{proposition}[Cone reduction argument]\label{prop:conereduction}
Let  $\lambda_i \in \R$ for $i\in \{1,\dots, m \}$, $g$ be a polynomial, and let 
$E_1 \subset \{1,\dots,m\}$, 
$E_2 \subset \{1,\dots,m\} \setminus E_1$, and $E_3 \subset \{1,\dots,m\} \setminus (E_1 \cup E_2)$.
Assume that 
\begin{itemize}
	\item[(a)] for $ l(i) = 1 + \sum_{j=1}^{i-1} f(j)$ the following polynomial relation holds:
 \begin{align*}
	\sum_{j \in E_3 } \lambda_j \tilde{\chi}_{l(j)l(j)} = g\bigg(\sum_{i \in E_1} \lambda_i \tilde\chi_{l(i)l(i)} \bigg) + \sum_{i \in E_2} \lambda_i \tilde\chi_{l(i)l(i)},
\end{align*}
\item[(b)] $V_{k_1} \cap V_{k_2} = \{0 \}$ for all $k_1, k_2 \in E_3$ with $ k_1 \neq k_2$,
\item[(c)] $V_i \cap V_j = \{ 0 \}$ for all $i\in E_1\cup E_2$ and all $j \in E_3$.
\end{itemize} 
Then, there exists a constant $C>0$ (depending on $\sup_{k \in  \{1,\dots,m\} } \vert \lambda_k \vert$, $\sup_{k \in  E_1 } \Vert \tilde \chi_{l(k)l(k)}\Vert_{L^\infty}$ and $g$) and some $M>0$ (depending on $g$) such that for 
all $\kappa_{l(k)} \in (0,1)$, $k \in E_1 \cup E_2 \cup E_3$,
 all $\mu_{l(i)}>0$, $i \in E_1 \cup E_2$,
and $\gamma_{l(j)} \geq \mu_{l(j)}$, $j \in E_3$, with
$\mu_{l(j)} \defas M \max_{i \in E_1 } \mu_{l(i)} \kappa_{l(i)}$
 it holds that   
\begin{align*}
& \quad \ \sum_{j \in E_3} \vert \lambda_j \vert \|\tilde{\chi}_{l(j)l(j)} - m_{l(j),\kappa_{l(j)}, \mu_{l(j)}}(D) \tilde{\chi}_{l(j)l(j)}\|_{L^2}\\
 & \leq C \sum_{j \in E_3} \|\tilde{\chi}_{l(j)l(j)} - m_{l(j),\kappa_{l(j)}, \gamma_{l(j)}}(D) \tilde{\chi}_{l(j)l(j)}\|_{L^2} \\
 & \quad + C \sum_{i \in E_1} \| \tilde\chi_{l(i)l(i)} - m_{l(i), \kappa_{l(i)}, \mu_{l(i)}}(D) \tilde\chi_{l(i)l(i)} \|_{L^2}^{1/2} \\
 & \quad + C \sum_{i \in E_2} \| \tilde\chi_{l(i)l(i)} - m_{l(i), \kappa_{l(i)}, \mu_{l(i)}}(D) \tilde\chi_{l(i)l(i)} \|_{L^2}  .
\end{align*}
\end{proposition}

Let us comment on this result. Proposition \ref{prop:conereduction} will play a central role for us in reducing the size of the possible cones of phase-space concentration. We will apply it with $E_3$ consisting of a single index (such that (b) is trivially satisfied), but state it in a more general form for the sake of completeness. The condition (c) is crucial, as it ensures that none of the cones associated with the support of elements in $E_1$ and $E_2$ has non-trivial overlap with the ones in $E_3$. This then allows us to reduce the size of the support of the cones associated with the $E_3$ indices iteratively. Here the linear components play a similar role as in \cite[Corollary 2.7]{RTTZ25} in that they do not contribute to the new cone support associated with $E_3$ ($\mu_{l(i)}$ is a free parameter for $i \in E_2$), while the ones from $E_1$ are ``thickened up'' by the nonlinear relation encoded by $g$ through convolution. This thickening up, then determines the size of the new support of the reduced $E_3$ cones as illustrated in the definition of the new cut-off index $\mu_{l(j)}$.

\begin{proof}

The proof follows along the same lines as that of \cite[Corollary 2.7]{RTTZ25} with the main difference that in the present article our cones are not necessarily cones around \emph{one-dimensional} but around potentially \emph{higher-dimensional} vector spaces. Due to the orthogonality of the involved vector spaces, the proof strategy however remains essentially unchanged.
 As in \cite{RTTZ25}, we proceed in several steps.

\emph{Step 1: A first nonlinear relation.}  The first step is similar to \cite[Lemma 2.6]{RTTZ25} and \cite[Lemma~4.6]{RT23}.
Given a polynomial of the form $g(x) = \sum_{k = 0}^N c_k x^k$ and values $a,b \in \R$, we have
\begin{align*}
	g(a) - g(b) = (a-b) G(a,b) \defas (a-b) \sum_{k=0}^{N} c_k \sum_{l=1}^k a^{k-l} b^{l-1}.
\end{align*} 
By Hölder's inequality (with powers $3/2$ and $3$) and the interpolation inequality for Lebesgue spaces (with $L^3 \subset L^2 \cap L^6$) we infer for maps $f_1 , f_2 \colon \mathbb{T}^2 \to \R$ that
\begin{align*}
  \|g(f_2 ) - g(f_1)\|_{L^2} \leq \|f_2  - f_1 \|_{L^{3}} \| G (f_2 , f_1 ) \|_{L^{6}} \leq \|f_2  - f_1 \|_{L^{2}}^{1/2} \|f_2  - f_1 \|_{L^6}^{1/2}  \| G (f_2 , f_1 ) \|_{L^{6}} .
\end{align*}
Then, the choices
 $f_1 = \sum_{i \in E_1} \lambda_i m_{l(i), \kappa_{l(i)}, \mu_{l(i)}}(D) \tilde\chi_{l(i)l(i)}$ and
$f_2 = \sum_{i \in E_1} \lambda_i \tilde\chi_{l(i)l(i)} $, the triangle inequality, Hölder's inequality, \eqref{eq:decaymultiplay}, the uniform boundedness of $\tilde \chi_{l(j)l(j)}$, and the concavity of $\vert \cdot \vert^{1/2}$
lead to
\begin{align}\label{ineq:firstpolynomialrelation}
& \quad \ \Bigg\|g \bigg(\sum_{i \in E_1} \lambda_i \tilde\chi_{l(i)l(i)} \bigg) + \sum_{i \in E_2} \lambda_i  \tilde\chi_{l(i)l(i)} - g\bigg(\sum_{i \in E_1} \lambda_i m_{l(i), \kappa_{l(i)}, \mu_{l(i)}}(D) \tilde\chi_{l(i)l(i)} \bigg) \notag\\
& \qquad \qquad \qquad \qquad - \sum_{i \in E_2} \lambda_i m_{l(i), \kappa_{l(i)},\mu_{l(i)}}(D)  \tilde\chi_{l(i)l(i)} \Bigg\|_{L^2}  \notag \\
&\leq C \sum_{i \in E_1} \| \tilde\chi_{l(i)l(i)} - m_{l(i), \kappa_{l(i)}, \mu_{l(i)}}(D) \tilde\chi_{l(i)l(i)} \|_{L^2}^{1/2} + C \sum_{i \in E_2} \| \tilde\chi_{l(i)l(i)} - m_{l(i), \kappa_{l(i)}, \mu_{l(i)}}(D) \tilde\chi_{l(i)l(i)} \|_{L^2} ,
\end{align}
where we note that the constant depends on $\sup_{i \in E_1 \cup E_2} \vert\lambda_i \vert$, $g$ and $\sup_{i \in E_1} \Vert \tilde \chi_{l(i)l(i)} \Vert_{L^\infty}$.

\emph{Step 2: A triangle inequality argument.}
Proceeding as in the proof of \cite[Corollary 2.7]{RTTZ25},
the triangle inequality, the relation in (a), and \eqref{ineq:firstpolynomialrelation} imply that
\begin{align}\label{eq:Step2}
& \quad \ \Bigg\|  \sum_{j \in E_3} \lambda_j  m_{l(j),\kappa_{l(j)}, \gamma_{l(j)}}(D) \tilde{\chi}_{l(j)l(j)} -  g\bigg(\sum_{i \in E_1} \lambda_i m_{l(i),\kappa_{l(i)},\mu_{l(i)}}(D) \tilde{\chi}_{l(i)l(i)} \bigg) \notag \\
& \qquad \qquad \qquad \qquad - \sum_{i \in E_2} \lambda_i m_{l(i), \kappa_{l(i)},\mu_{l(i)}}(D)  \tilde\chi_{l(i)l(i)}\Bigg\|_{L^2} \notag \\
& \leq C \sum_{j \in E_3} \|  m_{l(j),\kappa_{l(j)}, \gamma_{l(j)}}(D) \tilde\chi_{l(j)l(j)} -  \tilde\chi_{l(j)l(j)} \|_{L^2}  \notag\\
& \quad  +  \Bigg\| \sum_{j \in E_3}  \lambda_j \tilde\chi_{l(j)l(j)} - g\bigg(\sum_{i \in E_1} \lambda_i m_{l(i), \kappa_{l(i)},\mu_{l(i)}}(D) \tilde{\chi}_{l(i)l(i)}\bigg) - \sum_{i \in E_2} \lambda_i m_{l(i), \kappa_{l(i)},\mu_{l(i)}}(D)  \tilde\chi_{l(i)l(i)} \Bigg\|_{L^2}\notag \\
& \leq C \sum_{j \in E_3}  \| m_{l(j),\kappa_{l(j)}, \gamma_{l(j)}}(D) \tilde\chi_{l(j)l(j)} -   \tilde\chi_{l(j)l(j)} \|_{L^2} + C \sum_{i \in E_1} \| \tilde\chi_{l(i)l(i)} - m_{l(i), \kappa_{l(i)}, \mu_{l(i)}}(D) \tilde\chi_{l(i)l(i)} \|_{L^2}^{1/2} \notag \\
& \quad + C \sum_{i \in E_2} \| \tilde\chi_{l(i)l(i)} - m_{l(i), \kappa_{l(i)}, \mu_{l(i)}}(D) \tilde\chi_{l(i)l(i)} \|_{L^2} ,
\end{align}
where the constant also depends on $\sup_{i \in E_3} \vert\lambda_i \vert$.

\emph{Step 3: Disjoint supports and a nonlinear relation in Fourier space.}
In our final step, we apply the estimate from Step 2 as well as further disjointness of the supports of the relevant multipliers to deduce the desired estimate.

To this end, let us consider the orthogonal complement $   ( \bigoplus_{i\in E_1 \cup E_2 } V_i)^\perp \subset \R^d $  and denote by $\xi'$ the components of $\xi \in \R^d$ which are orthogonal to $\bigoplus_{i\in E_1 \cup E_2 } V_i$.
Moreover, let $m_{\{|\xi'|\geq \mu_{l(j)}\}}(D) f \defas \mathcal{F}^{-1} \chi_{\{|\xi'|\geq \mu_{l(j)}\}} \hat{f}$ denote the Fourier multiplier associated with a smoothed out version of the characteristic function $\chi_{\{|\xi'|\geq \mu_{l(j)} \}}$ with $m_{\{|\xi'|\geq \mu_{l(j)} \}}(\xi) = 0$ for $\vert \xi' \vert \leq \mu_{l(j)}/2$. 
As $\gamma_{l(j)} \geq \mu_{l(j)}$ for $j \in E_3$, and the estimate
$ \vert m_{l(j),\kappa_{l(j)}, \mu_{l(j)}}(\xi) - m_{l(j),\kappa_{l(j)}, \gamma_{l(j)}}(\xi) \vert \leq \chi_{\{|\xi'|\geq \mu_{l(j)}\}}(\xi) m_{l(j),\kappa_{l(j)}, \gamma_{l(j)}}(\xi)  $  
yields
\begin{align}
\label{eq:triange_disjoint_FSupp}
\begin{split}
& \quad \ \sum_{j \in E_3} \|\lambda_j \tilde{\chi}_{l(j)l(j)} - \lambda_j  m_{l(j),\kappa_{l(j)}, \mu_{l(j)}}(D) \tilde{\chi}_{l(j)l(j)}\|_{L^2}\\
&\leq  \sum_{j \in E_3}  \|\lambda_j \tilde{\chi}_{l(j)l(j)} - \lambda_j m_{l(j),\kappa_{l(j)}, \gamma_{l(j)}}(D) \tilde{\chi}_{l(j)l(j)}\|_{L^2}\\
& \quad 
+  \sum_{j \in E_3}  \|m_{\{|\xi'|\geq \mu_{l(j)}\}}(D) \lambda_j  m_{l(j),\kappa_{l(j)}, \gamma_{l(j)}}(D) \tilde{\chi}_{l(j)l(j)}\|_{L^2}.
\end{split}
\end{align}

We continue by discussing the second contribution on the right-hand side of \eqref{eq:triange_disjoint_FSupp}. 
We claim that there exists some $M>0$ (depending on $g$) such that for $\mu_{l(j)} \defas M \max_{i \in E_1 } \mu_{l(i)} \kappa_{l(i)}$ it holds that
\begin{align}\label{eq:fouriersup}
m_{\{|\xi'|\geq \mu_{l(j)}\}}(D)    g\bigg(\sum_{i \in E_1} \lambda_i m_{l(i),\kappa_{l(i)},\mu_{l(i)}}(D)  \tilde\chi_{l(i)l(i)} \bigg)  = 0.
\end{align}  
Indeed, as $g$ is a polynomial, the support of $g(\sum_{i \in E_1} \lambda_i m_{l(i),\kappa_{l(i)},\mu_{l(i)}}(D) \tilde{\chi}_{l(i)l(i)})$ can be estimated by a (finite) Minkowski sum of the support of $\sum_{i \in E_1} \lambda_i m_{l(i),\kappa_{l(i)},\mu_{l(i)}}(D) \tilde{\chi}_{l(i)l(i)}$ that depends on the degree of the polynomial $g$. Since the support of $m_{l(i),\kappa_{l(i)},\mu_{l(i)}}(D)  \tilde{\chi}_{l(i)l(i)}$ is given by a cone around $V_{i}$ with the largest radius $\kappa_{l(i)} \mu_{l(i)}$, see Figure~\ref{fig:defcone}, which is orthogonal to the $\xi'$ directions, the property $m_{\{|\xi'|\geq \mu_{l(j)} \}}(\xi) = 0$ for $\vert \xi' \vert \leq \mu_{l(j)}/2$ implies that we can find some $M>0$ such that the choice of $\mu_{l(j)}$ gives \eqref{eq:fouriersup}. 
As a consequence of the assumptions (c) and (b) and \eqref{eq:fouriersup}  we find that
\begin{align}\label{eq:support_orthogonal}
 &\quad \ \sum_{j \in E_3}  \| m_{\{|\xi'|\geq \mu_{l(j)} \}}(D) \lambda_j m_{l(j),\kappa_{l(j)}, \gamma_{l(j)}}(D) \tilde{\chi}_{l(j)l(j)}\|_{L^2} \notag \\
 & \leq  \Bigg\| m_{\{|\xi'|\geq \mu_{l(j)} \}}(D)  \Bigg( \sum_{j \in E_3}  \lambda_j m_{l(j),\kappa_{l(j)}, \gamma_{l(j)}}(D) \tilde{\chi}_{l(j)l(j)} - \sum_{i \in E_2} \lambda_i m_{l(i),\kappa_{l(i)},\mu_{l(i)}}(D)  \tilde\chi_{l(i)l(i)}   \Bigg) \Bigg\|_{L^2}  \notag \\
& \leq \Bigg\| m_{\{|\xi'|\geq \mu_{l(j)} \}}(D) \Bigg( \sum_{j \in E_3} \lambda_j  m_{l(j),\kappa_{l(j)}, \gamma_{l(j)}}(D) \tilde{\chi}_{l(j)l(j)} \notag  \\
& \qquad \qquad \qquad \qquad  -  g\bigg(\sum_{i \in E_1} \lambda_i m_{l(i),\kappa_{l(i)},\mu_{l(i)}}(D) \tilde{\chi}_{l(i)l(i)}\bigg)  - \sum_{i \in E_2} \lambda_i m_{l(i),\kappa_{l(i)},\mu_{l(i)}}(D)  \tilde\chi_{l(i)l(i)}  \Bigg)\Bigg\|_{L^2} \notag \\
& \leq \Bigg\|  \sum_{j \in E_3} \lambda_j  m_{l(j),\kappa_{l(j)}, \gamma_{l(j)}}(D) \tilde{\chi}_{l(j)l(j)} -  g\bigg(\sum_{i \in E_1} \lambda_i m_{l(i),\kappa_{l(i)},\mu_{l(i)}}(D) \tilde{\chi}_{l(i)l(i)} \bigg) \notag \\
& \qquad \qquad \qquad \qquad -\sum_{i \in E_2} \lambda_i m_{l(i),\kappa_{l(i)},\mu_{l(i)}}(D)  \tilde\chi_{l(i)l(i)}  \Bigg\|_{L^2}.
\end{align}
Eventually, the assumption (b), \eqref{eq:triange_disjoint_FSupp}, \eqref{eq:support_orthogonal}, and \eqref{eq:Step2} imply that  
\begin{align*}
& \quad \ \sum_{j \in E_3} \vert \lambda_j \vert \|\tilde{\chi}_{l(j)l(j)} - m_{l(j),\kappa_{l(j)}, \mu_{l(j)}}(D) \tilde{\chi}_{l(j)l(j)}\|_{L^2}\\
&\leq   C \sum_{j \in E_3} \|\tilde{\chi}_{l(j)l(j)} - m_{l(j),\kappa_{l(j)}, \gamma_{l(j)}}(D) \tilde{\chi}_{l(j)l(j)}\|_{L^2} \\
 &\quad  + \Bigg\|  \sum_{j \in E_3} \lambda_j  m_{l(j),\kappa_{l(j)}, \gamma_{l(j)}}(D) \tilde{\chi}_{l(j)l(j)} -  g\bigg(\sum_{i \in E_1} \lambda_i m_{l(i),\kappa_{l(i)},\mu_{l(i)}}(D) \tilde{\chi}_{l(i)l(i)} \bigg) \\
& \qquad \qquad \qquad \qquad  -\sum_{i \in E_2} \lambda_i m_{l(i),\kappa_{l(i)},\mu_{l(i)}}(D)  \tilde\chi_{l(i)l(i)}  \Bigg\|_{L^2} \\
& \leq C \sum_{j \in E_3} \|\tilde{\chi}_{l(j)l(j)} - m_{l(j),\kappa_{l(j)}, \gamma_{l(j)}}(D) \tilde{\chi}_{l(j)l(j)}\|_{L^2} + C \sum_{i \in E_1} \| \tilde\chi_{l(i)l(i)} - m_{l(i), \kappa_{l(i)}, \mu_{l(i)}}(D) \tilde\chi_{l(i)l(i)} \|_{L^2}^{1/2} \notag \\
& \quad + C \sum_{i \in E_2} \| \tilde\chi_{l(i)l(i)} - m_{l(i), \kappa_{l(i)}, \mu_{l(i)}}(D) \tilde\chi_{l(i)l(i)} \|_{L^2}  .
\end{align*}
This concludes the proof.
\end{proof}

\section{Proof of the lower bounds}\label{sec:proofs}
 
In this section, we present the proofs for the lower bound results from Theorems \ref{mainth:3d}-\ref{mainth:mdperiodic}. We will rely on the Fourier perspective developed in the previous section. A key difficulty will consist in the different treatment of the degenerate and non-degenerate symmetrized rank-one connections.

As a brief motivation, we explain the derivation of the lower bounds from Theorem~\ref{thm:CC14} via our Fourier analytical approach.
Recall that we view $\tilde \chi = \chi - F $
 as a $\mathbb{T}^d$-periodic function.
Given a boundary datum $F = \diag (\alpha,0)$ for $K_{1,1}$ and $\alpha \in (0,1)$ and $F = \diag(\alpha, - \alpha)$ for $K_{1,2}$, respectively,
any phase indicator satisfies
	\begin{align}\label{ineq:positivitysimple}
		\Vert \tilde \chi_{11} \Vert_{L^2}^2 = \Vert \chi_{11} - \alpha \Vert_{L^2}^2 \geq \min \{ \vert 1-\alpha \vert^2, \vert \alpha \vert^2 \} > 0 
	\end{align}
	as $\chi_{11}$ only attains the values $0$ and $1$. Thus, by Plancherel's formula, Lemma~\ref{lem:lowfrequencies1}, and Lemma~\ref{lem:highfrequencyoutsidecones}
	there exist constants $c>0$ and $C>0$ such that for any $\mu_1 >1 $ it holds that
	\begin{align*}
		& \quad \ 0 < c \leq \Vert \tilde \chi_{11} \Vert_{L^2}^2   \leq  	\sum_{ \{ \xi \in \Z^2 : \, \vert \xi_1 \vert \leq \mu_1 \} } \vert \F \tilde \chi_{11}  (\xi) \vert^2 +	 \sum_{ \{ \xi \in \Z^2 : \, \vert \xi \vert \geq \mu_1 \} } \vert \F \tilde \chi_{11} (\xi) \vert^2  \\
		& \leq C \begin{cases}
		\mu_1^2 ( \Eelinf{2}{p}^{2/p}  + \Eelinf{2}{p}^{1/p} ) +\mu_1^{-1} ( \Esurf{2} + \Per ) & \ \text{if } 1  \in S_1, \\
		 \mu_1^2 \Eelinf{2}{p}^{2/p} +\mu_1^{-1} ( \Esurf{2} + \Per )  & \  \text{if } 1 \in S_2 . 
		\end{cases}  
	\end{align*} 
	The previous estimate holds uniformly among all phase indicators. By \eqref{eq:minprob}, \eqref{eq:Dirichletsetandperiodic}, 
	and Young's inequality (in the form $ab \leq \delta^{-1} \tilde p^{-1} a^{\tilde p} + \delta \tilde q^{-1} b^{\tilde q} $) we obtain (for a potentially smaller $c>0$)
			\begin{align}\label{ineq:firsttimeyoung}
		  c (1 -  \mu_1^{-1}  \Per) \leq C \begin{cases}
		(\mu_1^p + \mu_1^{2p}+\epsilon^{-1} \mu_1^{-1})  \Energyinf{2}{p}{\dir}     & \ \text{if } 1  \in S_1, \\
			(\mu_1^p +\epsilon^{-1} \mu_1^{-1})  \Energyinf{2}{p}{\dir}   & \  \text{if } 1 \in S_2 . 
		\end{cases}  
	\end{align} 
	The choices $\mu_1 \sim \epsilon^{-1/(2p+1)}$ and $\mu_1 \sim \epsilon^{-1/(p+1)}$ in the former and latter case, respectively, lead to the desired scalings
	\begin{align}\label{ineq:firstorderp}
			\Energyinf{2}{p}{\dir} \gtrsim \begin{cases}
			 \epsilon^{2p/(2p+1)} & \text{for the wells } K_{1,1} \text{ with } F = \diag (\alpha,0) \text{ for } \alpha \in (0,1), \\
			 \epsilon^{p/(p+1)} &  \text{for the wells } K_{1,2} \text{ with } F = \diag (\alpha,-\alpha) \text{ for } \alpha \in (0,1) .  
		\end{cases}
	 \end{align}

In what follows below, we seek to argue similarly as for the sets $K_{1,1}$ and $K_{1,2}$ also for the higher order sets $K_m$ and for the case of boundary data of higher lamination order. In addition to the arguments from above, this will require the use of the cone reduction result from Proposition \ref{prop:conereduction}. As the argument for the general set $K_m$ is notationally rather heavy due to the presence of degenerate and non-degenerate symmetrized rank-one directions, we begin by introducing these ideas for the three-well problem (Theorem \ref{mainth:3d}) in Section \ref{sec:threewell}, then turn to the four-well problem both qualitatively and quantitatively (Theorem \ref{mainth:4d}) in Section \ref{sec:fourwell} and finally present the general $d$-dimensional argument (Theorems \ref{mainth:md} and \ref{mainth:mdperiodic}) in Section \ref{sec:mwell}.

\subsection{The three-well problem} 
\label{sec:threewell}
In this section we prove Theorem~\ref{mainth:3d}.   
To obtain a lower scaling bound for a Dirichlet datum of \emph{arbitrary} order one could follow the strategy which is
 presented in \eqref{ineq:positivitysimple}--\eqref{ineq:firstorderp}
 by analyzing a variable $\tilde \chi_{jj}$ which oscillates.
 A finer analysis, however, shows that these bounds are suboptimal if the order of the datum is larger than one. 
More precisely, for higher order data, one expects oscillations between more than two wells on different scales. More surface energy is needed to reduce the elastic energy, which suggest that the optimal scaling tends to zero slower than the bounds proved in \eqref{ineq:firstorderp}.
In order to capture this, we employ the technique presented in Section~\ref{sec:Fourier5} (Proposition \ref{prop:conereduction}), which relates the oscillations of different components with each other.
In essence, this technique describes quantitatively that oscillations with a well of higher order is already influenced by oscillations between wells of lower order.

\begin{proof}[Proof of Theorem~\ref{mainth:3d}] We split the proof into several steps and discuss the cases of first and second order boundary data and the settings with degenerate and non-degenerate symmetrized rank-one directions separately.

\emph{Step 1 (Proof of (i)):}
Let us first address the case \eqref{eq:choiceofmatrices}.
Due to \eqref{eq:3D_wells}, the boundary conditions are prescribed by the matrix $ F = \diag(\alpha, 0,0)$ for $\alpha \in (0,1)$.
Our goal is to follow the lines of the discussion above \eqref{ineq:firstorderp}.
The essential difference is that $\chi_{11}$ can attain the values $0$, $1/2$, and $1$, i.e.,
for $\alpha = 1/2$, there could be a set with positive Lebesgue measure such that $\tilde \chi_{11} = \chi_{11} - \alpha = 0$ and the estimate in \eqref{ineq:positivitysimple} becomes subtler.
This set, however, needs to be small for small elastic energy as the zero boundary conditions of $u^*_2$ for the optimal $u^*$ in \eqref{eq:Dirichletsetandperiodic} imply that there exists a constant $c>0$ such that
	\begin{align*}
		\Vert \tilde \chi_{11} \Vert_{L^2}^2 & \geq \int_{\{ \chi_{11} \in \{0,1\} \} } \vert \chi_{11} - \alpha \vert^2 \di x =  \min \{ \vert 1-\alpha \vert^2, \vert \alpha \vert^2 \} (\vert \Omega \vert - \vert \{ \chi_{11} = 1/2 \} \vert) \\
& = c \left(\vert \Omega \vert - \int_{\{ \chi_{22} = 1 \}} \chi_{22}  \di x \right) = c \left(\vert \Omega \vert - \int_{\Omega} \chi_{22} - \partial_2 u^*_2 \di x  \right) \\
& \geq c \left(\vert \Omega \vert - \vert \Omega \vert^{(p-1)/p}  \Eelinf{3}{p}^{1/p} \right) .
	\end{align*}

As both Lemma~\ref{lem:lowfrequencies1} and Lemma~\ref{lem:highfrequencyoutsidecones} do not depend on the third dimension and $\mu_1$ becomes large, we conclude as in \eqref{ineq:firstorderp}.
Part (i) in the case \eqref{eq:secondcase} follows along similar lines, where we have to replace $\tilde \chi_{22}$ by $\tilde \chi_{33}$ in the previous calculation.

\emph{Step 2 (Proof of (ii) in the case \eqref{eq:choiceofmatrices}):}
	According to \eqref{eq:3D_wells}, the boundary conditions are prescribed by the matrix $ F = \diag(1/2, \alpha, -\alpha)$ for $\alpha \in (0,1)$.
	To derive the optimal scaling, we start by analyzing the variable $\tilde \chi_{22}$, which captures the coarsest scale of oscillations. (One could also estimate $\tilde \chi_{33}$ in this case.)
	As $\chi_{22}$ only attains the values $0$ and $1$ it holds that
	\begin{align}
		\Vert \tilde \chi_{22} \Vert_{L^2}^2 \geq \int_\Omega \vert \chi_{22} - \alpha \vert^2 \di x =  \min \{ \vert 1-\alpha \vert^2, \vert \alpha \vert^2 \}  . \label{eq:3dproof1W1}
	\end{align}
	By the construction of the Fourier multiplier in the discussion below \eqref{eq:defcones}, the multiplier $m_{2,\kappa_2, \mu_2}(D)$ for $\kappa_2 \in (0,1)$ and $\mu_2 > 1$ concentrates in a cone around the union of the vector spaces generated by $e_2 + e_3$ and $e_3 - e_2$. 
	The Pythagorean theorem yields $m_{2,\kappa_2, \mu_2}(D) m_{\{|\xi_2|\geq 4\mu_2\}}(D) \tilde \chi_{22} = 0$, where
	 $m_{\{|\xi_2|\geq 4 \mu_2\}}(D) f:= \mathcal{F}^{-1} \chi_{\{|\xi_2|\geq 4 \mu_2\}}(k) \hat{f}$ denotes the Fourier multiplier associated with a smoothed out version of the characteristic function $\chi_{\{|\xi_2|\geq 4 \mu_2 \}}$ with value zero for $\vert \xi_2 \vert < 2\mu_2 $. 
Thus, Plancherel's formula, the triangle inequality, \eqref{eq:decaymultiplay}, and Lemma~\ref{lem:lowfrequencies1} (as $2 \in S_2$) imply that  
	\begin{align}
		& \quad \ \Vert \tilde \chi_{22} \Vert_{L^2}  \leq  \Vert \tilde \chi_{22} - m_{2,\kappa_2, \mu_2}(D)\tilde \chi_{22} +  m_{2,\kappa_2, \mu_2}(D) \tilde \chi_{22}- m_{2,\kappa_2, \mu_2}(D) m_{\{|\xi_2|\geq 4\mu_2\}}(D) \tilde \chi_{22}\Vert_{L^2}  \notag \\
		&\leq \Vert  \tilde \chi_{22} - m_{2,\kappa_2, \mu_2}(D) \tilde \chi_{22} \Vert_{L^2} + C \mu_2  \Eelinf{3}{p}^{1/p} . \label{eq:3dproof2W1}
	\end{align} 
	Next, we notice that $\tilde \chi_{22} = g(\tilde \chi_{11}) $ for $g(t) = -4t^2 + 1 - \alpha$. By \eqref{def:vectorspaces} it holds that $V_1 \cap V_2 = \{0 \}$.
	Hence, we can apply Proposition~\ref{prop:conereduction} (with $ \gamma_1 = \mu_1 \geq \mu_2$) and derive for any $\mu_1>0$ and $\kappa_1 \in (0,M^{-1})$ with $\mu_2 = M \kappa_1 \mu_1 \leq \mu_1 $ that
	\begin{align}
		& \quad \ \Vert  \tilde \chi_{22} - m_{2,\kappa_2, \mu_2}(D) \tilde \chi_{22} \Vert_{L^2} \notag \\
		& \leq C   \Vert  \tilde \chi_{11} - m_{1,\kappa_1,\mu_1} (D) \tilde \chi_{11} \Vert_{L^2}^{1/2} + C \Vert \tilde \chi_{22} - m_{2,\kappa_2,\mu_1}(D) \tilde \chi_{22} \Vert_{L^2}  . \label{eq:3dproof3W1}
	\end{align}
	Using Corollary~\ref{cor:lowfrequencyelastic} and Lemma~\ref{lem:highfrequencyoutsidecones}, 
 we obtain that
\begin{align}
\begin{split}
\label{eq:intermediateW1}
 \Vert  \tilde \chi_{11} - m_{1,\kappa_1,\mu_1} (D) \tilde \chi_{11} \Vert_{L^2}^{1/2}
& \leq   C \kappa_1^{-1}  \Eelinf{3}{p}^{1/(2p)} + C  \mu_1^{-1/4} ( \Esurf{3}^{1/4} + \Per^{1/4} ),\\
 \Vert \tilde \chi_{22} - m_{2,\kappa_2,\mu_1}(D) \tilde \chi_{22} \Vert_{L^2}  
 & \leq  C \kappa_2^{-1}  \Eelinf{3}{p}^{1/p}+ C  \mu_1^{-1/2} ( \Esurf{3}^{1/2} + \Per^{1/2} ).
 \end{split}
\end{align}
	
Hence, by combining \eqref{eq:3dproof1W1}--\eqref{eq:intermediateW1} and employing the concavity of $\vert \cdot \vert^{1/2}$, we find that
	\begin{align}
			&  \quad \ \min \{ \vert 1-\alpha \vert, \vert \alpha \vert \}  \notag \\
			& \leq   C \mu_2  \Eelinf{3}{p}^{1/p} + C \kappa_2^{-1}  \Eelinf{3}{p}^{1/p} + C  \mu_1^{-1/2} ( \Esurf{3}^{1/2} + \Per^{1/2} ) \notag \\
			& \quad \  +    C \kappa_1^{-1}  \Eelinf{3}{p}^{1/(2p)} + C  \mu_1^{-1/4} ( \Esurf{3}^{1/4} + \Per^{1/4} ) .
	\end{align}
		 In order to obtain the same powers of the energy on the right-hand side, we employ Young's inequality, where the arising constant can be absorbed by the left-hand side, similarly to \eqref{ineq:firsttimeyoung}. 
	 Recalling \eqref{eq:minprob} and \eqref{eq:Dirichletsetandperiodic}, there exists a constant $c>0$ such that
	\begin{align}\label{eq:optimizethefirst}
			  c (1 - \mu_1^{-1}  \Per) \leq   C (\mu_2^p + \kappa_1^{-2p} +\kappa_2^{-p} +  \epsilon^{-1} \mu_1^{-1} ) \Energyinf{3}{p}{\dir} .
	\end{align}
	 	To derive the optimal scaling, we choose the prefactors on the right-hand side	such that they enjoy the same scaling in $\epsilon$, i.e., $\mu_2^p \sim \kappa_1^{-2p} \sim \kappa_2^{-p} \sim \epsilon^{-1} \mu_1^{-1}$.
	As $\mu_2$ is determined by the relation $\mu_2 = M \kappa_1 \mu_1$, we find that
	$(\kappa_1 \mu_1)^p \sim \kappa_1^{-2p}$ and thus $\kappa_1 \sim \mu_1^{-1/3}$. Using this relation, our goal is to choose $\mu_1$ such that
	$\mu_1^{2p/3} \sim  \epsilon^{-1} \mu_1^{-1} $, which leads to the choice
$\mu_1 \sim \epsilon^{-3/(2p+3)}$. 
By using these relations in \eqref{eq:optimizethefirst}, we obtain a constant $c>0$ such that for $\epsilon>0$ sufficiently small it holds that
		\begin{align*}
			  c  \leq   \epsilon^{-2p/(2p+3)}\Energyinf{3}{p}{\dir} .
	\end{align*}
This shows Theorem~\ref{mainth:3d} in the case \eqref{eq:choiceofmatrices}.
 
\emph{Step 3 (Proof of (ii) in the case \eqref{eq:secondcase}):}
The strategy is the same as in Step 2, however, the choice of the parameters will be a different one.
Recalling \eqref{eq:3D_wells}, the boundary conditions are given by $F = \diag (1/2, 1/2, \alpha)$ for $\alpha \in (0,1)$
and $\tilde \chi_{33}$ attains the values $1-\alpha$ and $\alpha$.
Thus,
	similarly to  \eqref{eq:3dproof1W1} and \eqref{eq:3dproof2W1}, Plancherel's formula, the triangle inequality, \eqref{eq:decaymultiplay}, and Lemma~\ref{lem:lowfrequencies1} (as $3 \in S_1$) imply that there exists a constant $c>0$ such that  
	\begin{align}
		 c \leq  \Vert \tilde \chi_{33} \Vert_{L^2} \leq \Vert  \tilde \chi_{33} - m_{3,\kappa_3, \mu_3}(D) \tilde \chi_{33} \Vert_{L^2} + C \mu_3 ( \Eelinf{3}{p}^{1/p} + \Eelinf{3}{p}^{1/(2p)} ). \label{eq:3dproof2W13}
	\end{align} 
By using the relation $ \tilde \chi_{33} = g(\tilde \chi_{11}) $ for $g(t) = -4t^2 + 1 - \alpha$ and the fact that $V_1 \cap V_2 = \{ 0 \}$,
	we apply Proposition~\ref{prop:conereduction}  and derive for any $\mu_1>0$ and $\kappa_1 \in (0,M^{-1})$ with $\mu_3 = M \kappa_1 \mu_1$ that
	\begin{align}
		& \quad \ \Vert  \tilde \chi_{33} - m_{3,\kappa_3, \mu_3}(D) \tilde \chi_{33} \Vert_{L^2} \notag \\
		& \leq C    \Vert  \tilde \chi_{11} - m_{1,\kappa_1,\mu_1} (D) \tilde \chi_{11} \Vert_{L^2}^{1/2} + C \Vert \tilde \chi_{33} - m_{3,\kappa_3,\mu_1}(D) \tilde \chi_{33} \Vert_{L^2}  . \label{eq:3dproof3W114}
	\end{align}
	 Thus, Corollary~\ref{cor:lowfrequencyelastic}, Lemma~\ref{lem:highfrequencyoutsidecones}, \eqref{eq:3dproof2W13}, \eqref{eq:3dproof3W114}, and the concavity of $\vert \cdot \vert^{1/2}$ yield
	\begin{align}
			  c &\leq C \mu_3 ( \Eelinf{3}{p}^{1/p} + \Eelinf{3}{p}^{1/(2p)} ) \notag \\
			&\quad \ + C \kappa_3^{-2}  \Eelinf{3}{p}^{1/p} + C  \mu_1^{-1/2} ( \Esurf{3}^{1/2} + \Per^{1/2} ) \notag \\
			& \quad \  + C \kappa_1^{-1/2}  \Eelinf{3}{p}^{1/(2p)} + C  \mu_1^{-1/4} ( \Esurf{3}^{1/4} + \Per^{1/4} )  .
	\end{align}
		Similarly to \eqref{eq:optimizethefirst2}, we use Young's inequality with suitable coefficients to obtain a constant $c>$ such that
	\begin{align}\label{eq:optimizethefirst2}
			  c (1 - \mu_1^{-1}  \Per) \leq   C (\mu_3^p+ \mu_3^{2p} + \kappa_3^{-2p} +\kappa_1^{-p} +  \epsilon^{-1} \mu_1^{-1} ) \Energyinf{3}{p}{\dir} .
	\end{align}
	By the choices $\kappa_3 \sim \mu_3^{-1}$, $\mu_3 \sim \kappa_1^{-1/2}$, $ \kappa_1 \sim \mu_1^{-2/3} $, and $\mu_1 \sim \epsilon^{-3/(2p+3)}$
we conclude the proof of Theorem~\ref{mainth:3d} in the case  \eqref{eq:secondcase}.
\end{proof}

\begin{remark}[Fourier transformation for Dirichlet boundary conditions]\label{rem:Fourierseriesvstrans}
	If one prefers using the (continuous) Fourier transformation, the phase indicator $\chi\colon \Omega \to K_3$ needs to be extended suitably
outside of $\Omega$ to ensure $L^2(\R^d)$-integrability. An extension
by zero would enforce the additional constraint $g(0) = 0$ in the proof above, which was discussed for instance in \cite[equation~(25)]{RT23a} for a nucleation problem.
Our choice of wells is not compatible with the latter constraint and this extension thus reduces the class of wells one can consider by a naive adaptation of the presented technique.
The Fourier \emph{series} instead allows for a periodic extension of the phase indicator.
The trade-off is, however, that one has to choose domains on which the displacement can be extended periodically without paying elastic energy.
As we also prove lower bounds in periodic settings, see Section~\ref{sec:periodicboundary}, the domain $\Omega = (0,1)^d$ hence serves as a convenient choice.
\end{remark}

\subsection{The four-well problem}
\label{sec:fourwell}

As a final motivation on the structure of the general lower scaling bound argument, we discuss the four-well problem from Theorem \ref{mainth:4d}. In particular, as a guideline for our quantitative argument, we present the proof of an analogous qualitative statement (Proposition \ref{prop:stress-free}) first.

\subsubsection{The stress-free case}

We begin by discussing the stress-free case. This serves as an instructive motivation for the structure of the desired result and can also be viewed as giving central hints on the necessary form of the coercivity bounds in the elastic energy which were deduced in Lemma \ref{lem:elasticenergybound}. Compared to the case with stress, we will deduce a slightly more precise result in the stress-free setting. 

\begin{proposition}[The stress-free case]
\label{prop:stress-free}
	Let $u \in W^{1, \infty}(\Omega, \mathbb{R}^{4})$ be a solution to $\sym (\nabla u) \in K_3 $, where $ K_3 $ be determined by the matrices in \eqref{eq:Kin4d}.
	Then, the following holds:
	\begin{itemize}
		\item $\nabla u=A_2$, $\nabla u = A_3 $ or $\nabla u \in\left\{A_0, A_1\right\}$ a.e.\ in $\Omega$.
		\item If $\nabla u \in\left\{A_0, A_1\right\}$ a.e.\ in $\Omega$, then $u$ is (locally) a simple laminate.	
	\end{itemize}
\end{proposition}

\begin{proof}
	The proof is divided into two steps. In the first one, we analyze the influence of the (linear) Saint-Venant compatibility conditions whereas the (nonlinear) relation between the wells is exploited in step 2.

	\emph{Step 1:}
	We define $e \defas \sym (\nabla u)$. As $\Omega$ is simply connected and $e$ is a diagonal matrix,
	the Saint-Venant compatibility conditions of strains imply that $e$ satisfies
	\begin{align}\label{eq:compatibilitystrain}
		\partial_{ij} e_{kk} = 0 \quad \text{for all } i \neq j \neq k \quad \text{and} \quad \partial_{ii} e_{jj} + \partial_{jj} e_{ii} = 0 \quad \text{for all } i \neq j
	\end{align}
	in a distributional sense.
	Using the first condition, there exist functions $b_{ij}$ with $i \neq j$, depending on two variables, such that
	\begin{align*}
		e_{ii}(x) = \sum_{j\neq i} b_{ij} (x_i,x_j) \quad\text{for all } i \in \{1,2,3,4\} . 
	\end{align*}
	As $e_{22} + e_{33} = 0$ and hence $\partial_{k} (e_{22} + e_{33}) = 0$ for $k = 1,4$, we have
	$ \partial_k (b_{2k} (x) + b_{3k}(x) ) = 0$ for $k = 1,4$. As this identity holds for any $x_2$ and $x_3$,
	the derivatives must be functions of $x_1$ and $x_4$, respectively.
	This implies that there exist functions $\tilde b_{2k}$, $\tilde b_{3k}$ and $\hat b_{2k}$, $\hat b_{3k}$ such that
	  $b_{2k}(x) = \tilde b_{2k}(x_k) + \hat b_{2k}(x_2)$ and $b_{3k}(x) = \tilde b_{3k}(x_k) + \hat b_{3k}(x_3)$ for $k = 1,4$.
	By replacing $b_{23}$ and $b_{32}$ with $b_{23} + \hat b_{21} + \hat b_{24} $ and $b_{32} + \hat b_{31} + \hat b_{34}$, respectively,
	we can assume that $b_{ik}(x)=b_{ik}(x_k)$ for $i = 2,3$ and $k = 1,4$.

	By taking variations in $x_2$ and $x_3$, $e_{22} + e_{33} = 0$ yields $b_{23} + b_{32} = const$. Thus,
	the second condition in \eqref{eq:compatibilitystrain} for $i = 2$ and $j = 3$ reveals that
	$\partial_{22} b_{32} - \partial_{33} b_{32} = 0$ distributionally.
	 The formula of d'Alembert provides the existence of one-dimensional functions $\tilde F_{32}$ and $\hat F_{32}$ such that 
	we have
	\begin{align*}
		e_{22}(x) = b_{21}(x_1) + \tilde F_{32}(x_2+x_3) + \hat F_{32} (x_2 - x_3) + b_{24}(x_4), \\
		e_{33}(x) = b_{31}(x_1) - \tilde F_{32}(x_2+x_3) - \hat F_{32} (x_2 - x_3) + b_{34}(x_4).
	\end{align*}
	Notice that $e_{22} + e_{33} = 0$ also implies that $b_{2k} = - b_{3k} + const $ for $k = 1,4$.

Moreover, the second condition in \eqref{eq:compatibilitystrain} for $j = 2,3$ and $i = 1,4$ implies that
\begin{align*}
	\partial_{ii} b_{ji} (x_i) = - \partial_{jj}  b_{ij} (x_i,x_j),
\end{align*}
implying the representation
\begin{align*}
	b_{ij} (x_i,x_j) = c_{ij} (x_i) p_{ij}^{(2)}( x_j) + \tilde c_{ij}(x_i )p_{ij}^{(1)}( x_j) + \hat c_{ij}(x_i) \quad \text{for }j = 2,3 \text{ and }i = 1,4,
\end{align*}
where $p_{ij}^{(k)}( x_j)$ are polynomials of degree $k = 1,2$ and $c_{ij}$, $\tilde c_{ij}$, and $\hat c_{ij}$ are functions only depending on $x_i$.
By taking variations of $x_j$ and the fact that the functions $e_{11}$ and $e_{44}$ only attain discrete values we find that
$p_{ij}^{(k)} =  const $ for $k =1,2$, and we obtain the representations
\begin{align*}
	e_{11} (x) = b_{12} (x_1) + b_{13} (x_1) + b_{14} (x_1,x_4), \qquad
	e_{44} (x) = b_{41} (x_1,x_4) + b_{42} (x_4) + b_{43} (x_4).
\end{align*}

Using the second equation in \eqref{eq:compatibilitystrain} for $i = 1,4$ and $j =2$ once again, we find that
\begin{align*}
	\partial_{ii} b_{2i} (x_i) = - \partial_{22} e_{ii} = 0,
\end{align*}
implying that $b_{2i}(x_i)$ are polynomials of degree $2$. Taking variations of $x_i$ and the fact that $e_{ii}$ only attain
discrete values imply that $b_{2i}(x_i)$ need to be constants, and $e_{22}$ and $e_{33}$ only depend on $x_2$ and $x_3$.
Up to a redefinition of the functions above, we collect the following representations  
\begin{align}\label{eq:dependencies}
	e_{11} (x) &=  b_{14} (x_1,x_4), &&
		e_{22}(x) =  \tilde F_{32}(x_2+x_3) + \hat F_{32} (x_2 - x_3) \notag \\
				e_{33}(x) &=  - \tilde F_{32}(x_2+x_3) - \hat F_{32} (x_2 - x_3)  , && e_{44} (x) = b_{41} (x_1,x_4) . 
\end{align}

\emph{Step 2:}
	Using the definition of $A_i$ in \eqref{eq:Kin4d}, one can check that
	\begin{align}	 
	e_{22} = 4 (e_{11} + \frac{1}{2 \sqrt{2}} e_{44}) (1 - (e_{11} + \frac{1}{2 \sqrt{2}} e_{44}))
 \quad \text{ and } \quad e_{44} = -4 (e_{22}-1) e_{22}  .\label{eq:polynomialdependency}
	\end{align}
The second equation and the fact that $e_{44}$ depends only on $x_1$ and $x_2$ implies that
$e_{44}$ needs to be constant. Moreover, the first polynomial relation in \eqref{eq:polynomialdependency} shows that $\tilde F_{32}$ and $\hat F_{32}$ need to be constant functions, and hence $e_{22}$ is a constant.
Thus, if $e_{44} = 1$, we have $e = A_3$ a.e.\ in $\Omega$.
Moreover, if $e_{44} = 0 $ and $e_{22} = 1$, we have $e = A_2$ a.e.\ in $\Omega$.
Eventually, if $e_{44} = 0 $ and $e_{22} = 0$, we have $e \in \{ A_0, A_1 \}$.
 In this case, the first relation in \eqref{eq:polynomialdependency} reveals that $ e_{11}(1-e_{11}) = 0$, and the second condition in \eqref{eq:compatibilitystrain}
 for $i = 1$ and $j = 4$ yields $\partial_{44} e_{11} = 0$. Proceeding as before, we find that $e_{11}$ only depends on $x_1$, implying that $u$ is locally a simple laminate.
\end{proof}

As a preparation for the quantitative argument, let us draw some conclusions from the proof of Proposition \ref{prop:stress-free}. 
The first step represents the qualitative analogoue of Lemma~\ref{lem:Fourierchar} and Lemma~\ref{lem:elasticenergybound} which transform
the Saint-Venant compatibility conditions and the compatibility directions in the symmetrized lamination convex hull into Fourier variables.
In particular, the dependencies in \eqref{eq:dependencies} are reflected in the multipliers in \eqref{eq:conclusionFourierchar}.
Then, in Step 2, the relation between the wells is taken into account.
Here, we \emph{strongly} benefit from the presence of a non-degenerate compatibility direction between degenerate directions:
The different dependencies in \eqref{eq:dependencies} and the first polynomial relation in \eqref{eq:polynomialdependency} allow us to exclude dependencies of the strain on certain variables.
Contrary, in the case of consecutive degenerate compatibility directions
such a polynomial relation still holds, where both the left- and right-hand side of the relation depend on the \emph{same} coordinates. The proof therefore becomes more difficult which
is also reflected in the more complicated lamination convex hull in Proposition~\ref{prop:topoprop}.

\subsubsection{The quantitative case}\label{sec:fourwellquant}

With the qualitative argument for the four-well problem in hand, we now turn to its quantification. Here we follow the structure of the qualitative argument on a quantitative level.

\begin{proof}[Proof of Theorem~\ref{mainth:4d}]       
Recall that the boundary conditions are prescribed by the matrix $F = \diag(
	1/2, 1/2, -1/2, \alpha )$ for $\alpha \in (0,1)$.
	As $\chi_{44}$ only attains the values $0$ and $1$, it holds that
	\begin{align}
		\Vert \tilde \chi_{44} \Vert_{L^2} = \Vert \chi_{44} - \alpha \Vert_{L^2} \geq \min \{ \vert 1-\alpha \vert, \vert \alpha \vert \} > 0. \label{eq:4dproof1W1}
	\end{align}
	Using Plancherel's formula and Lemma~\ref{lem:lowfrequencies1}, we get for $\kappa_4\in (0,1)$ and $\mu_4>1$  that  
	\begin{align}
		\Vert \tilde \chi_{44} \Vert_{L^2} &\leq \Vert  \tilde \chi_{44} - m_{4,\kappa_4, \mu_4}(D) \tilde \chi_{44} \Vert_{L^2} +  C \mu_4 ( \Eelinf{4}{p}^{1/p} + \Eelinf{4}{p}^{1/(2p)} ) . \label{eq:4dproof2W1} 
	\end{align} 
	 We now employ the cone reduction argument and observe the polynomial relations
	 \begin{align}
	 \label{eq:poly_rel_quant}
		\tilde \chi_{22} = g( (\tilde \chi_{11} + F_{11}) + \frac{1}{2 \sqrt{2}} (\tilde \chi_{44} + F_{44})) - F_{22} \quad \text{ and } \quad \tilde \chi_{44} = g(\tilde \chi_{22} + F_{22}) - F_{44} 
	 \end{align}
	 for $g(t) = 4 t (1-t)$.
Using the second polynomial relation from \eqref{eq:poly_rel_quant}, we apply Proposition~\ref{prop:conereduction} ($ V_2 \cap V_3 = \{ 0 \}$) for the parameters
	 $\mu_2 = \frac{\mu_4}{M\kappa_2} > \mu_4 $, $\gamma \geq \mu_4$, and find that
	\begin{align*}
		& \quad \ \Vert  \tilde \chi_{44} - m_{4,\kappa_4, \mu_4}(D) \tilde \chi_{44} \Vert_{L^2} \\
		& \leq C   \Vert  \tilde \chi_{22} - m_{2,\kappa_2,\mu_2} (D) \tilde \chi_{22} \Vert_{L^2}^{1/2} + C \Vert \tilde \chi_{44} - m_{4,\kappa_4,\gamma}(D) \tilde \chi_{44} \Vert_{L^2}.
	\end{align*}
	A combination of the first polynomial relation from \eqref{eq:poly_rel_quant} with Proposition~\ref{prop:conereduction} (for $\gamma \geq \max \{ \mu_2, \mu_4 \}$ and multipliers $m_{4,\kappa_1, \mu_1}(D)$ and $ m_{1,\kappa_1, \mu_1}(D)$ such that $\mu_2 =M \kappa_1 \mu_1 $ as $ V_2 \cap V_3 = \{ 0 \} = V_1 \cap V_2$) and the concavity of $\vert \cdot \vert^{1/2}$ leads to 
	\begin{align}
		& \quad \ \Vert  \tilde \chi_{44} - m_{4,\kappa_4, \mu_4}(D) \tilde \chi_{44} \Vert_{L^2} \notag \\
		& \leq   C  \Vert \tilde \chi_{11} - m_{1,\kappa_1,\mu_1}(D) \tilde \chi_{11} \Vert_{L^2}^{1/4} + C  \Vert \tilde \chi_{44} - m_{4,\kappa_1, \mu_1}(D) \tilde \chi_{44} \Vert_{L^2}^{1/4} \notag \\
		& \quad  + C \Vert  \tilde \chi_{22} - m_{2,\kappa_2,\gamma} (D) \tilde \chi_{22} \Vert_{L^2}^{1/2} + C \Vert \tilde \chi_{44} - m_{4,\kappa_4,\gamma}(D) \tilde \chi_{44} \Vert_{L^2}. \label{eq:4dproof3W1}
	\end{align}
	 Using \eqref{eq:4dproof1W1}--\eqref{eq:4dproof3W1} and Young's inequality with a small constant
	  to get rid of the roots in \eqref{eq:4dproof3W1} and \eqref{eq:4dproof2W1},
there exists a constant $c>0$ such that
\begin{align*}
	c & \leq  C  (  \mu_4^p + \mu_4^{2p}  ) \Eelinf{4}{p}   + C  \Vert \tilde \chi_{11} - m_{1,\kappa_1,\mu_1}(D) \tilde \chi_{11} \Vert_{L^2} + C  \Vert \tilde \chi_{44} - m_{4,\kappa_1, \mu_1}(D) \tilde \chi_{44} \Vert_{L^2} \notag \\
		& \quad  + C \Vert  \tilde \chi_{22} - m_{2,\kappa_2,\gamma} (D) \tilde \chi_{22} \Vert_{L^2} + C \Vert \tilde \chi_{44} - m_{4,\kappa_4,\gamma}(D) \tilde \chi_{44} \Vert_{L^2}.
\end{align*}
	Using Lemma~\ref{lem:highfrequencyoutsidecones} and Corollary~\ref{cor:lowfrequencyelastic} ($\kappa_i \in (0,1)$), we find that  
\begin{align*}
	c & \leq  C  (  \mu_4^p + \mu_4^{2p}  )  \Eelinf{4}{p}   + C (\gamma^{-1/2} + \mu_1^{-1/2}) ( \Esurf{ 4 }^{1/2} + \Per^{1/2} ) \\ & \quad +  C ( \kappa_4^{-2} + \kappa_2^{-1} +\kappa_1^{-2}  )  \Eelinf{ 4 }{p}^{1/p}  .
\end{align*}
This relation holds for any $\chi$ and thus also for the infimum in \eqref{eq:minprob}. As $\mu_4 > 1$, Young's inequality yields 
\begin{align*}
	c -  (\gamma^{-1} +\mu_1^{-1})  \Per \leq  C  \big(    \mu_4^{2p} + \kappa_4^{-2p} + \kappa_2^{-p} + \kappa_1^{-2p} + \epsilon^{-1}( \gamma^{-1} + \mu_1^{-1}) \big)    \Energyinf{4}{p}{\dir} .
\end{align*}
To derive the optimal scaling, the prefactors of the right-hand side should enjoy the same scaling in $\epsilon$, i.e.,
$\mu_4 \sim \kappa_4^{-1} \sim \kappa_2^{-1/2} \sim \kappa_1^{-1}$.
Hence, the relations $\mu_4 \sim \kappa_2 \mu_2$ and $\mu_2 \sim \kappa_1 \mu_1$ give $
\mu_4^4 \sim \mu_1$, which combined with
 the choice $\gamma = \mu_1$ leads to
\begin{align*}
	c -   \mu_1^{-1} \Per \leq  C  \big(    \mu_1^{p/2} +\epsilon^{-1} \mu_1^{-1} \big)    \Energyinf{4}{p}{\dir} .
\end{align*}
 By choosing $\mu_1 \sim \epsilon^{-2/(p+2)}$, the perimeter can be absorbed by the left-hand side. This concludes the proof. 
\end{proof}

\subsection{The $m+1$-well problem}
\label{sec:mwell}
This section lifts the technique from Subsection~\ref{sec:fourwellquant} to the general setting, aiming to prove Theorems~\ref{mainth:md} and \ref{mainth:mdperiodic}. 
We recall the notation for the effective index $ l(i) = 1 + \sum_{j=1}^{i-1} f(j)$ for $i \in \{1,\dots,m\}$. 
In what follows we will make use of the following nonlinear relations:

\begin{lemma}
\label{lem:poly}
For all $i \in (S_1 \cup S_2) \setminus \{1\} $ there exist polynomials $Q$ and $g$ such that it holds that 
\begin{align*}
\tilde \chi_{l(i) l(i)} = Q(\sqrt{5} \tilde \chi_{l(i-1) l(i-1)} + \sqrt{3}\tilde \chi_{l(i+1)l(i+1)}) = g ( \tilde \chi_{l(i-1)l(i-1) } ) - \sum_{k = i+1}^m 2^{- k - i } \tilde  \chi_{l(k)l(k)},
\end{align*}
where we set $\tilde \chi_{l(m+1)l(m+1)} = 0$. Here the polynomials $Q$ and $g$ depend on $F$. 
\end{lemma}

Let us comment on the role of these polynomial relations. As in the articles \cite{RT22,RT23a,RTTZ25} we use the polynomial relations from Lemma \ref{lem:poly} to relate the components of the phase indicator and, hence, to iteratively decrease the possible areas of phase-space concentration. Contrary to the earlier articles \cite{RT22,RT23a} without gauge invariances and also contrary to the setting of the geometrically linearized cubic-to-tetragonal phase transformation from \cite{RT24}, in our setting with gauges, we are confronted with \emph{two} types of rank-one connections -- the degenerate and the non-degenerate ones. As it will turn out that they require \emph{different} conical phase-space localizations (i.e., localizations in cones with \emph{different} opening angles), we will have to use both nonlinear identities from Lemma \ref{lem:poly}.

Indeed, as discussed in the next sections, the first polynomial relation given by the polynomial $Q$, is particularly suited to entries $i \in S_1$, as then by construction $i-1, i+1 \in S_2$. In particular, this ensures that the cones which enter in the application of the cone reduction result (Proposition \ref{prop:conereduction}) are disjoint. The second relation, given by the polynomial $g$, in turn is particularly suited for indices $i\in S_2$. Although the index $i-1$ can then either be an element of $S_1$ or $S_2$, it still allows us to work with disjoint cones in the nonlinear relation. Similarly, also the linear contributions have disjoint cones of support with respect to the cone associated with the index $i$. Combined with Proposition  \ref{prop:conereduction} this will then eventually allow us to deduce the desired lower scaling bounds (see the proof of Theorem \ref{mainth:md} in the next section).

\begin{proof}
 We first note that it suffices to prove the relations for $\chi$ instead of $\tilde \chi$ due to a linear transformation of the polynomial.

The existence of the claimed polynomials follows from the discrete structure of our wells and their explicit construction. 
According to \eqref{eq:constructiondim}, the values of the matrices $A_j$ develop as follows: 
\begin{center}
\begin{tabular}{|c||c| c| c| c| c| c|}
\hline
matrix & $l(i-1)$ & $l(i)$ & $l(i+1)$ &$l(i+2)$ & $\dots $ & $ l(m)$\\ 
\hline
\hline
$A_{i-1}$ &0 & 0 & 0 & 0 & $\dots $ & $ 0$\\
\hline
$A_{i}$ & 1 & 0 & 0& 0& $\dots $ & $ 0$\\
\hline
$A_{i+1}$ & $\frac{1}{2}$ & 1 & 0 & 0 & $\dots $ & $ 0$\\
\hline
$A_{i+2}$ & $\frac{1}{2}$ & $\frac{1}{2}$ & 1& 0& $\dots $ & $ 0$\\
\hline
$A_{i+3}$ & $\frac{1}{2}$ & $\frac{1}{2}$ & $\frac{1}{2}$& 1 & $\dots $ & $ 0$\\
\hline
$\vdots$ & $\vdots$ & $\vdots$ & $\vdots $& $\vdots$ & $\ddots $ & $ 0$\\
\hline
$A_{m}$ & $\frac{1}{2}$ & $\frac{1}{2}$ & $\frac{1}{2}$& $\frac{1}{2}$ &$\frac{1}{2}$ & $ 1$\\
\hline
\end{tabular}
\end{center}
Now, forming the values $\sqrt{5}\chi_{l(i-1) l(i-1)} + \sqrt{3}\chi_{l(i+1)l(i+1)}$ for $i < m$ we hence obtain as a necessary condition for the existence of the polynomial $Q$ the following conditions
\begin{align*}
Q(0) = 0, \ Q(\sqrt{5}) = 0, \ Q(\sqrt{5}/2) = 1,\ Q(\sqrt{5}/2 + \sqrt{3}) = \frac{1}{2}, \ Q(\sqrt{5}/2 + \sqrt{3}/2) = 1/2.
\end{align*}
This can be satisfied by interpolation with a fourth order polynomial.

Similarly, setting $g(t) = 4 t (1-t)$, we obtain the second polynomial identity for $i<m$, as $g$ only attains the two values $0,1$ on the possible entries of the matrices $A_j$. As observed from the above table, these values are then corrected appropriately by the linear term.

Finally, for $i = m$, we set $Q(t) = g(t) = 4 t (1-t)$,
concluding the proof. 
\end{proof}

\subsubsection{Dirichlet boundary conditions}

In this section, we combine the previously deduced auxiliary results into the proof of Theorem \ref{mainth:md}. Compared to the frequency space arguments from previous works as in \cite{RT22, RT23a, RTTZ25} a major new difficulty consists in the presence of higher order laminates which may be given by combinations of degenerate and non-degenerate symmetrized rank-one conditions. In particular, this will force us to work with both of the polynomials from Lemma \ref{lem:poly} in the application of the cone reduction arguments from Proposition \ref{prop:conereduction}. Moreover, these more complicated polynomial relations will not allow us to immediately reduce the size of the cones from one step to the next, but will require a more complex iterative reduction of the relevant scales. In particular, in what follows, contrary to earlier works, we will consider cones with the two different opening angles $\kappa_{S_1}, \kappa_{S_2}>0$ associated with indices $i \in S_1$ and $i\in S_2$, respectively.

\begin{figure}[t] 
\includegraphics[width =0.7 \textwidth]{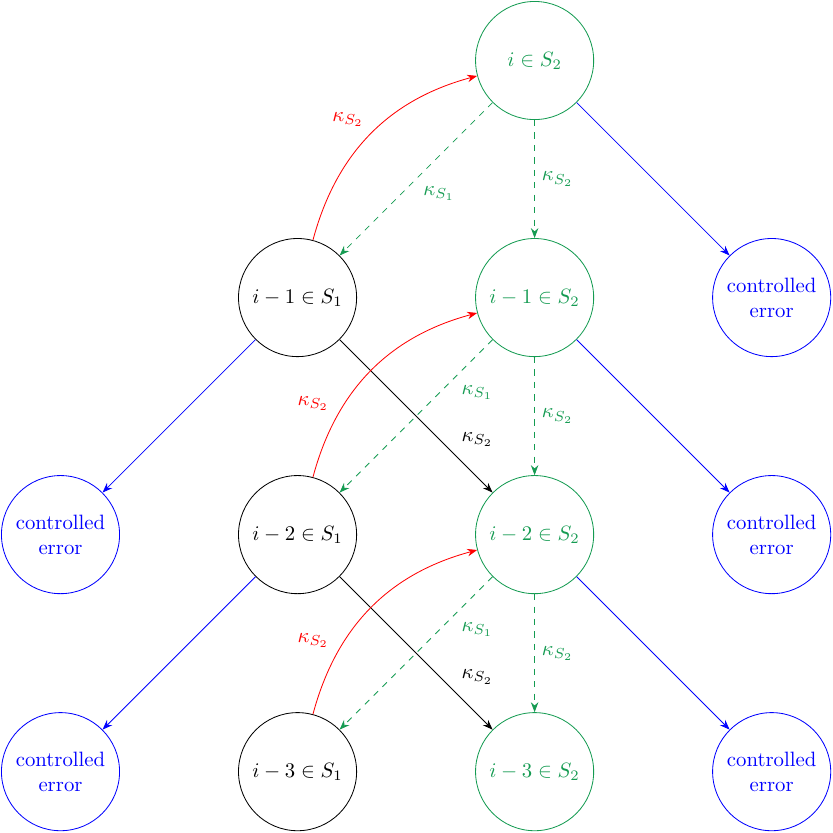}
\caption{A schematic illustration of the iteration argument for the general lower bound. 
The circles represent cones associated with the $i$th index, $i \in \{1 ,\dots , l(\ell) \}$, where
cones corresponding to two compatibility directions (green circles) are controlled via Case 2, and the ones in black via Case 1, i.e.,
black circles indicate cones around degenerate compatibility directions.
More precisely, the Fourier coefficients outside of the cones associated to $S_1$ (green circles) are controlled via
a cone corresponding either to degenerate compatibility directions (black circle with lower index) or again to non-degenerate directions (green circle with lower index), and controlled error terms (blue circles).
The newly arising cones are truncated at higher frequencies, where the new effective truncation parameter is (roughly) the length of the $i$th cone divided by $\kappa_{S_1}$ or $\kappa_{S_2}$, respectively.
Although all paths in the picture can arise, we note that for a given $f \in \mathcal{S}$, the iteration follows only one of the dashed arrows.  
In Case 1 the cones are truncated via non-degenerate ones, where the preceding (black arrow) and \emph{following} index (red arrow) arises.
 The red arrows need a separate discussion in order to ensure that the resulting contributions remain controlled.
}
\label{fig:iteration}
\end{figure} 

\begin{proof}[Proof of Theorem~\ref{mainth:md}] 
	Let $\ell \in \{1,\dots,m\}$, $F \in K^{(\ell)}_m \setminus K^{(\ell-1)}_m$. 
	To derive the optimal scaling, we analyze the variable $\tilde \chi_{l(\ell)l(\ell)}$ which captures the coarsest scale of oscillations. We recall that $l(\ell) = 1 + \sum_{j=1}^{\ell-1} f(j)$ represents the effective index.

	The proof is divided into four steps. A uniform lower bound on the $L^2$-norm of $\tilde \chi_{l(\ell)l(\ell)}$ is provided in Step 1.
	While the boundary conditions are exploited in Step 2, Step 3 addresses the iterative cone reduction argument. 
	The precise scaling behaviour is then determined in Step 4.

\emph{Step 1: A lower bound using the boundary conditions.}
	If $\ell = m$, $\tilde \chi_{l(\ell)l(\ell)}$ only attains the values $1- F_{l(m)l(m)}$ and $-  F_{l(m)l(m)}$ implying that there exists a constant such that
	\begin{align*}
		\Vert \tilde \chi_{l(\ell)l(\ell)} \Vert_{L^2}^2 & \geq c > 0 .
	\end{align*}
If $\ell < m$,
 $\chi_{l(\ell)l(\ell)}$ can attain the values $0$, $1/2$, and $1$ and the argument is more involved due to the potential boundary datum $F_{l(\ell)l(\ell)} = 1/2$. As $u^*_{l(k)}$ attains zero boundary conditions for $k \in \{ \ell +1, \dots, m \}$ for the optimal $u^*$ in \eqref{eq:Dirichletsetandperiodic}, it holds that
 \begin{align}\label{eq:estimatephaseindhigh}
			 \sum_{k = \ell +1}^m \vert \{ \chi_{l(k)l(k)} = 1 \} \vert  \leq  \sum_{k = \ell +1}^m \int_{\Omega} \chi_{l(k)l(k)} - \partial_{l(k)} u^*_{l(k)} \di x   \leq C   \Eelinf{m+1}{p}^{1/p}.
 \end{align}
 Thus, there exists a constant $c>0$ such that
	\begin{align*}
		\Vert  \tilde \chi_{l(\ell)l(\ell)} \Vert_{L^2}^2 & \geq \int_{\{  \chi_{l(\ell)l(\ell)}  \in \{0,1\} \} } \vert  \chi_{l(\ell)l(\ell)}  - F_{l(\ell)l(\ell)} \vert^2 \di x = c (\vert \Omega \vert - \vert \{ \chi_{l(\ell)l(\ell)} = 1/2 \} \vert) \\
				&= c \left(\vert \Omega \vert - \sum_{k = \ell +1}^m \vert \{ \chi_{l(k)l(k)} = 1 \} \vert \right) \geq c \left( 1 -    \Eelinf{m+1}{p}^{1/p} \right) .
	\end{align*}

	\emph{Step 2: Combining Step 1 with the low frequency bound.}
Using Step 1, Plancherel's formula and Lemma~\ref{lem:lowfrequencies1}, we derive for $\mu_{l(\ell)}>1$  that  
	\begin{align} \label{eq:dproof2new} 
		 c  & \leq C \Eelinf{m+1}{p}^{1/(2p)} + \Vert  \tilde \chi_{l(\ell)l(\ell)} - m_{{l(\ell)},\kappa_{l(\ell)}, \mu_{l(\ell)}}(D) \tilde \chi_{l(\ell)l(\ell)} \Vert_{L^2}  \notag   \\
		 & \quad + \| m_{{l(\ell)},\kappa_{l(\ell)}, \mu_{l(\ell)}}(D) \tilde \chi_{l(\ell)l(\ell)}  \|_{L^2} \notag \\
		& \leq C \Eelinf{m+1}{p}^{1/(2p)} + \Vert  \tilde \chi_{l(\ell)l(\ell)} - m_{{l(\ell)},\kappa_{l(\ell)}, \mu_{l(\ell)}}(D) \tilde \chi_{l(\ell)l(\ell)} \Vert_{L^2}  \notag   \\
		& \quad + C \mu_{l(\ell)} \delta_{1 f(\ell)}  ( \Eelinf{m+1}{p}^{1/p} + \Eelinf{m+1}{p}^{1/(2p)} ) + C \mu_{l(\ell)}  \delta_{2 f(\ell)} \Eelinf{m+1}{p}^{1/p}  .  
	\end{align}

\emph{Step 3: Iterative cone reduction using the nonlinear relations between the entries.}
	We will employ Proposition~\ref{prop:conereduction} iteratively (and reduce the cones $\ell-1$ times) in \emph{two} different versions. To this end, let $i \in  \{2,\dots, \ell\}$.

\emph{Case 1.}
If $i \in S_1$, we deduce from Lemma~\ref{lem:poly} that there exists a polynomial $Q$ such that $\tilde \chi_{l(i)l(i)} = Q (\sqrt{5}\tilde \chi_{l(i-1)l(i-1)} + \sqrt{3} \tilde \chi_{l(i+1)l(i+1)}  )$ with
$\tilde \chi_{l(m+1)l(m+1)} = 0$. 
As there does not exist $i \in \{1,\dots, m-1\}$ such that $f(i)= f(i+1) = 1$, we have $V_i \cap V_{i-1} = \{0 \}$ and $V_i \cap V_{i+1} = \{ 0 \}$ by \eqref{def:vectorspaces}, and the assumptions of Proposition~\ref{prop:conereduction} are verified (for $E_1 = \{i-1,i+1\}$, $E_2 = \emptyset$, $E_3 = \{ i \}$ with $\lambda_i = 1$, and $\lambda_j = 0$ for $j \notin E_1 \cup E_3$).
More precisely, for any $i \in S_1$ 
the choices (in the notation of Proposition~\ref{prop:conereduction}) $\gamma = \gamma_{l(i)} \geq \mu_{l(i)} = \alpha_i $, $\kappa_{S_1} = \kappa_{l(i)}$, $\kappa_{S_2} = \kappa_{l(i-1)} = \kappa_{l(i+1)}$
with $ \alpha_i = M \kappa_{S_2}  \alpha_{i-1} \leq \alpha_{i-1} $ 
lead to
\begin{align}\label{ineq:coneredcase1}
& \quad \ \|\tilde{\chi}_{l(i)l(i)} - m_{l(i),\kappa_{S_1}, \alpha_i }(D) \tilde{\chi}_{l(i)l(i)}\|_{L^2} \notag \\
 & \leq C  \|\tilde{\chi}_{l(i)l(i)} - m_{l(i),\kappa_{S_1}, \gamma}(D) \tilde{\chi}_{l(i)l(i)}\|_{L^2} \notag \\
 & \quad + C \| \tilde\chi_{l(i-1)l(i-1)} - m_{l(i-1), \kappa_{S_2} , \alpha_{i-1} }(D) \tilde\chi_{l(i-1)l(i-1)} \|_{L^2}^{1/2} \notag \\
 & \quad + C (1-\delta_{ m i }) \| \tilde\chi_{l(i+1)l(i+1)} - m_{l(i+1), \kappa_{S_2} , \alpha_{i-1} }(D) \tilde\chi_{l(i+1)l(i+1)} \|_{L^2}^{1/2}.
\end{align}

\emph{Case 2.}
If $i \in S_2$, we use the second polynomial relation from Lemma~\ref{lem:poly}.
Due to \eqref{def:vectorspaces}, the assumptions of Proposition~\ref{prop:conereduction} are verified (for $E_1 = \{i-1\}$, $E_2 = \{i+1, \dots,m \}$, and $E_3 = \{ i \}$).
Thus, we derive
for any $i \in S_2$ and
the choices $\gamma = \gamma_{l(i)} \geq \mu_{l(i)} = \alpha_i $, $\kappa_{S_2} = \kappa_{l(i)}$, $\kappa_{ S_{f(i-1)} } = \kappa_{l(i-1)}$
with $ \alpha_i  = M \kappa_{ S_{f(i-1)} }   \alpha_{i-1} \leq \alpha_{i-1} $,
$\gamma = \mu_{l(k)}$ for $r \in \{i+1, \dots, m \}$
that
\begin{align}\label{ineq:coneredcase2}
& \quad \ \|\tilde{\chi}_{l(i)l(i)} - m_{l(i),\kappa_{S_2}, \alpha_i }(D) \tilde{\chi}_{l(i)l(i)}\|_{L^2} \notag \\
 & \leq C  \|\tilde{\chi}_{l(i)l(i)} - m_{l(i),\kappa_{S_2}, \gamma}(D) \tilde{\chi}_{l(i)l(i)}\|_{L^2} \notag \\  
 & \quad + C  \| \tilde\chi_{l(i-1)l(i-1)} - m_{l(i-1), \kappa_{ S_{f(i-1)} } , \alpha_{i-1} }(D) \tilde\chi_{l(i-1)l(i-1)} \|_{L^2}^{1/2} \notag \\
 & \quad + C \sum_{r \in \{i+1,\dots, m\} \cap S_1}   \| \tilde\chi_{l(r)l(r)} - m_{l(r), \kappa_{S_1}, \gamma}(D) \tilde\chi_{l(r)l(r)} \|_{L^2} \notag \\
 & \quad +  C \sum_{r \in \{i+1,\dots, m\} \cap S_2}   \| \tilde\chi_{l(r)l(r)} - m_{l(r), \kappa_{S_2}, \gamma}(D) \tilde\chi_{l(r)l(r)} \|_{L^2} .
\end{align}
 
Our goal is to apply 
\eqref{ineq:coneredcase1} and \eqref{ineq:coneredcase2}
iteratively for $i \in \{2, \dots, \ell \}$ such that the coefficients of $\tilde \chi_{l(\ell)l(\ell)}$ outside of the cone with length $\sim \alpha_\ell $
can be controlled with coefficients at significantly higher frequencies with threshold $\sim \alpha_1$. As $\alpha_1 \gg \alpha_\ell$, this has the advantage that less surface energy needs to be used, leading to a sharper bound.
 Compared to \cite{RTTZ25} this iteration is more intricate due to the different opening angles $\kappa_{S_1}$ and $\kappa_{S_2}$,
which we explain in more detail by also providing a schematic visualization in
Figure~\ref{fig:iteration}:
We note that all terms on the right hand sides of \eqref{ineq:coneredcase1} and \eqref{ineq:coneredcase2}  with cut-off length at scale $\gamma$ are ``good'' terms as we may always choose $\gamma = \alpha_1$. Moreover, it will also turn out in our iteration argument below, that the contributions involving $l(i-1)$ and $\alpha_{i-1}$ are ``good'' contributions in the sense that 
 the ``correct'' length scales emerge.

It thus remains to discuss the last term in \eqref{ineq:coneredcase1} (red arrow in Figure~\ref{fig:iteration}) for 
 $i < \ell$ (the phase indicator $\tilde \chi_{l(\ell+1)l(\ell+1)}$ will be controlled with elastic energy, see \eqref{ineq:specialcase2} below), 
and to verify that also this is a ``good'' contribution in the sense that the length scales obtained in an iteration of this term are controlled by the previous ones.
To this end, we first recall that as in \eqref{ineq:coneredcase1} it holds that $i \in S_1$, we have that $i+1 \in S_2$. Hence, in our next iteration step we first apply \eqref{ineq:coneredcase2}. In this application the only contribution on the right hand side of \eqref{ineq:coneredcase2} which does not already involve $\gamma$  as the cut-off length scale (and is thus a ``good contribution'' in the sense outlined above) 
 then again corresponds to the index $l(i)$ with $i \in S_1$. 
 In a second step, we then again apply \eqref{ineq:coneredcase1}. Combined, this leads to the introduction of a new parameter $\tilde \alpha_{i-1}$ that satisfies the relation 
$\alpha_{i-1} = M^2 \kappa_{S_1} \kappa_{S_2} \tilde \alpha_{i-1}$. More precisely, using the concavity of the root functions, we find that
\begin{align*}
	 & \quad \| \tilde\chi_{l(i+1)l(i+1)} - m_{l(i+1), \kappa_{S_2} , \alpha_{i-1} }(D) \tilde\chi_{l(i+1)l(i+1)} \|_{L^2} \\
 & \leq   C \| \tilde\chi_{l(i-1)l(i-1)} - m_{l(i-1), \kappa_{S_2} , \tilde \alpha_{i-1}}(D) \tilde\chi_{l(i-1)l(i-1)} \|_{L^2}^{1/4} \notag \\
 & \quad + C   \| \tilde\chi_{l(i+1)l(i+1)} - m_{l(i+1), \kappa_{S_2} , \tilde \alpha_{i-1} }(D) \tilde\chi_{l(i+1)l(i+1)} \|_{L^2}^{1/4} \notag \\
 & \quad + C \sum_{r \in \{i+1,\dots, m\} \cap S_1}   \| \tilde\chi_{l(r)l(r)} - m_{l(r), \kappa_{S_1}, \gamma}(D) \tilde\chi_{l(r)l(r)} \|_{L^2} \notag \\
 & \quad +  C \sum_{r \in \{i+1,\dots, m\} \cap S_2}   \| \tilde\chi_{l(r)l(r)} - m_{l(r), \kappa_{S_2}, \gamma}(D) \tilde\chi_{l(r)l(r)} \|_{L^2} \notag \\
 & \quad + C  \|\tilde{\chi}_{l(i)l(i)} - m_{l(i),\kappa_{S_1}, \gamma}(D) \tilde{\chi}_{l(i)l(i)}\|_{L^2}^{1/2} .
\end{align*}
In particular, by iterating this estimate, up to changing $C>0$ and adding fractional powers of type $2^{-j+1}$ for $j \in \{1, 2\ell -1\}$, we can choose $\tilde \alpha_{i-1} = M^{2 (1-\ell)} \kappa_{S_1}^{-\ell+1} \kappa_{S_2}^{-\ell+1}  \alpha_{i-1} \gg \alpha_{i-1} $.  
The iterative application of the estimates \eqref{ineq:coneredcase1} and \eqref{ineq:coneredcase2} for $i \in \{2, \dots, \ell \}$ gives the relation
\begin{align}\label{eq:iterativerelation}
	\alpha_\ell = M^{\ell-1} \kappa_{S_1}^{k_{\ell-1}}  \kappa_{S_2}^{\ell -1 - k_{\ell-1}}  \alpha_1
\end{align}
with $k_{\ell-1}$ defined in \eqref{eq:kr}. Moreover, we have ensured that $\tilde \alpha_{i-1} \gg \alpha_1$,
 which together with the choice $\gamma = \alpha_1$ leads to 
\begin{align}\label{ineq:afterconered}
	& \quad \ \Vert  \tilde \chi_{l(\ell)l(\ell)} - m_{{l(\ell)},\kappa_{l(\ell)}, \alpha_\ell}(D) \tilde \chi_{l(\ell)l(\ell)} \Vert_{L^2} \notag \\
	& \leq C  \sum_{i \in \{1,\dots, m\} \cap S_1}  \sum_{j=1}^{ 2\ell }  \Vert  \tilde \chi_{l(i)l(i)} - m_{{l(i)},\kappa_{S_1}, \alpha_1}(D) \tilde \chi_{l(i)l(i)} \Vert_{L^2}^{2^{-j+1}} \notag \\
	& \quad  + C  \sum_{i \in \{1,\dots, m\} \cap S_2}  \sum_{j=1}^{ 2\ell }   \Vert  \tilde \chi_{l(i)l(i)} - m_{{l(i)},\kappa_{S_2}, \alpha_1}(D) \tilde \chi_{l(i)l(i)} \Vert_{L^2}^{2^{-j+1}},
\end{align}
where we added some additional terms on the right-hand side to achieve a compact notational relation.

\emph{Step 4: Optimization and conclusion.}
For $i \in \{ \ell +1, \dots , m\}$, the triangle inequality,  \eqref{eq:decaymultiplay}, Plancherel's formula, and \eqref{eq:estimatephaseindhigh} yield
\begin{align}\label{ineq:specialcase2}
	& \quad \ \Vert  \tilde \chi_{l(i)l(i)} - m_{{l(i)},\kappa_{S_1}, \alpha_1}(D) \tilde \chi_{l(i)l(i)} \Vert_{L^2}^2 \leq C \Vert \tilde \chi_{l(i)l(i)}  \Vert_{L^2}^2 \leq C \sum_{k = \ell +1}^m \vert \{ \chi_{l(k)l(k)} = 1 \} \vert \notag \\
	& \leq C   \Eelinf{m+1}{p}^{1/p}.
\end{align}
As in the previous proofs we combine
 \eqref{eq:dproof2new} and \eqref{ineq:afterconered}, eliminate the roots in \eqref{ineq:afterconered} via Young's inequality, employ Lemma~\ref{lem:highfrequencyoutsidecones}, Corollary~\ref{cor:lowfrequencyelastic}, and \eqref{ineq:specialcase2}
 such that there exist constants $c, C>0$ satisfying
 	\begin{align*}
		 c - \alpha_1^{-1/2} \Per^{1/2}  & \leq C \Eelinf{m+1}{p}^{1/(2p)} +C \epsilon^{-1/2} \alpha_1^{-1/2}  \Energyinf{m+1}{p}{\dir}^{1/2}  \\
		& \quad +  C (\kappa_{S_1}^{-2} + \kappa_{S_2}^{-1})  \Eelinf{m+1}{p}^{1/p}   \\
		& \quad + C \alpha_\ell \delta_{1 f(\ell)}  ( \Eelinf{m+1}{p}^{1/p} + \Eelinf{m+1}{p}^{1/(2p)} )   \\
		&\quad + C \alpha_\ell \delta_{2 f(\ell)} \Eelinf{m+1}{p}^{1/p} .
	\end{align*} 
As this inequality holds for any $\tilde \chi$, we apply Young's inequality once again, use that $\alpha_\ell> 1$ and derive that
 	\begin{align*}
		 c - \alpha_1^{-1/2} \Per^{1/2}  & \leq  C (\alpha_\ell^{2p} \delta_{1 f(\ell)} + \alpha_\ell^p  \delta_{2 f(\ell)} +\kappa_{S_1}^{-2p} + \kappa_{S_2}^{-p} +\epsilon^{-1} \alpha_1^{-1})  \Energyinf{m+1}{p}{\dir} .
	\end{align*} 
We choose every parameter on the right-hand side to have the same scaling in $\epsilon$. In view of
\eqref{eq:iterativerelation} and  the relation $k_\ell = k_{\ell -1} + \delta_{1f(\ell)}$, we thus get
$\alpha_1 \sim \alpha_\ell \kappa_{S_1}^{-2({\ell -1 - k_{\ell-1}})-{k_{\ell-1}}} \sim 
  \kappa_{S_1}^{-2\ell + k_{\ell} } $.
This choice leads to the inequality 
  	\begin{align*}
		 c - \alpha_1^{-1/2} \Per^{1/2}  & \leq  C  (\alpha_1^{-2p / (-2\ell  + k_{\ell})} +\epsilon^{-1} \alpha_1^{-1})  \Energyinf{m+1}{p}{\dir} .
	\end{align*} 
The choice $\alpha_1 \sim \epsilon^{ (- 2\ell + k_\ell ) / (2p +2 \ell - k_\ell) }$ results in the desired scaling and concludes the proof.
\end{proof}

We remark that the above argument shows that due to the presence of the two types of symmetrized rank-one directions -- in contrast to the setting without gauges -- we here work with two types of cones -- one with opening angle $\kappa_{S_1}$ and one with $\kappa_{S_2}$. Our above arguments show that $0<\kappa_{S_2}<\kappa_{S_1}<1$.

\subsubsection{Periodic boundary conditions}

Finally, as a last result on lower bounds, we now turn to the proof of Theorem~\ref{mainth:mdperiodic}. The main difference compared to the Dirichlet setting here consists in the treatment of the low frequencies where instead of using the Dirichlet datum, we rely on the discreteness of the periodic Fourier transform. Technically, this gives rise to the shift in the scaling behaviour, as expected from upper bound considerations which do not require the outermost branching level.

\begin{proof}[Proof of Theorem~\ref{mainth:mdperiodic}] 
	Let $\ell \in \{1,\dots,m\}$, $F \in K_m^{(\ell)} \setminus K_m^{(\ell-1)} $.
The proof follows along the lines of Theorem~\ref{mainth:md} with the essential difference being that we cannot use Dirichlet boundary conditions for the control of low frequencies.
Due to the similarities, we only refer to the relevant changes.

In particular, we choose some \emph{fixed} $\mu_{l(\ell)}<1$, use \eqref{ineq:foursym} for the zero frequency and replace the estimate \eqref{eq:dproof2new} by 
	\begin{align} \label{eq:dproof2newperi} 
		 0 &< c     \leq C \Eelinf{m+1}{p}^{1/(2p)} + \Vert   m_{{l(\ell)},\kappa_{l(\ell)}, \mu_{l(\ell)}}(D) \tilde \chi_{l(\ell)l(\ell)} \Vert_{L^2} \notag \\
		 & \qquad \qquad \qquad \qquad \qquad \qquad + \Vert  \tilde \chi_{l(\ell)l(\ell)} - m_{{l(\ell)},\kappa_{l(\ell)}, \mu_{l(\ell)}}(D) \tilde \chi_{l(\ell)l(\ell)} \Vert_{L^2}    \\
		 		    & \leq C (\Eelinf{m+1}{p}^{1/(2p)} + \Eelinf{m+1}{p}^{1/p}  )  + \Vert  \tilde \chi_{l(\ell)l(\ell)} - m_{{l(\ell)},\kappa_{l(\ell)}, \mu_{l(\ell)}}(D) \tilde \chi_{l(\ell)l(\ell)} \Vert_{L^2} . \notag  
	\end{align} 
The iterative application of Proposition~\ref{prop:conereduction} leads to the same relation as in \eqref{eq:iterativerelation}
\begin{align}\label{eq:iterativerelationper}
	\mu_{l(\ell)} = M^{\ell-1} \kappa_{S_1}^{k_{\ell-1}}  \kappa_{S_2}^{\ell -1 - k_{\ell-1}}  \alpha_1
\end{align}
 and to the similar optimization problem 
 	\begin{align*}
		  c - \alpha_1^{-1/2} \Per^{1/2}    & \leq  C (\kappa_{S_1}^{-2p} + \kappa_{S_2}^{-p} +\epsilon^{-1} \alpha_1^{-1})  \Energyinf{m+1}{p}{\per} .
	\end{align*} 
We choose every parameter on the right-hand side to have the same scaling in $\epsilon$.
As the parameter $\mu_{l(\ell)}<1$ is fixed, \eqref{eq:iterativerelationper} implies that
$ \kappa_{S_1} \sim \alpha_1^{-1/ (2\ell -2 - k_{\ell-1})}   $, leading to
  	\begin{align*}
		  c - \alpha_1^{-1/2} \Per^{1/2}  & \leq  C  (\alpha_1^{2p / (2\ell -2 - k_{\ell-1})}+\epsilon^{-1} \alpha_1^{-1})  \Energyinf{m+1}{p}{\dir} .
	\end{align*} 
The choice $\alpha_1 \sim \epsilon^{ (- 2\ell + k_{\ell-1}+2 ) / (2p +2 \ell - k_{\ell-1} - 2) }$ results in the desired scaling and concludes the proof.
\end{proof}

\section{Upper bounds}
\label{sec:upper}  
In this section, we provide the proof of the upper scaling bounds in the model setting of \eqref{eq:choiceofmatrices}, see Proposition~\ref{prop:Uppercylinder} and Theorem~\ref{thm:mainresultupper}. To this end, this section is divided into three parts. In the first subsection, we prove Proposition~\ref{prop:Uppercylinder} and state 
auxiliary lemmata that  will help us to provide a construction in three dimensions for second order Dirichlet data. The latter is achieved in the second subsection.
 Here, we focus on a suitably rotated cube such that the compatibility directions appearing in the lamination convex hull match the basis vectors of the rotated cube.
The advantage of such a rotated domain is that it allows us to reduce the problem to two-dimensional branching constructions,
which have been presented explicitly in \cite[Section~2]{CC15} and \cite[Section~3.1]{RT23}, by means of a strategy presented in \cite[Section~6.2]{RT23}.
Eventually, in the last subsection, we discuss optimal upper bound constructions for higher order data \emph{heuristically}, 
 giving evidence of the optimality of the lower bounds in Theorem~\ref{mainth:md} and Theorem~\ref{mainth:mdperiodic}.

\subsection{Auxiliary constructions}
We begin by recalling the change of the elastic energy under coordinate transformations, including their relation to the problem without the gauge group ${\rm Skew}(d)$. Here and in what follows, it will be convenient to track also the domain dependence of the energy contributions. Relying on the equivalence result from \cite[Theorem 3, Proposition 4.3 and Lemma 4.4]{RTTZ25}, we further drop the phase indicator in what follows below and work with surface energies which directly measure the total variation of the deformation. 

To this end, given some set $\omega \subset \R^d$, a deformation $v \in W^{1,p}_{\loc}(\R^d, \R^d)$, $p \in [1,\infty)$ and some finite set $K \subset \R^{d \times d}_{\rm sym}$, 
 we define the elastic, surface, and singularly perturbed energy on \emph{$\omega$} by

			\begin{align*}
				E_{\rm el}^{(p)}(v,\omega) & \defas \int_{\omega} \dist^p( \sym (\nabla v), K) \di x,  && E_{\rm surf}(\omega) \defas   \Vert D^2 v \Vert_{\rm TV (\omega)}, \\
				 E_\epsilon^{(p)}(v;\omega) & \defas E_{\rm el}^{(p)}(v;\omega) + \epsilon E_{\rm surf}(\omega),  && E_{\rm el}^{(p)}(v;\omega)  \defas E_{\rm el}^{(p)}(v;\omega).
			\end{align*}

		Moreover, if we seek to emphasize the dependence on the prescribed boundary data $F \in \R^{d\times d}$ and domain $\omega \subset \R^d$, we also write $E_{\epsilon,F}^{(p)}(v; \omega;K)$. 
With this notation, we record the following transformation behaviour.

\begin{lemma}[Transformation of the problem]\label{lem:coordinatetransform}
	\begin{itemize}
		\item[(i)] Let $K \subset \R^{d \times d}_{\sym}$ be a set of wells, let $\omega \subset \R^d$, $F \in \R^{d\times d}$ and let $u\colon\omega \rightarrow \R^d$. Let $S \in GL_+(d)$ and set $\hat K = S K S^T$, $\hat F = S F S^T$, $\hat{\omega}\defas S \omega$ and $\hat{u} \defas S u(S^T \cdot)\colon \hat{\omega} \rightarrow \R^d$.  Let $s\colon(0,\infty) \rightarrow [0,\infty]$.
		Then the following equivalences hold
		\begin{align*}
		E_{\epsilon,F}^{(p)}(u; \omega;K) \lesssim s(\epsilon) \Leftrightarrow E_{\epsilon,\hat{F}}^{(p)}(\hat{u}; \hat{\omega}; \hat{K}) \lesssim s(\epsilon),\\
		E_{\epsilon,F}^{(p)}(u; \omega;K) \gtrsim s(\epsilon) \Leftrightarrow E_{\epsilon,\hat{F}}^{(p)}(\hat{u}; \hat{\omega}; \hat{K}) \gtrsim s(\epsilon).	
		\end{align*}
		The constants in the equivalences only depend on $F, F^{-1}$.
		\item[(ii)] Let $K = \{ A,B \} $ with $A = a \odot b$ for $a,b \in \R^d$ and $B = 0$. 		
		Given an upper bound for the two-well problem without the gauge group ${\rm Skew}(d)$ with the wells $ \tilde K = \{ a \otimes b, 0 \} $ or $ \tilde K = \{ b \otimes a, 0 \} $, then the upper bound (up to a skew-symmetric linear transformation) is also admissible for the symmetrized two-well problem with wells $K$, in particular, 
\begin{align*}
\Energyinf{2}{p}{\dir} \lesssim  \int\limits_{\Omega} \dist^p(\nabla u, \tilde{K}) dx + \epsilon \| D^2 u \|_{TV(\Omega)}  ,
\end{align*}		
where $F = \lambda A$ and $\tilde{F} = \lambda a \otimes b$ or $\tilde{F} = \lambda b\otimes a$ for $\lambda \in (0,1)$.
		
		\end{itemize}
\end{lemma}
The proof is rather elementary and follows from a suitable change of variables.
\begin{proof}[Proof of Lemma~\ref{lem:coordinatetransform}]
 Given a function $  u \colon \omega \to \R^d$ on a set $\omega \subset \R^d$, we define its transformation by
$ \hat u \colon S\omega \to \R^d$ with $ \hat u (y) = S  u (S^T y)$. Using the chain rule, we find that
\begin{align*}
	\dist (\sym (\nabla \hat u (y)) , \hat K) = \dist \big( S ( \sym (\nabla  u (S^Ty) )  ) S^T , S K S^T \big) \sim\dist (  \sym (\nabla  u (S^T y) )  ,  K ) ,
\end{align*}
where the inequality of the last term depends on the norm of $S$ and $S^{-1}$. By the definition of $ \hat u $ the boundary datum changes to $\hat F$.
Similarly, this linear transformation implies that $\Vert D^2  u \Vert_{TV(\omega)} \sim  
  \Vert D^2 \hat u \Vert_{TV(S\omega)} $, which yields (i).

Consider the boundary datum $ F = \lambda a \odot b $ for $\lambda \in (0,1)$. Consider the upper bound 
construction $v \colon \omega \to \R^d$ for the gradient case with boundary datum $ \lambda a \otimes b$ and wells $\tilde K = \{ a \otimes b, 0 \}$.
Setting
 $u(x) \defas v(x) - \lambda \skw (a \otimes b) x$, the correct boundary datum is attained, and we have
\begin{align*}
	\dist ( \sym (\nabla u) , \{ A \odot B,0 \} ) \leq \dist (\nabla  v   ,  \{ A \otimes B ,0 \}) .
\end{align*}
As $\Vert D^2 v \Vert_{TV(\omega)} = 
  \Vert D^2 u \Vert_{TV(S\omega)} $, we conclude (ii).  
\end{proof}

\begin{remark}[Upper bounds satisfying an $\epsilon^{2/3}$-scaling]
	In the spirit of Lemma~\ref{lem:coordinatetransform} one can show
	that upper bound constructions for the symmetrized two-well problem that guarantee an $\epsilon^{2/3}$-scaling can be deduced from constructions for the two-well problem without the gauge group ${\rm Skew}(d)$, up to a suitable linear transformation of the domain and the boundary datum. 
	Indeed, given matrices $A$ and $B$ which are symmetrized rank-one connected, we have to find a suitable linear transformation $C \in GL_+(d)$ and a function $v\colon C\omega \to \R^d$ such that $E_\epsilon^{(2)}(v;C\omega)\lesssim \epsilon^{2/3}$.
As $A$ and $B$ are symmetrized rank-one connected, there exist $a,b \in \R^d$ such that $A-B = a \odot b$. 
By the linear transformation $Ax$, we can thus assume that $A = a \odot b$ and $B = 0$. Due to Lemma~\ref{lem:coordinatetransform}(ii), it suffices to find
an upper bound for the two-well problem without the gauge group ${\rm Skew}(d)$ with wells $\tilde K = \{ a \otimes b, 0 \} $.
Given an upper bound $v\colon \omega \to \R^d$ for the two-well problem with wells $\{ c \otimes d , 0 \}$,
we find matrices $C, D \in GL_+(d)$ with $Ca = c$ and $Db = d$ and set $u\colon D^T \omega \to \R^d$ as
$u(x) = C^{-1} v( D^{-T} x)$. Using the relation $C^{-1} (c \otimes d) D^{-T}  =  a \otimes b$, we conclude the statement as in the proof of Lemma~\ref{lem:coordinatetransform}(i).
\end{remark}

The following lemma collects existing constructions in two dimensions, which we will rely on in the three-dimensional setting. As observed by \cite{CC14}, 
we distinguish between two cases as the scaling for the two-well problem depends on the number of compatibility directions.
\begin{lemma}[First order constructions in two dimensions]\label{lem:firstorder2d}
Let $p \in [1,+\infty)$, $0<L, H\le1$, $R=[0,L]\times[0,H]$, and $N\in\N$ such that $N>\frac{4L}{H}$. 
\begin{itemize}
	\item[(i)] \emph{Two (non-degenerate) symmetrized rank-one directions:}
Let $$A = \begin{pmatrix}
0 & -1 \\
-1 & 0 \\
\end{pmatrix} \quad \text{and} \quad B = \begin{pmatrix}
0 & 1 \\
1 & 0 \\
\end{pmatrix}.$$  
Then there exists $u \in W^{1,\infty}_0( R;\R^2)$  (depending on $\epsilon$) with $\nabla u\in BV(R;\R^{2\times 2})$  such that
\begin{align}\label{eq:eltwodirectiontwowell}
	E_{\rm el}^{(p)}(u;R) \leq \int_R \dist^p \left(\nabla u ,  \left\{ \begin{pmatrix}
	0 & 0 \\
	-2 & 0 \\
	\end{pmatrix},  \begin{pmatrix}
	0 & 0 \\
	2 & 0 \\
	\end{pmatrix} \right\}  \right) \di x \lesssim   \frac{L^{p+1}}{N^p H^{p-1}}
\end{align}
and 
\begin{align}\label{eq:surftwodirectiontwowell}
	 \Vert D^2 u\Vert_{TV(R)} \lesssim H N.
\end{align}
 Here $E_{\rm el}^{(p)}(u;R)$ denotes the elastic energy of the geometrically linearized problem with respect to the wells given by $A,B$. 
\item[(ii)] \emph{One (degenerate) symmetrized rank-one direction:}
Let $A = \diag (1, 0)$ and $B = \diag(-1,0)$. Then there exists $u\in W_0^{1,\infty}(R;\R^2)$ with $\nabla u\in BV(R;\R^{2\times 2})$ such that
\begin{equation} \label{eq:onedirectiontwowell}
E_{\rm el}^{(p)}(u;R) \lesssim   \frac{L^{2p+1}}{H^{2p-1} N^{2p}}.
\end{equation}
Similarly, the estimate \eqref{eq:surftwodirectiontwowell} holds.
\end{itemize}
\end{lemma}

We mention that an optimization of $N$ in the cases (i) and (ii) leads to the scaling $E_\epsilon^{(2)}(u;R) \sim \epsilon^{2/3}  $ and  $E_\epsilon^{(2)}(u;R) \sim \epsilon^{4/5}  $, respectively.
The proof of these (well-known) upper bounds requires careful and technical constructions and can be found in the literature \cite{CC15}. 
Still, as some structural results of these constructions will be needed in the three-dimensional setting,
we recall some of the to us most relevant aspects.

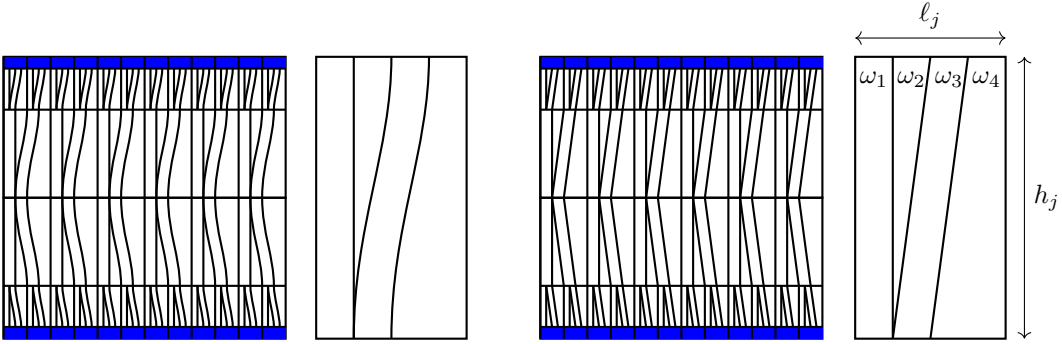
\begin{figure}[h] 
\centering 
\hspace{0.03\textwidth}
\begin{minipage}[b]{0.212\textwidth}
  \centering
\begin{tikzpicture}[scale=0.3125]

  \begin{scope}
	
\draw[thick] (0,0) rectangle (12,6);
\fill[blue] (0,5.5) rectangle (12,6);

\draw[thick] (0,5.5) -- (12,5.5); 
\draw[thick] (0,3.75) -- (12,3.75);

\foreach \i in {0,2,...,10} {
    \draw[thick]  (\i, 3.75) -- (\i, 0);
}
\foreach \i in {0,1,...,11} {
    \draw[thick]  (\i, 5.5) -- (\i, 3.75);
}
\foreach \i in {0,1,...,11} {
    \draw[thick]  (\i, 6) -- (\i, 5.5);
}

\foreach \i in {0,2,...,10} {
   \draw[thick] (\i+0.5,0) .. controls (\i+0.5,1.5) and (\i+0.5,2.5) .. (\i+0.5,3.75);
   \draw[thick] (\i+0.5,0) .. controls (\i+0.5,1.5) and (\i+1,2.5) .. (\i+1,3.75);  
    \draw[thick] (\i+1,0) .. controls (\i+1,1.5) and (\i+1.5,2.5) .. (\i+1.5,3.75);

}

\foreach \i in {0,1,...,11} {
    \draw[thick] (\i+0.25,3.75) -- (\i+0.25,5.5); 
    
   \draw[thick] (\i+0.25,3.75) .. controls (\i+0.25,4.25) and (\i+0.5,5) .. (\i+0.5,5.5);  
    \draw[thick] (\i+0.5,3.75) .. controls (\i+0.5,4.25) and (\i+0.75,5) .. (\i+0.75,5.5);  
  
}

\end{scope}

\begin{scope}[yscale=-1]  
 
\draw[thick] (0,0) rectangle (12,6);
\fill[blue] (0,5.5) rectangle (12,6);
 
\draw[thick] (0,5.5) -- (12,5.5); 
\draw[thick] (0,3.75) -- (12,3.75);  
 
\foreach \i in {0,2,...,10} {
    \draw[thick]  (\i, 3.75) -- (\i, 0);
}
\foreach \i in {0,1,...,11} {
    \draw[thick]  (\i, 5.5) -- (\i, 3.75);
}
\foreach \i in {0,1,...,11} {
    \draw[thick]  (\i, 6) -- (\i, 5.5);
}
 
\foreach \i in {0,2,...,10} {
   \draw[thick] (\i+0.5,0) .. controls (\i+0.5,1.5) and (\i+0.5,2.5) .. (\i+0.5,3.75); 
   \draw[thick] (\i+0.5,0) .. controls (\i+0.5,1.5) and (\i+1,2.5) .. (\i+1,3.75); 
   \draw[thick] (\i+1,0) .. controls (\i+1,1.5) and (\i+1.5,2.5) .. (\i+1.5,3.75); 
}

\foreach \i in {0,1,...,11} {
    \draw[thick] (\i+0.25,3.75) -- (\i+0.25,5.5); 
    \draw[thick] (\i+0.25,3.75) .. controls (\i+0.25,4.25) and (\i+0.5,5) .. (\i+0.5,5.5); 
    \draw[thick] (\i+0.5,3.75) .. controls (\i+0.5,4.25) and (\i+0.75,5) .. (\i+0.75,5.5); 
}

\end{scope}

\end{tikzpicture}
 
\end{minipage}
\hfill
\begin{minipage}[b]{0.24\textwidth}
  \centering
 \begin{tikzpicture}[scale=1]

  \begin{scope}  
\draw[thick] (0,0) rectangle (2,3.75);
 
\foreach \i in {0} {
   \draw[thick] (\i+0.5,0) .. controls (\i+0.5,1.5) and (\i+0.5,2.5) .. (\i+0.5,3.75); 
   \draw[thick] (\i+0.5,0) .. controls (\i+0.5,1.5) and (\i+1,2.5) .. (\i+1,3.75); 
   \draw[thick] (\i+1,0) .. controls (\i+1,1.5) and (\i+1.5,2.5) .. (\i+1.5,3.75); 
}

\end{scope}

\end{tikzpicture}
 
\end{minipage}
\hfill
\begin{minipage}[b]{0.24\textwidth}
  \centering
\begin{tikzpicture}[scale=0.3125]

  \begin{scope} 
\draw[thick] (0,0) rectangle (12,6);
\fill[blue] (0,5.5) rectangle (12,6);
 
\draw[thick] (0,5.5) -- (12,5.5);  
\draw[thick] (0,3.75) -- (12,3.75);  

\foreach \i in {0,2,...,10} {
    \draw[thick]  (\i, 3.75) -- (\i, 0);
}
\foreach \i in {0,1,...,11} {
    \draw[thick]  (\i, 5.5) -- (\i, 3.75);
}
\foreach \i in {0,1,...,11} {
    \draw[thick]  (\i, 6) -- (\i, 5.5);
}

\foreach \i in {0,2,...,10} {
   \draw[thick] (\i+0.5,0) -- (\i+0.5,3.75);  
   \draw[thick] (\i+0.5,0) -- (\i+1,3.75);  
    \draw[thick] (\i+1,0) -- (\i+1.5,3.75);   
 
}

\foreach \i in {0,1,...,11} {
    \draw[thick] (\i+0.25,3.75) -- (\i+0.25,5.5);  
   \draw[thick] (\i+0.25,3.75) -- (\i+0.5,5.5);  
    \draw[thick] (\i+0.5,3.75) -- (\i+0.75,5.5);

}

\end{scope}

\begin{scope}[yscale=-1]  
 
\draw[thick] (0,0) rectangle (12,6);
\fill[blue] (0,5.5) rectangle (12,6);
 
\draw[thick] (0,5.5) -- (12,5.5);  
\draw[thick] (0,3.75) -- (12,3.75);  
 
\foreach \i in {0,2,...,10} {
    \draw[thick]  (\i, 3.75) -- (\i, 0);
}
\foreach \i in {0,1,...,11} {
    \draw[thick]  (\i, 5.5) -- (\i, 3.75);
}
\foreach \i in {0,1,...,11} {
    \draw[thick]  (\i, 6) -- (\i, 5.5);
}
 \foreach \i in {0,2,...,10} {
   \draw[thick] (\i+0.5,0) -- (\i+0.5,3.75); 
   \draw[thick] (\i+0.5,0) -- (\i+1,3.75); 
   \draw[thick] (\i+1,0) -- (\i+1.5,3.75); 
}

\foreach \i in {0,1,...,11} {
    \draw[thick] (\i+0.25,3.75) -- (\i+0.25,5.5); 
    \draw[thick] (\i+0.25,3.75) -- (\i+0.5,5.5); 
    \draw[thick] (\i+0.5,3.75) -- (\i+0.75,5.5); 
}

\end{scope}

\end{tikzpicture}
 
\end{minipage}
\hfill
\begin{minipage}[b]{0.24\textwidth}
  \centering
 \begin{tikzpicture}[scale=1]
  \begin{scope}  
\draw[thick] (0,0) rectangle (2,3.75); 
\foreach \i in {0} {
   \draw[thick] (\i+0.5,0) -- (\i+0.5,3.75); 
   \draw[thick] (\i+0.5,0) -- (\i+1,3.75); 
   \draw[thick] (\i+1,0) -- (\i+1.5,3.75); 
}

 \node at (0.25,3.45) {$\omega_1$};
\node at (0.75,3.45) {$\omega_2$};
\node at (1.25,3.45) {$\omega_3$};
\node at (1.75,3.45) {$\omega_4$};

\draw[<->] (2.25,0) -- (2.25,3.75) node[midway,right] {$h_j$};  
\draw[<->] (0,4) -- (2,4) node[midway,above] {$\ell_j$};

\end{scope}

\end{tikzpicture}

\end{minipage} 
  \caption{ Visualization of the constructions in Lemma~\ref{lem:firstorder2d} for $H = L = 1$. Whereas the interfaces in Lemma~\ref{lem:firstorder2d}(i) can be straight lines, the interfaces in (ii) necessarily need to be curved. The picture shows $j_0 =2$ layers of building blocks, where the blue part corresponds to the cut-off layer. In Lemma~\ref{lem:firstorder2d}(i), we have $\partial_1   u_2 = -2$ on $\omega_1$ and $\omega_3$, and on $\omega_2$ and $\omega_4$, it holds that $\partial_1   u_2 = 2$. }
    \label{fig:branchingpicture}
\end{figure}

\begin{proof}[Proof of Lemma~\ref{lem:firstorder2d}]
\emph{Step 1: Argument for (i).}
As $A$ and $B$ are symmetrized rank-one connected with $A-B = - 2 e_1 \odot e_2 $, we wish to show the statement by resorting to Lemma~\ref{lem:coordinatetransform}(ii).
By \cite[Lemma~3.2]{RT23}, there exists a function $v \in W^{1,\infty}_0 ( R; \R^2) $ with $\nabla v\in BV(R;\R^{2\times 2})$ satisfying
\begin{align*} 
	\int_{R}  \dist (\nabla v,  \tilde K )^p \di x  \lesssim \frac{L^{p+1}}{N^{p} H^{p-1}} \qquad \text{ and } \quad   \Vert  D^2 v \Vert_{{\rm TV} (R)} \lesssim   H N,
\end{align*}
where 
$\tilde K \defas \{ \diag (1/2,0) , \diag (-1/2,0) \}$. (The construction allows one to estimate the elastic and surface energy separately.)
Notice that $$ S \tilde K = \frac{1}{2} \{ - e_2 \otimes e_1 , e_2 \otimes e_1 \} \text{ for } S = \begin{pmatrix}
	0 & 1 \\
	-1 & 0 \\
\end{pmatrix}.$$
Hence, setting $u (x_1,x_2) = 2 S v (x_1,x_2) + e_2 \otimes e_1 x$,
we can conclude the statement by Lemma~\ref{lem:coordinatetransform}(ii) and a further linear transformation.

\emph{Step 2: Argument for (ii).}
The claim from (ii) is proved for $H = 1$, $L = 1$, and $p = 2$ in \cite[Section~4]{RRTT24}, which itself is based on \cite[Section~2]{CC15}. A straight-forward adaptation (similarly to \cite[Lemma~2.3]{RT23}) shows the lemma for a general $H$, $L$, and $p \in [1,+\infty)$.
For later reference, we highlight the following bounds on the associated deformation: The key observation from \cite{CC15} is that in the case of only one rank-one connection, by a genuinely vectorial argument, it is possible to ensure that all components of the strain tensor attain the desired values from $K= \{A,B\}$ (apart from a cut-off layer), except for the $\partial_2u_2$ component in which the energy concentrates. In order to ensure this, the $\partial_1 u_1$ component is prescribed according to    the ${e}_{11}$-entries in $K$ and $u_1$ is obtained by integrating the equation $\partial_2 u_1 + \partial_1 u_2=0$ while preserving the boundary data. In this process, both $u_1$ and $u_2$ are non-trivial and  
 additionally satisfy the following bounds
\begin{align}\label{eq:boundsfunctionepsilon45}
	\Vert \partial_2 u_1 \Vert_{L^q(R)}^q + \Vert \partial_1 u_2 \Vert_{L^q(R)}^q \leq C \frac{L^{q+1}}{N^q H^{q-1}} \quad \text{ and } \quad \Vert \partial_2 u_2 \Vert_{L^q(R)}^q \leq C \frac{L^{2q+1}}{H^{2q-1} N^{2q}}.
\end{align}
 The constant $C>0$ depends on $q$, however it can be chosen uniformly among values in the interval $[1,q^*]$ for $q^* > p $. 

For later purposes, we give more details on the structure of $u$ from Lemma~\ref{lem:firstorder2d}, which is also visualized in Figure~\ref{fig:branchingpicture}.  
 The upper bounds in Lemma~\ref{lem:firstorder2d} are defined by a branching construction.
 Firstly, they satisfy the symmetry condition $u( z_1, z_2 - H/2)  = u( z_1, \vert z_2 - H/2 \vert) $ for $z \in R$.
The competitors are defined by attaching building blocks $\omega_{j,k}$, where $\omega_{j,k} \subset [0,L]\times[\frac{H}{2},H]$ are given by 
\begin{align} 
\omega_{j,k}:= 
(k\ell_j,y_j)+[0,\ell_j]\times[0,h_j]  \quad \text{ for }  j \in \{0,\dots,j_0 \}, \ j_0 \in \N  \text{ and } k\in \{0,\dots, N2^j-1\} . \label{eq:buildingblock} 
\end{align}
Here, for some fixed $\theta \in (2^{-2q^*/(2q^*-1)}, 2^{-1})  $, the variables are given by
\begin{align} 
\ell_j \defas \frac{L}{2^j N},
\quad y_j \defas H-\frac{H}{2}\theta^j \quad \text{and} \quad
h_j \defas y_{j+1}-y_j=\theta^j\frac{H(1-\theta)}{2}. \label{eq:blocklengthdef} 
\end{align} 
Then, $j_0$ is chosen in such a way that on the remaining part of $R$, i.e.,
\begin{align}
	\bigcup_{k = 0} ^{N 2^{j_0}-1} (k\ell_{j_0},y_{j_0+1})+[0,\ell_{j_0}]\times \left(\Big[0,\frac{H}{2}\theta^{j_0+1}\Big] \cup  \Big[H-\frac{H}{2}\theta^{j_0+1},H\Big] \right) ,  \label{eq:cutoffbuildingblock} 
\end{align}
 the elastic and surface energy are sufficiently small, see e.g.\ \cite[Proof of Lemma~3.2]{RT23} or \cite[Lemma~4.4]{RRTT24} for details. 
 Finally, we mention that the building blocks can be subdivided into four blocks that generate the branching pattern, see Figure~\ref{fig:branchingpicture}. 
\end{proof}

Having presented the upper bound constructions in two dimensions, we make use of these results to provide a
preliminary bound in three dimensions with Dirichlet data on the full boundary.
Here, we highlight once again that 
Dirichlet boundary conditions have to be satisfied on every side of the domain, which is the crucial challenge in our setting. While three-dimensional analogues of two-dimensional constructions have been deduced in settings without gauges (see, for instance, \cite{RT23}) and
 directly transfer to the setting with gauges in the presence of \emph{two} symmetrized rank-one connections, in the setting with only \emph{one} symmetrized rank-one connection, new difficulties arise. The central difficulty here is that the deformation is genuinely vectorial with \emph{two} non-zero components already in two dimensions. This provides substantial new challenges, which will be addressed in the following lemma.

\begin{lemma}[Preliminary construction for a two-well problem in three dimensions]\label{lem:2well3d}
	Let $\omega = [0,L] \times [0,H_1] \times [0,H_2]$ and $A,B \in K$ such that
	$A= e_1 \odot e_1$ and $B = -A$ and $1 \geq H_2 > H_1$, and assume that 
$N >  4L H_1^{-1} $. 
	Then there exists $u\in W_0^{1,\infty}(\omega;\R^3)$ with $\nabla u\in BV(\omega;\R^{3\times 3})$ such that
\begin{equation} \label{eq:3donedirectiontwowell}
E_\epsilon^{(p)}(u;\omega) \lesssim L^{p+2} H_1^{-p+1 } N^{-p-1}  +  \frac{H_2 L^{2p+1}}{H_1^{2p-1} N^{2p}} +\epsilon  \left( H_1 H_2 N  +  L^{2}     N^{-1} \right) . 
\end{equation}
Moreover, $u$ satisfies the bound
\begin{align}\label{ineq:uniformboundpartial2}
	 \Vert \partial_2 u \Vert_{L^p(\omega)}^p  \lesssim \frac{H_2 L^{p+1}}{H_1^{p-1} N^{p}}.  
\end{align}
\end{lemma}

Let us comment on this result by focusing on the case $p=2$.
Notice that due to the first term on the right-hand side of \eqref{eq:3donedirectiontwowell} the optimal choice $N \sim \epsilon^{-1/4} $ leads to the scaling $E_\epsilon^{(2)} \sim \epsilon^{3/4}$, which is \emph{better} than $\epsilon^{2/3}$ but slower than $\epsilon^{4/5}$.
 However, as we will show in the proof of Theorem \ref{thm:mainresultupper} below, 
the different scalings on the right-hand side of Lemma \ref{lem:2well3d} still allow for
an optimal scaling for second order data in the case if the finer branching construction can only be formed by one compatibility direction, at least if $p \in [1,3]$, see Proposition~\ref{prop:upperboundcons} below.
This is surprising, but \emph{consistent} with the construction of \cite[Section~6~and~Lemma~6.1]{RT23}, where the elastic energy depends on the full gradient. 
Therein, for branching constructions of higher order, a constant extension to the third dimension with a naive cut-off was sufficient for $p \in [1, 2]$ whereas in our setting such a construction leads to a scaling which is too slow. 

 We remark that, the limitations in $p$ to small values does not come as a complete surprise. In general, it is known that for $p=1$ simple laminates (without branching) and branching structures share the same scaling behaviour in the absence of gauges (see, for instance, the discussion in the introduction of \cite{RT23}), while for all larger values of $p>1$ branching is preferable. Hence, it can be hoped that for small values of $p$, in a suitable neighbourhood of $p=1$, also concatenations of optimal and non-optimal structure may suffice to determine an overall optimal scaling behaviour. 
\begin{proof}
	The proof is divided into two steps. Here the competitor is constructed in step 1, the relevant energy estimates are provided in the second step.

	\emph{Step 1 (Construction of the competitor):}
	Our goal is to rely on the rotation argument in \cite[Section~6.2]{RT23}
	such that we can  reduce the problem to the constructions in two dimensions as stated in Lemma~\ref{lem:firstorder2d}.
	We decompose $[0,H_1] \times [0,H_2]$ into four subsets $D_i$, $i = 1, \dots, 4$, visualized in Figure~\ref{fig:decompositionofH1H2}.
\begin{figure}
	\centering
	\begin{tikzpicture}[thick, scale= 1]
  \draw (0, 0)   to[out=0, in=180] (2, 0);
  \draw (0, 5)   to[out=0, in=180] (2, 5);
  \draw (0, 0) to[out=90, in=-90] (0, 5); 
 \draw (2, 0) to[out=90, in=-90] (2, 5); 

  \draw  [dash pattern={on 2pt off 2pt on 2pt off 2pt}] (0, 0) to (1, 1);
  \draw  [dash pattern={on 2pt off 2pt on 2pt off 2pt}] (1, 1) to (2, 0);

   \draw  [dash pattern={on 2pt off 2pt on 2pt off 2pt}] (0, 5) to (1, 4);
  \draw  [dash pattern={on 2pt off 2pt on 2pt off 2pt}] (1, 4) to (2, 5);

    \draw  [dash pattern={on 2pt off 2pt on 2pt off 2pt}] (1, 4) to (1, 1);

 \node[below] at (1,0) {$H_1$};
\node[left] at (0,2.5) {$H_2$};

\node at (1,0.5) {$D_1$}; 

\node at (1,4.5) {$D_2$}; 

\node at (0.5,2.5) {$D_3$}; 

\node at (1.5,2.5) {$D_4$};

\end{tikzpicture}
\caption{Decomposition of $[0,H_1] \times [0,H_2]$, which denotes a subset of the $x_2$-$x_3$-plane.}
\label{fig:decompositionofH1H2}
\end{figure}
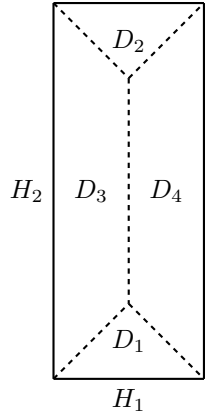
More precisely, we define
\begin{align*}
	D_1 & \defas \left\{ \left\vert x_2 - \frac{H_1}{2} \right\vert < \left\vert x_3 - \frac{H_1}{2} \right\vert \right\},  && D_2 \defas \left\{ \left\vert x_2 - \frac{H_1}{2} \right\vert < \left\vert (H_2 - x_3) - \frac{H_1}{2} \right\vert \right\}, \\
D_3 & \defas \left[0,\frac{H_1}{2}\right] \times [0,H_2] \setminus ( D_1 \cup D_2 ) ,  && D_4  \defas \left[\frac{H_1}{2}, H_1\right] \times [0,H_2] \setminus ( D_1 \cup D_2 ).
\end{align*}
Let $\tilde u\colon  (0,L) \times (0,H_1) \to \R^2$ be defined via Lemma~\ref{lem:firstorder2d}(ii), and denote $v\colon  (0,L) \times (0,H_1) \times (0,H_2) \to \R^3$ by 
\begin{align}\label{eq:def3dfirstcase}
	\begin{cases}
		\begin{pmatrix}
	v_1 (x)\\
	v_3 (x) \\
\end{pmatrix}
=    
\tilde u \big(  x_1  , x_3  \big)  & \quad   \text{if } x \in [0,L] \times (D_1 \cup D_2) , \\
\begin{pmatrix}
	v_1 (x)\\
	v_2 (x) \\
\end{pmatrix}
=    
\tilde u \big(  x_1  , x_2  \big)  & \quad   \text{if } x \in [0,L] \times (D_3 \cup D_4),
	\end{cases}
\end{align}
and set the remaining component to be zero. The function attains zero  Dirichlet boundary conditions due to the fact that $\tilde u$ satisfies zero Dirichlet boundary conditions and its symmetry property.
However, a further cut-off argument is required for
$v_2$ and $v_3$ as, a priori, no continuity is ensured along the interior interfaces of $[0,L] \times \partial (D_1 \cup D_2)$.
 We define these cut-off functions for $v_2$ and $v_3$ along the interface at $D_1$ by
\begin{align*}
	\theta^{(2)}_s (x_2,x_3) &\defas \varphi \left( \big(   \rho(x_2,x_3) - ( \vert x_3 - \frac{H_1}{2} \vert ) \big) s^{-1}  \right) \quad  \text{ and} \\
	\theta^{(3)}_s (x_2,x_3) &\defas \varphi \left( \big(    \rho(x_2,x_3)  - ( \vert x_2  - \frac{H_1}{2} \vert   ) \big) s^{-1}  \right),
\end{align*}
where 
$\rho (x_2,x_3) \defas \max \{ \vert x_2 - H_1/2 \vert, |x_3 - H_1/2 \vert \} $, 
$\varphi \in C^\infty(\R) $ satisfies $\varphi(0) = 0$ and $\varphi(z) = 1$ for $z \geq 1$, and $s >0$ is the length of the cut-off to be fixed below.
 The definition via $\rho$ ensures that $\theta^{(3)}_s = 0 $ on $D_3 \cup D_4$ and $\theta^{(2)}_s = 0$ on $D_1$. 
The cut-off at the interface of $[0, L ] \times \partial D_2$ can be defined by symmetry, and in what follows it suffices to provide the relevant energy estimates at the interface $ [0,L] \times \partial D_1$. The deformation
\begin{align}\label{def:constcut-off}
	u(x) \defas \begin{pmatrix}
		v_1 (x) \\
		\theta^{(2)}_s (x_2,x_3) v_2 (x) \\
		\theta^{(3)}_s (x_2,x_3) v_3 (x) \\
	\end{pmatrix}
\end{align}
therefore attains zero boundary conditions and is Lipschitz continuous with $\nabla u \in BV$.

\emph{Step 2.1 (Auxiliary computations):}
A computation of the gradient of $u$ yields
\begin{align}\label{eq:gradient2}
	\nabla u(x) = \begin{cases}
		  \begin{pmatrix}
		 \partial_1 \tilde u_1 & \partial_2 \tilde u_1 & 0 \\
		 \theta^{(2)}_s \partial_1 \tilde u_2 & \theta^{(2)}_s  \partial_2 \tilde u_2 + \partial_2 (\theta^{(2)}_s )  \tilde u_2  & \partial_3 (\theta^{(2)}_s )  \tilde u_2  \\
		  0 & 0 & 0 \\
		\end{pmatrix}    &   \quad   \text{if } x \in [0,L] \times (D_3 \cup D_4) ,\\
  \begin{pmatrix}
		 \partial_1 \tilde u_1 &  0 &    \partial_2 \tilde u_1  \\
		  0 & 0 & 0 \\
		  		 \theta^{(3)}_s  \partial_1 \tilde u_2 &  \partial_2 (\theta^{(3)}_s )    \tilde u_2  &  \theta^{(3)}_s  \partial_2 \tilde u_2+ \partial_3 (\theta^{(3)}_s )    \tilde u_2  \\
		\end{pmatrix}  & \quad   \text{if } x \in [0,L] \times (D_1 \cup D_2) . \\
	\end{cases}
\end{align}
For later purposes, we compute the symmetric part of $\nabla u$ in the second case of \eqref{eq:gradient2}, which is given by
\begin{align}\label{eq:cutoffestimate}
	\sym (\nabla u) &=  \begin{pmatrix}
		 \partial_1 \tilde u_1 &  0 &   \frac{1}{2} ( \partial_2 \tilde u_1 +  \partial_1 \tilde u_2 ) \\
		  0 & 0 & 0 \\
		  	 \frac{1}{2} ( \partial_2 \tilde u_1 +  \partial_1 \tilde u_2 ) & 0 &  \theta^{(3)}_s   \partial_2 \tilde u_2  \\
		\end{pmatrix} \notag \\ & \qquad  +  \begin{pmatrix}
		 0 &  0 &     \frac{1}{2}(\theta^{(3)}_s -1) \partial_1 \tilde u_2  \\
		  0 & 0 & \frac{1}{2}  \partial_2 (\theta^{(3)}_s )    \tilde u_2  \\
		  		  \frac{1}{2}(\theta^{(3)}_s -1) \partial_1 \tilde u_2  & \frac{1}{2}  \partial_2 (\theta^{(3)}_s )    \tilde u_2  &   \partial_3 (\theta^{(3)}_s )    \tilde u_2  \\
		\end{pmatrix}  .
\end{align}
Given a function $f \colon \R^2 \to \R_+$, a change of variables on the cut-off part $C^{(i)} \defas [0,L] \times \{ \theta_s^{(i)} \in (0,1) \}$, $i \in \{2,3\}$,   yields
\begin{align}\label{eq:auxestimate}
	\int_{C^{(i)}} f(x_1, x_j) \di x = \int_0^L \int_{\{ \theta_s^{(i)} \in (0,1) \}}  f(x_1,x_j) \di x_1 \di x_2 \di x_3 \lesssim s \int_0^L \int_0^{H_1} f(x_1,x_j) \di x_1 \di x_j
\end{align}
for $j \in \{2,3\}$.

\emph{Step 2.2 (Estimates on the cut-off layer):}
Special care is needed for the behaviour of $\sym (\nabla u)$ in the cut-off part $C^{(i)} $, $i \in \{2,3\}$.

\emph{Step 2.2.1 (The elastic part):}
By symmetry it suffices to provide estimates for the second case in \eqref{eq:gradient2}.
The first term on the right-hand side of \eqref{eq:cutoffestimate} enjoys a good control, and can be estimated as in \eqref{eq:goodpartsym} below.
Next, we compute the elastic contribution of the second term on the right-hand side of \eqref{eq:cutoffestimate}
by controlling the $L^p$-norm over the cut-off layer $C^{(i)}$, $i \in \{2,3\}$. 
 The uniform bound of $\theta_s^{(j)}$, \eqref{eq:auxestimate}, and \eqref{eq:boundsfunctionepsilon45} imply that 
\begin{IEEEeqnarray}{r;c;l}
\label{eq:estimatecutoff1}
	\int_{C^{(i)}} \vert \theta_s^{(3)} -1 \vert^p \vert \partial_1 \tilde  u_2 \vert^p \di x 
&   \leq & 	\int_{C^{(i)}} \vert \partial_1 \tilde  u_2 \vert^p \di x 
\stackrel{\eqref{eq:auxestimate}}{\leq} Cs \int\limits_{0}^{L} \int\limits_{0}^{H_1} |\partial_1 \tilde{u}_2|^p \di x  \notag	\\
&	\stackrel{\eqref{eq:boundsfunctionepsilon45}}{\leq} & C s  L^{p+1} H_1^{-p+1} N^{-p}.  
\end{IEEEeqnarray}
Moreover, using that the derivative of the cut-off function is controlled by $s^{-1}$,  \eqref{eq:auxestimate}, the fundamental theorem of calculus, and \eqref{eq:boundsfunctionepsilon45}
imply that
\begin{IEEEeqnarray}{r;c;l}\label{eq:estimatecutoff2}
	\int_{C^{(i)}} \vert \partial_i (\theta^{(3)}_s ) \vert^p \vert  \tilde u_2 \vert^p \di x 
	 & {\overset{ \eqref{eq:auxestimate}}{\lesssim}} &  s^{-p} s   \int_0^L \int_0^{H_1}  \vert \tilde  u_2 \vert^p \di x \notag \\
	&  \lesssim & s^{-p+1} \int_0^L \int_0^{H_1}  \left\vert \int_0^{x_2} \partial_2 \tilde  u_2 (x_1,r) \di r \right\vert^p \di x_2 \di x_1  \notag \\
		 & {\overset{\text{Jensen}}{\lesssim}} & s^{-p+1}   H_1^{p} \int_0^L \int_0^{H_1} \left\vert \partial_2 \tilde  u_2 (x_1,r) \right\vert^p \di r    \di x_1  \notag \\
			  &  {\overset{\eqref{eq:boundsfunctionepsilon45}}{\lesssim}} & s^{-p+1} L^{2p+1}   H_1^{-p+1} N^{-2p} .  
\end{IEEEeqnarray}
An optimization shows that the best $s$ in 
\eqref{eq:estimatecutoff1} and \eqref{eq:estimatecutoff2} is

\begin{equation}
\label{eq:s}
  s \sim L   N^{-1}  .
  \end{equation} 
  
For that choice the right-hand side of the latter inequalities becomes
$ L^{p+2} H_1^{-p+1 } N^{-p-1} $.

\emph{Step 2.2.2 (The surface energy):}
Our goal is to control the surface energy in terms of the surface energy for $v$ that is generated by $D^2 \tilde u$.
By Fubini's theorem it suffices to control the total variation of the distributional derivative of the components in \eqref{eq:gradient2}, where the cut-off function appears.
The other parts enjoy a good control, as shown in \eqref{eq:surface} below.

On the one hand, we need to control the total variation of $D (\theta_s^{(3)} \partial_i \tilde u_2)$ for $i = 1,2$.
 By using the $L^\infty$-bound of the cut-off function, \eqref{eq:auxestimate}, 
\eqref{eq:boundsfunctionepsilon45}, and \eqref{eq:surftwodirectiontwowell} combined with Fubini's theorem
 we find that 
  \begin{align}\label{eq:3dsurface1}
	&\quad \ \Vert D (\theta_s^{(3)} \partial_1 \tilde u_2 ) \Vert_{TV([0,L] \times D_1 )} + \Vert D (\theta_s^{(3)} \partial_2 \tilde u_2 ) \Vert_{TV([0,L] \times D_1 )}\notag \\
	 & \leq \Vert \theta_s^{(3)} \Vert_{L^\infty([0,L] \times D_1)}  \Vert D^2 \tilde  u \Vert_{TV([0,L] \times D_1)}  +  \Vert \nabla \theta_s^{(3)} \Vert_{L^\infty(C^{(i)})} (\Vert  \partial_1 \tilde u_2 \Vert_{L^1(C^{(i)})} + \Vert  \partial_2 \tilde u_2 \Vert_{L^1(C^{(i)})} )  \notag  \\
	 & \lesssim    H_1^2 N + s^{-1} s \frac{L^{2}}{N}  + s^{-1} s \frac{L^{3}}{N^2 H_1}     \leq H_1 H_2 N  +  L^{2}        N^{-1} + L^3 H_1^{-1} N^{-2}  \lesssim H_1 H_2 N  +  L^{2}     N^{-1}   ,   
 \end{align} 
 where we additionally used that $H_1 < H_2$ and the assumption on the size of $N$. 

On the other hand, we need to control the total variation of $ D (\partial_i \theta_s^{(3)}  \tilde u_2 ) $ for $i = 2,3$.
In view of \eqref{eq:auxestimate} and \eqref{eq:boundsfunctionepsilon45}, a similar computation to \eqref{eq:estimatecutoff2} and \eqref{eq:3dsurface1} implies that
  \begin{align}\label{eq:3dsurface2}
	&\quad \ \Vert D (\partial_i \theta_s^{(3)}  \tilde u_2 ) \Vert_{TV([0,L] \times D_1 )} \notag \\ 
	& \leq \Vert \partial_i \theta_s^{(3)} \Vert_{L^\infty(C^{(i)})}  \Vert  \nabla \tilde u_2 \Vert_{L^1(C^{(i)})}  +  \Vert \nabla^2 \theta_s^{(3)}\Vert_{L^\infty(C^{(i)})} \Vert  \tilde u_2 \Vert_{L^1(C^{(i)})}   \notag \\
		& \lesssim   L^3 H_1^{-1} N^{-2} + L^2   N^{-1} + s^{-2} s L^3  N^{-2} \notag  \\
		&  \lesssim  L^2   N^{-1} , 
 \end{align}
 where in the third line we used the estimates for the two components of $\nabla \tilde{u}_2$ and the fundamental theorem for the estimate of $\tilde{u}_2$ and the last estimate follows by the choice of $s$, see \eqref{eq:s}, and the assumption on $N$.

\emph{Step 2.3 (Control of remaining bulk part):}
As every relevant contribution of the cut-off has been estimated in the previous step, we need to provide estimates on the remaining bulk part.
Using \eqref{eq:cutoffestimate}, \eqref{eq:onedirectiontwowell}, and Fubini's theorem we find that for $i \in \{2,3\}$
\begin{align}\label{eq:goodpartsym} 
	\int_{[0,L] \times (D_1 \setminus C^{(i)}) } \dist^p( \sym (\nabla u ; K)) \di x \lesssim  \frac{H_1 L^{2p+1}}{H_1^{2p-1} N^{2p}} .
\end{align}
Similarly, we find that
\begin{align}\label{eq:goodpartsym2} 
	\int_{[0,L] \times (D_2 \setminus C^{(i)}) } \dist^p( \sym (\nabla u ; K)) \di x \lesssim   \frac{H_2 L^{2p+1}}{H_1^{2p-1} N^{2p}}.
\end{align}
As before, we can estimate the surface energy by using the chain rule and Fubini's theorem, i.e.,
\begin{align} \label{eq:surface}
	\Vert D^2 v \Vert_{TV(\omega)} \lesssim \left( \Vert D^2 \tilde u \Vert_{TV((0,L) \times (0,H_1))}  +  \vert  \partial ( (0,L) \times \partial D_1 )  \vert  \right) \lesssim H_2 H_1 N  .
\end{align} 

\emph{Step 3 (Conclusion):}
A combination of the estimates \eqref{eq:estimatecutoff1} (for $s = L  N^{-1} $), \eqref{eq:estimatecutoff2}, \eqref{eq:3dsurface1}, \eqref{eq:3dsurface2},
\eqref{eq:goodpartsym}, \eqref{eq:goodpartsym2}, \eqref{eq:surface}, and $H_2 H_1^{-1} \geq 1$,  lead to the control
\begin{equation*}
E_\epsilon^{(p)}(u;\omega) \lesssim  L^{p+2} H_1^{-p+1 }  N^{-p-1}  +  \frac{H_2 L^{2p+1}}{H_1^{2p-1} N^{2p}} +\epsilon  (H_2 H_1 N + L^2 N^{-1} ) . 
\end{equation*}
 It remains to provide the argument for \eqref{ineq:uniformboundpartial2}. 
In view of \eqref{eq:gradient2}, the estimates \eqref{eq:estimatecutoff1} and \eqref{eq:estimatecutoff2} (for $s = L  N^{-1} $) control the cut-off part, whereas the estimate \eqref{eq:boundsfunctionepsilon45} controls the bulk part (similarly to \eqref{eq:goodpartsym} and \eqref{eq:goodpartsym2}). This implies that
 \begin{align*}
	\Vert \partial_2 u \Vert_{L^p(\omega)}^p & \lesssim \int_{\omega} \vert \partial_2 \tilde u_1 \vert^p + \vert \theta^{(2)}_s  \partial_2 \tilde u_2 \vert^p + \vert \partial_2 (\theta^{(2)}_s )  \tilde u_2 \vert^p + \vert \partial_2 (\theta^{(3)}_s )    \tilde u_2 \vert^p \di x  \\
	 & \lesssim \frac{H_2 L^{p+1}}{H_1^{p-1} N^{p}} +  \frac{H_2 L^{2p+1}}{H_1^{2p-1} N^{2p}} +  L^{p+2}  H_1^{-p+1 }    N^{-p-1} ,
\end{align*} 
where the second term of the right-hand side is controlled by the first one due to
  $N >  4L H_1^{-1} $. The latter assumption on $N$ and $1 < H_2 H_1^{-1}$ also yields 
\begin{align*}
	L^{p+2}  H_1^{-p+1 }    N^{-p-1}  < L H_1^{-1} H_2 L^{p+1}  H_1^{-p+1 }    N^{-p-1} < H_2 L^{p+1}  H_1^{-p+1 }  N^{-p}.
\end{align*} 
This concludes the argument for \eqref{ineq:uniformboundpartial2}, and hence the proof.
\end{proof} 
The two-dimensional branching construction with $\epsilon^{4/5}$-scaling in two dimensions enforces both components $\tilde u_1$ and $\tilde u_2$ of the competitor to be nonzero.
Hence, to ensure continuity, a further cut-off at the interfaces visualized in Figure~\ref{fig:decompositionofH1H2} is required once the three-dimensional rotation argument as in \cite[Section~6]{RT23} has been applied. The cut-off enforces a competition between the terms on the right-hand side of
\eqref{eq:estimatecutoff1} and \eqref{eq:estimatecutoff2}, which
in the case $p >1$, leads to a worse scaling than $N^{-2p}$ in \eqref{eq:3donedirectiontwowell}. 

 In order to avoid the cut-off one might rotate the direction of the branching continuously, which is briefly discussed in the following proposition for a cylinder. 
\begin{proposition}[Optimal first order constructions in cylinders]\label{prop:optimalfirst-order}
	Consider the cylinder
	$\omega = (0,1) \times B_1(0)$,  the wells $A,B \in K$ such that
	$A= e_1 \odot e_1$ and $B = -A$, and assume that 
$N >  4$. 
	Then there exists $u\in W_0^{1,\infty}(\omega;\R^3)$ with $\nabla u\in BV(\omega;\R^{3\times 3})$ such that
\begin{equation} \label{ineq:optimalitycylinder}
E_\epsilon^{(p)}(u;\omega) \lesssim   N^{-2p} +\epsilon  N   .
\end{equation}
 In particular, $E_\epsilon^{(p)}(u;\omega)  \lesssim \epsilon^{\frac{2p}{2p+1}}$. 
\end{proposition}
 We emphasize that Proposition \ref{prop:optimalfirst-order} particularly implies Proposition \ref{prop:Uppercylinder} by invoking the equivalence results from Section 4 in \cite{RTTZ25}. 

\begin{proof}
Given a cylindrical domain 
	$\omega = (0,1) \times B_1(0)$ and the two-dimensional branching construction $\tilde u\colon R \to \R^2$ from Lemma~\ref{lem:firstorder2d}(ii) with $R =  (0,1) \times (0,1/2) $,
let us parametrize $B_1(0)$ via the radius $r$ and angle $\phi$, where $B_1(0) \cap \{\phi = 0\} \subset \R \times \{0\}.$
We then define $u\colon  (0,1) \times B_1(0) \to \R^3$ by 
\begin{align}\label{eq:definitionin3dpolar}
		\begin{pmatrix}
	u_1 (x) \\
	u_2 (x) \\
	u_3 (x) \\
\end{pmatrix} \defas e_1
	\tilde u_1 (x_1, r) + \begin{pmatrix}
	0 \\
	\cos(\phi) \\
	\sin(\phi) \\
	\end{pmatrix} \tilde u_2(x_1,r) .
\end{align}  
This function attains zero Dirichlet boundary conditions due to the fact that $\tilde u$ satisfies zero Dirichlet boundary conditions. Notice that a cut-off as in \eqref{def:constcut-off} is \emph{not} needed.
 The identities 
\begin{align*}
	\frac{\partial r}{\partial x_2} = \cos(\phi), \quad \frac{\partial r}{\partial x_3} = \sin(\phi), \quad \frac{\partial \phi}{\partial x_2} = -\frac{ \sin (\phi)}{r}, \quad \text{and} \quad \frac{\partial \phi}{\partial x_3} = \frac{ \cos (\phi)}{r},
\end{align*}
imply that the gradient in cylindrical coordinates reads as
\begin{align*}
	\nabla u_i = \begin{pmatrix}
		\partial_1 u_i & \frac{\partial u_i}{\partial r} \cos(\phi) - r^{-1}\frac{\partial u_i}{\partial \phi} \sin(\phi) & \frac{\partial u_i}{\partial r} \sin(\phi) + r^{-1}\frac{\partial u_i}{\partial \phi} \cos(\phi)
	\end{pmatrix}   \text{ for } i = 1,2,3.
\end{align*}
Thus, the definition in \eqref{eq:definitionin3dpolar} leads to 
\begin{align}\label{eq:gradientpolar}
	\nabla u = \begin{pmatrix}
		\partial_1 \tilde u_1 & \partial_2 \tilde u_1 \cos(\phi) & \partial_2 \tilde u_1 \sin(\phi)    \\
				\partial_1 \tilde u_2 \cos(\phi) & \partial_2 \tilde u_2 \cos^2(\phi)  + r^{-1} \tilde u_2 \sin^2(\phi) &  (\partial_2 \tilde u_2 - r^{-1} \tilde u_2 ) \cos(\phi)\sin(\phi)  \\
						\partial_1 \tilde u_2 \sin(\phi) & (\partial_2 \tilde u_2 - r^{-1} \tilde u_2 ) \cos(\phi)\sin(\phi) & \partial_2 \tilde u_2 \sin^2(\phi) + r^{-1} \tilde u_2 \cos^2(\phi)   \\
	\end{pmatrix}.
\end{align} 
Hence, the uniform bound of $\sin$ and $\cos$, a change of coordinates and \eqref{eq:boundsfunctionepsilon45} yield
\begin{align}\label{ineq:energyboundrot}
	E_\epsilon^{(p)}(u;\omega) 
	&\lesssim  E_\epsilon^{(p)}(\tilde u; R ) + \int_0^1  \int_0^{1/2} \vert  r^{-1} \tilde u_2 (x_1,r) \vert^p r  \di x_1 \di r .    
\end{align}
As $\tilde u_2(x_1,0) = 0$ for a.e.\ $x_1$, the fundamental theorem of calculus and Jensen's inequality imply that
\begin{align*} 
	r^{-p+1} \vert \tilde u_2 (x_1,r) \vert^p = r r^{-p} \left\vert \int_0^r  \partial_2 \tilde u_2 (x_1,s) \di s \right\vert^p \leq \int_0^r \vert \partial_2 \tilde u_2 (x_1,s) \vert^p \di s . 
\end{align*}
Then, \eqref{eq:boundsfunctionepsilon45} implies that the right-hand side of \eqref{ineq:energyboundrot} enjoys the expected scaling.
 Using the chain rule in \eqref{eq:gradientpolar}, the absolutely continuous part of $D^2 u$ can be controlled via
\begin{align*}
	\vert \nabla^2 u \vert &\lesssim r^{-1} ( \vert \partial_2 \tilde u_1 + \partial_1 \tilde u_2 \vert + \vert \partial_2 \tilde u_2 \vert )  + \vert \nabla^2 \tilde u \vert  + r^{-2} \vert \tilde u_2  \vert 
\end{align*}
a.e.\ in $\Omega$.
 By an approximation argument, see \cite[Theorem~3.9]{AFP00}, the chain rule lifts to the total variation and we obtain 
\begin{align}\label{ineq:chainruleBV}
	\Vert D^2 u \Vert_{TV(\omega)} &\leq C \left( \int_0^1 \int_0^{1/2} ( \vert \partial_2 \tilde u_1 + \partial_1 \tilde u_2 \vert + \vert \partial_2 \tilde u_2 \vert ) \di x_1 \di r \right. \notag \\
& \qquad  \left. + \int_0^1 \int_0^{1/2} r^{-1} \vert \tilde u_2 \vert \di x_1 \di r  + \Vert D^2 \tilde u \Vert_{TV((0,1) \times (0,1/2))} \right).
\end{align}
Given $0\leq s \leq r \leq 1/2$, indicator functions satisfy the identity $\chi_{[0,r]} (s) = \chi_{[s,1/2]} (r)$.
Thus, Fubini's theorem and Hölder's inequality with coefficients $q$ and $\tilde q$ imply that 
\begin{align*}
	  \int_0^1 \int_0^{1/2} r^{-1} \vert \tilde u_2 \vert \di r \di x_1   & \leq	\int_0^1 \int_0^{1/2}  \vert \partial_2 \tilde u_2 (x_1, s)  \vert \int_s^{1/2} r^{-1}   \di r \di s   \di x_1 \\ 
		 			& \leq	C \Vert \partial_2 \tilde u_2 \Vert_{L^1} + 	  \Vert \partial_2 \tilde u_2 \Vert_{L^q} \left(\int_0^{1/2} \vert\log(s) \vert^{\tilde q} \di s \right)^{1/\tilde q}  .
	 \end{align*}
As $\log^{\tilde q}$ is integrable on $(0,1/2)$ for any $\tilde q \in [1,+\infty)$, we can use \eqref{eq:boundsfunctionepsilon45} for any $q > p $, \eqref{ineq:chainruleBV}, and \eqref{eq:surftwodirectiontwowell} such that
\begin{align*}
	\Vert D^2 u \Vert_{TV(\omega)} \lesssim N^{-2q} + N \lesssim N.
\end{align*}
 \end{proof}

\subsection{Second order Dirichlet data in the case of \eqref{eq:choiceofmatrices}}
In this subsection, we provide a rigorous upper bound construction in three dimensions in the context of second order Dirichlet boundary conditions.
To this end, we rely on the three-dimensional ansatz in \cite[Section~6.2]{RT23} once again and combine it with Lemma~\ref{lem:2well3d}.
\begin{proposition}\label{prop:upperboundcons}
Consider the setting of Theorem~\ref{mainth:3d}, and the set $K = K_{2,12}$ in \eqref{eq:choiceofmatrices}. Let $\Phi(z) \defas \frac{\sqrt{2}}{2} (z_1 + z_2, z_2 - z_1)^T$.
Consider $F = \diag(1/2, 1/2, -1/2) \in K^{(2)} \setminus K^{(1)}$ and the domain $\Omega \defas (0,1) \times \Phi^{-1}((0,1)^2)$.
Let $p \in [1,3]$.
Then for every $\epsilon\in(0,1)$ there exist $u \in W^{1,\infty}(\Omega;\R^3)$ such that $u (x)=F x$ for every $ x \in \partial \Omega$ and
\begin{equation}
E_\epsilon^{(p)}(u ;\Omega)\lesssim \epsilon^\frac{2p}{3+2p}.
\end{equation}
\end{proposition}
 We note that, as above, Proposition \ref{prop:upperboundcons} particularly implies Theorem~\ref{thm:mainresultupper} by invoking the equivalence results from Section 4 in \cite{RTTZ25}.

\begin{proof}[Proof of Proposition~\ref{prop:upperboundcons}]

According to Lemma~\ref{lem:coordinatetransform}(i) it suffices to provide an upper bound construction for
the set $\hat K = S K S^{-1}$ for $S \in SO(3)$ with $ Sz \defas (z_1, \frac{\sqrt{2}}{2} (z_2+z_3), \frac{\sqrt{2}}{2} (z_3-z_2))^T$ on the domain
$\hat \Omega \defas (0,1)^3$, where the wells and boundary datum are given as
\begin{align*}
	\hat K = \left\{  \diag (0,0,0), \diag (1,0,0), \begin{pmatrix}
	\frac{1}{2} & 0 & 0 \\
	0 & 0 & -1 \\
	0 & -1 & 0 \\
	\end{pmatrix} \right\} \quad \text{ and } \quad \hat F =\begin{pmatrix}
	\frac{1}{2} & 0 & 0 \\
	0 & 0 & -\frac{1}{2} \\
	0 & -\frac{1}{2} & 0 \\
	\end{pmatrix}.
\end{align*}
We follow the strategy of \cite[Section~6.2]{RT23} in constructing an upper bound competitor and proceed as in the proof of Lemma~\ref{lem:2well3d}.

	 To this end, we first provide a first order branching construction $u^{(1)}$ such that the boundary conditions are attained.
Let $\tilde u_{N_1} \colon  (0,1)^2 \to \R^2$ be defined via Lemma~\ref{lem:firstorder2d}(i) (with $L = 1$ and $H = 1$). Here, the parameter $N_1= \lceil \frac{1}{r_1} \rceil$ for $0 <r_1 = r_1(\epsilon) < 1$ is yet to be fixed.
We define $u^{(1)}$ as   $u_1^{(1)}(z)  = \frac{1}{2} z_1 $ and $u_2^{(1)}$ and $u_3^{(1)}$ by
\begin{align}\label{eq:definitionfirstorder}
\begin{pmatrix}
	u^{(1)}_2 (z_1,z_2,z_3) \\
	u^{(1)}_3 (z_1,z_2,z_3) \\
\end{pmatrix}
= - \frac{1}{2} \left( \tilde u_{N_1} \big(  z_2   , \rho(z_1,z_3)   \big)  \right) - \frac{1}{2} \begin{pmatrix}
	z_3 \\
	z_2 \\
\end{pmatrix}, 
\end{align}
 where $ \rho(z_1,z_3) \defas \max \{ \vert z_1 - \frac{1}{2} \vert,   \vert z_3 - \frac{1}{2} \vert  \} + \frac{1}{2}$.
Due to the construction of $\rho$, the fact that $\tilde u_{N_1}$ satisfies zero Dirichlet boundary conditions, and the symmetry condition of $\tilde u_{N_1}$ recalled in the argument for Lemma \ref{lem:firstorder2d}, $u^{(1)}$ attains the boundary condition on every face of the cubic domain.
 Geometrically, as the two-dimensional function $\tilde u_{N_1}$ branches out towards the bottom and top of $(0,1)^2$, the three-dimensional
construction $u^{(1)}$ branches out into the $e_1$- and the $e_3$-direction, i.e., the branching pattern is orthogonal to the vector $e_2$. 
In view of \eqref{eq:definitionfirstorder}, a computation of the gradient of $u^{(1)}$ yields
\begin{align}\label{eq:gradientofu1}
	\nabla u^{(1)}(z) = \begin{cases}
		\hat F   - \frac{1}{2} \begin{pmatrix}
		0 & 0 & 0 \\
		\partial_2  (\tilde u_{N_1})_1 (z_2,z_1) & \partial_1  ( \tilde u_{N_1})_1 (z_2,z_1)   & 0 \\
		\partial_2   (\tilde u_{N_1})_2 (z_2,z_1) & \partial_1   ( \tilde u_{N_1})_2 (z_2,z_1) & 0 \\
		\end{pmatrix}   & \text{if } \rho(z_1,z_3) = \rho(z_1,z_1) ,\\
		\hat F   - \frac{1}{2} \begin{pmatrix}
		0 & 0 & 0 \\
		0 & \partial_1  (\tilde u_{N_1})_1 (z_2,z_3)   & \partial_2  (\tilde u_{N_1})_1 (z_2,z_3)  \\
		0 & \partial_1(  \tilde  u_{N_1})_2 (z_2,z_3) & \partial_2   ( \tilde u_{N_1})_2 (z_2,z_3) \\
		\end{pmatrix}  & \text{if } \rho(z_1,z_3) = \rho(z_3,z_3). \\
	\end{cases}
\end{align}
Notice that in both cases the derivative $ \partial_1(  \tilde  u_{N_1})_2$ appears in the same entry.
Using the second inequality in \eqref{eq:eltwodirectiontwowell}, we thus derive that 
\begin{align}\label{eq:elasticfirstorder}
	 \int_{(0,1)^3} \dist^p \left(\sym (\nabla u^{(1)}), \tilde K \right) \di x  \lesssim N_1^{-p} \lesssim r_1^p  \text{ for } \tilde K \defas \left\{ \begin{pmatrix}
	\frac{1}{2} & 0 & 0 \\
	0 & 0 & 0 \\
	0 & 0 & 0 \\
	\end{pmatrix} ,  \begin{pmatrix}
	\frac{1}{2} & 0 & 0 \\
	0 & 0 & -1 \\
	0 & -1 & 0 \\
	\end{pmatrix}  \right\} .
\end{align}
Decomposing $\Omega$ into $E_1 \cup E_2 \defas (\Omega \cap \{ \vert z_1  -1/2 \vert > \vert z_3 - 1/2 \vert \}) \cup (\Omega \cap \{ \vert z_1  -1/2 \vert \leq \vert z_3 - 1/2 \vert \})$, 
on these sets $u^{(1)}$ is independent of $z_3$ or $z_1$, respectively.
Hence, the representation in \eqref{eq:gradientofu1} implies that the total variation of $D^2 u^{(1)}$ is controlled by the variation of $D^2 \tilde u_{N_1}$ and the interface $ \{ \vert z_1 - 1/2 \vert = \vert z_3 - 1/2 \vert \}$.
More precisely, by symmetry of this decomposition, \eqref{eq:surftwodirectiontwowell} and $r_1(\eps) \to 0 $ we have
\begin{align}\label{eq:surfacefirstorder}
	 \Vert D^2 u^{(1)} \Vert_{\rm TV(\Omega)} & \lesssim 
	 	 \Vert D[ (\chi_{E_1} +  \chi_{E_2}) \nabla u^{(1)}]  \Vert_{\rm TV(\Omega)} \leq 2 	 \Vert D [ \chi_{E_1} \nabla (u^{(1)}) ] \Vert_{\rm TV(\Omega)} \notag \\
	 & \leq 2 \Vert D \chi_{E_1}  \Vert_{\rm TV(\Omega)} \Vert \nabla (u^{(1)}) \Vert_{L^\infty} + \Vert D^2 \tilde u_{N_1} \Vert_{TV((0,1)^2)} \notag \\
	 &  \lesssim    \mathcal{H}^2( \Omega \cap \{ \vert z_1  -1/2 \vert > \vert z_3 - 1/2 \vert \}  ) +  \Vert D^2 \tilde u_{N_1} \Vert_{TV((0,1)^2)} \notag \\
	& \lesssim r_1^{-1} .
\end{align}
Analogously to \eqref{eq:buildingblock} and \eqref{eq:blocklengthdef}, we next divide $\Omega$ into the cells 
$$
\Omega_{j,k}=  \big\{ k\ell_j\leq z_2 \leq (k+1)\ell_j , \, y_j \leq \rho(z_1,z_3) \leq y_{j+1} \big\},
$$
where $\ell_j=\frac{1}{2^j N_1}$, $ y_j =1-\frac{\theta^j}{2}$ with $\theta \in (2^{-2q^*/(2q^*-1)}, 2^{-1})  $ for some $q^* > p$ and for $j \in \{0, \dots, j_0\}$, $k \in \{0,\dots, N2^j-1\}$ with $j_0$ such that $\ell_j < h_j = y_{j+1} - y_j $ for $j \leq j_0$.   
 Moreover, as in  the argument for Lemma \ref{lem:firstorder2d} and visualized in Figure~\ref{fig:branchingpicture}, we further subdivide the regions $\Omega_{j,k}$ into four blocks $\omega^1_{j,k}, \dots, \omega^4_{j,k}$.  
In what follows, the deformation on the blocks $\omega^2_{j,k}$ and $\omega^4_{j,k} $ will not be modified, while the blocks $\omega^1_{j,k}$ and $\omega^3_{j,k} $ play a further role in our construction
 as  ${\sym(\nabla u^{(1)})}$ is close to $\diag(1/2,0,0) \notin \hat K $. A careful inspection of the argument for Lemma \ref{lem:firstorder2d} reveals that 
these are given by
 \begin{align*}
\omega^1_{j,k} &= \Big\{ (z_1,z_2,z_3) \in\Omega_{j,k} : 0<z_2 - k\ell_j<\frac{\ell_j}{4}\Big\} \qquad  \text{and} \\
\omega^3_{j,k} &= \Big\{  (z_1,z_2,z_3)  \in\Omega_{j,k} : \frac{\ell_j}{4}+\frac{\ell_j (\rho(z_1,z_3) - y_j )}{4h_j}<z_2- k\ell_j<\frac{\ell_j}{2}+\frac{\ell_j (\rho(z_1,z_3) - y_j )}{4h_j}\Big\}.
 \end{align*}
As in these blocks too much elastic energy is paid, we wish to modify $u^{(1)}$ on these blocks by an inner branching construction via specifically chosen functions $\beta_{ijk} \colon \omega^i_{j,k} \to \R^3 $ for $i \in\{ 1,3\}$ with zero boundary conditions. 
More precisely, on these blocks we set $u = u^{(1)} + \beta_{ijk}$. Recall the definition of $\tilde K$ in \eqref{eq:elasticfirstorder}, set $\bar K = \{ \diag \big(1/2,0,0 \big) , \diag \big(-1/2,0,0 \big)  \}$, and estimate
\begin{align}\label{eq:helpforelasticenergy}
	\dist(\sym(\nabla u), \hat K) \leq  \dist(\sym(\nabla u^{(1)}) , \tilde K  ) +  \dist(\sym (\nabla \beta_{ijk}), \bar K).
\end{align}
A similar estimate can be derived for the surface energy.

\begin{figure}[t]
\centering
\begin{tikzpicture}[thick, scale= 7]

\draw[->] (0,0) -- (1.05,0) node[right] {$z_1$};
\draw[->] (0,0) -- (0,1.05) node[above] {$z_3$};
 
\def\yzero{0.5}
\def\yone{0.761905}  
\def\ytwo{0.88662}  
\def\ythree{0.94604}  
\def\yfour{0.97429}

\node[below,font=\scriptsize] at (\yzero,0) {$ y_0 $};
\node[below,font=\scriptsize] at (\yone,0) {$ y_1 $};
\node[below,font=\scriptsize] at (\ytwo,0) {$ y_2 $};
\node[below,font=\scriptsize] at (\ythree,0) {$ y_3 $}; 

\node[left,font=\scriptsize] at (0,\yzero) {$ y_0 $};
\node[left,font=\scriptsize] at (0,\yone) {$ y_1 $};
\node[left,font=\scriptsize] at (0,\ytwo) {$ y_2 $};
\node[left,font=\scriptsize] at (0,\ythree) {$ y_3 $};

\fill[black!50] (\yone,0.5) rectangle (\ytwo,\ytwo);
\fill[black!50] (\yone,{1-\ytwo}) rectangle (\ytwo,0.5);
\fill[pattern=north east lines] (\yone,0.5) rectangle (\ytwo,\ytwo);
\fill[pattern=north east lines] (\yone,{1-\ytwo}) rectangle (\ytwo,0.5);
 
\fill[black!40] (\ytwo,0.5) rectangle (\ythree,\ythree);
\fill[black!40] (\ytwo,{1-\ythree}) rectangle (\ythree,0.5);
\fill[pattern=north east lines] (\ytwo,0.5) rectangle (\ythree,\ythree);
\fill[pattern=north east lines] (\ytwo,{1-\ythree}) rectangle (\ythree,0.5);
 
\fill[black!30] (\ythree,0.5) rectangle (\yfour,\yfour);
\fill[black!30] (\ythree,{1-\yfour}) rectangle (\yfour,0.5);
\fill[pattern=north east lines]  (\ythree,0.5) rectangle (\yfour,\yfour);
\fill[pattern=north east lines]  (\ythree,{1-\yfour}) rectangle (\yfour,0.5);

\fill[black!45] (1-\yone,\yone) rectangle (\yone,\ytwo);
\fill[pattern=north west lines] (1-\yone,\yone) rectangle (\yone,\ytwo);
 
\fill[black!35] (1-\ytwo,\ytwo) rectangle (\ytwo,\ythree);
\fill[pattern=north west lines]  (1-\ytwo,\ytwo) rectangle (\ytwo,\ythree);
 
\fill[black!25] (1-\ythree,\ythree) rectangle (\ythree,\yfour);
\fill[pattern=north west lines]  (1-\ythree,\ythree) rectangle (\ythree,\yfour);

\fill[black!50] (1-\ytwo,0.5) rectangle (1-\yone,\ytwo);
\fill[black!50] (1-\ytwo,{1-\ytwo}) rectangle (1-\yone,0.5);
\fill[pattern=north east lines]  (1-\ytwo,0.5) rectangle (1-\yone,\ytwo);
\fill[pattern=north east lines]  (1-\ytwo,{1-\ytwo}) rectangle (1-\yone,0.5);
 
\fill[black!40] (1-\ythree,0.5) rectangle (1-\ytwo,\ythree);
\fill[black!40] (1-\ythree,{1-\ythree}) rectangle (1-\ytwo,0.5);
\fill[pattern=north east lines]  (1-\ythree,0.5) rectangle (1-\ytwo,\ythree);
\fill[pattern=north east lines]  (1-\ythree,{1-\ythree}) rectangle (1-\ytwo,0.5);
 
\fill[black!30] (1-\yfour,0.5) rectangle (1-\ythree,\yfour);
\fill[black!30] (1-\yfour,{1-\yfour}) rectangle (1-\ythree,0.5);
\fill[pattern=north east lines]  (1-\yfour,0.5) rectangle (1-\ythree,\yfour);
\fill[pattern=north east lines]  (1-\yfour,{1-\yfour}) rectangle (1-\ythree,0.5);

\fill[black!45] (1-\yone,1-\ytwo) rectangle (\yone,1-\yone);
\fill[pattern=north west lines] (1-\yone,1-\ytwo) rectangle (\yone,1-\yone);
 
\fill[black!35] (1-\ytwo,1-\ythree) rectangle (\ytwo,1-\ytwo);
\fill[pattern=north west lines] (1-\ytwo,1-\ythree) rectangle (\ytwo,1-\ytwo);
 
\fill[black!25] (1-\ythree,1-\yfour) rectangle (\ythree,1-\ythree);
\fill[pattern=north west lines] (1-\ythree,1-\yfour) rectangle (\ythree,1-\ythree);

\fill[black!60] (1-\yone,1-\yone) rectangle (\yone,\yone);
\fill[pattern= crosshatch] (1-\yone,1-\yone) rectangle (\yone,\yone);

\end{tikzpicture}
\caption{Projection of the domain onto the $z_1z_3$-plane -- schematic definition of the rectangular ``annuli''
$\{ y_j \leq \rho(z_1,z_3) \leq y_{j+1} \}$ for $j \in\{ 0,1,2,3\}$ and their decomposition into the 
sets $A_{m,j,k}$ for a fixed $k$ and $m = 1, \dots, 4$. }
\label{fig:decomposition}
\end{figure}
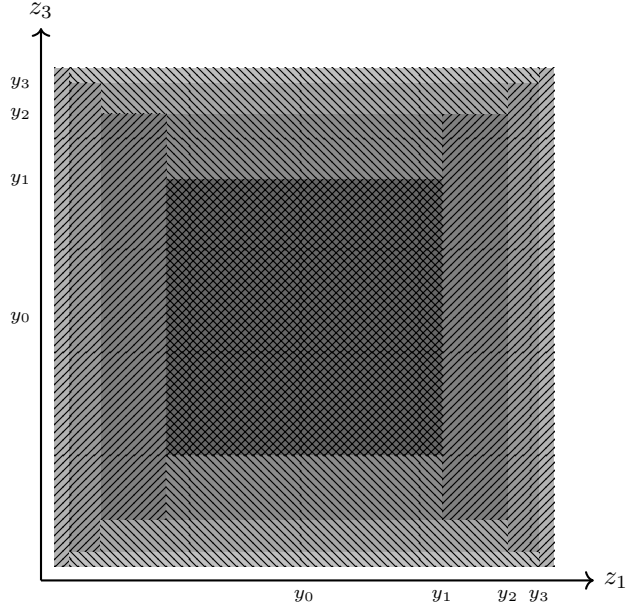

\emph{Definition of $\beta_{1jk}$ on $ \omega_{j,k}^1 $ by an inner branching construction:}  
$\omega^1_{j,k}$, $j > 1$, can be divided into four disjoint (but touching) cuboids via
\begin{align*}
	A_{1,j,k} &\defas \{ z_1 \in [ y_j , y_{j+1} ] , \  \vert z_3 - 1 /2 \vert + 1/2 \in [1/2, y_{j+1} ] \}    , \\
		A_{2,j,k} &\defas \{ \vert z_1 - 1/2 \vert + 1/2 \in [1/2, y_j ] , \   z_3  \in [ y_j , y_{j+1}  ] \} , \\  
			A_{3,j,k} &\defas \{ 1- z_1  \in [ y_j , y_{j+1} ] , \  \vert z_3 - 1/2 \vert + 1/2 \in [1/2, y_{j+1} ] \} ,    \\
		A_{4,j,k} &\defas \{ \vert z_1 - 1/2 \vert + 1/2 \in [1/2, y_j ] , \ 1 -  z_3  \in [ y_j , y_{j+1} ] \}.
\end{align*}
Notice that $\omega_{0,k}^1$ is already a cuboid and we set $A_{1,0,k} = A_{2,0,k} = A_{3,0,k} = A_{4,0,k} \defas \omega_{0,k}^1$.
For a schematic visualization we refer to Figure~\ref{fig:decomposition}. 

Let us first define $\beta_{ijk}$ on $ \omega_{j,k}^1 \cap A_{i,j,k}$ for $i \in \{ 1,\dots,4\}$ and $j \geq 0$, and recall the relations $\ell_j \sim \frac{r_1}{2^j }$ and  $h_j \sim \theta^j \sim y_{j+1} - y_j $, see \eqref{eq:blocklengthdef}.

 In particular, for $i \in \{ 2,4\}$ we apply Lemma~\ref{lem:2well3d}
for $L = 2 y_j - 1$, $H_1 =  \ell_j/4$, $H_2 = y_{j+1} - y_j $, and $N = N_j = \lceil 2^j  \theta^{-j} \varphi^j  r_2^{-1} \rceil $  (up to a suitable linear transformation), where $\varphi \in (\theta,\frac{1}{2})$.
Notice that Lemma~\ref{lem:2well3d} is applicable for $r_2 \lesssim r_1$ due to $N \geq 2^j r_2^{-1} \geq 2^j r_1^{-1} \geq  4 L H_1^{-1} $ and $H_1 = \ell_j/4 < h_j /4 = H_2/4 < H_2$. 

Thus,  \eqref{eq:3donedirectiontwowell} gives    
\begin{align}\label{eq:secondorder1}
  		& \quad \   \int_{\omega_{j,k}^1 \cap A_{i,j,k}} \dist^p(\sym (\nabla \beta_{ijk}), \tilde K)  \di x + \epsilon   \Vert D^2 \beta_{ijk} \Vert_{TV(\omega_{j,k}^1 \cap A_{i,j,k})} \notag \\
		   & \lesssim (\theta / \varphi)^{j (p+1)}  2^{-2j}  r_2^{p+1} r_1^{-p+1} +  (\theta / \varphi)^{2jp }  \theta^j 2^{-j} r_2^{2p} r_1^{1-2p}  + \epsilon \varphi^j  r_2^{-1} r_1 + \epsilon 2^{-j} \varphi^{-j} \theta^j  r_2  . 
\end{align}
 Below in \eqref{eq:collecteverything} we will see that the parameters $\theta$ and $\varphi$ have to be chosen such that the last line of \eqref{eq:secondorder1} multiplied with $2^j$ generates a converging series when summed up among all $j \in \N$.

Similarly, for $i \in\{ 1,3\}$, 
 we use Lemma~\ref{lem:2well3d}
for $L =  y_{j+1} - y_j $, $H_1 =  \ell_j/4$, $H_2 = 2 y_j - 1 $, and $N = N_j = \lceil (2\theta)^j r_2^{-1} \rceil$ (up to a suitable linear transformation).
 As before we have $N \geq (2\theta)^j r_2^{-1} \geq \theta^j 2^j r_1^{-1} \geq  4 L H_1^{-1} $ due to  $r_2 \lesssim r_1$. 
Hence, \eqref{eq:3donedirectiontwowell} implies that 
\begin{align}\label{eq:secondorder2}
  		& \quad \   \int_{\omega_{j,k}^1 \cap A_{i,j,k}} \dist^2(\sym (\nabla \beta_{ijk}), \tilde K)  \di x + \epsilon   \Vert D^2 \beta_{ijk} \Vert_{TV(\omega_{j,k}^1 \cap A_{i,j,k})}   \notag \\
		 &   \lesssim    2^{-2j} \theta^{j}        r_2^{p+1}   r_1^{-p+1 }   +  \theta^j   2^{-j}     r_2^{2p} r_1^{1-2p}   +\epsilon  \theta^j  r_2^{-1} r_1 + \epsilon  \theta^{j}  2^{-j} r_2    . 
\end{align}
Notice that the above construction also works for $j = 0$.

\emph{Definition of $\beta_{3jk}$ on $ \omega_{j,k}^3 $:} 
The construction can be obtained from the one above by a change of variables.
Notice that we have $f \colon  \omega_{j,k}^1 \to  \omega_{j,k}^3$
\begin{align*}
	f(z) = \begin{pmatrix}
		z_1 \\
		z_2 + \frac{\ell_j}{4h_j}  (\rho(z_1,z_3) - y_j ) \\
		z_3
	\end{pmatrix}
	 + \begin{pmatrix}
		0 \\
		\frac{\ell_j}{4} \\
		0 \\
	\end{pmatrix}.
\end{align*} 
Thus, we set
\begin{align*} 
\beta_{3jk}(z)=\beta_{1jk}\Big(z_1,z_2-\frac{\ell_j ( \rho(z_1,z_3) - y_j )}{4h_j}-\frac{\ell_j}{4},z_3\Big),
\end{align*}
which attains zero boundary conditions by construction.
A computation of the gradient reveals
\begin{align*} 
	\nabla \beta_{3jk}(z) = \begin{cases} 
		\nabla \beta_{1jk}(z_1,z_2-\frac{\ell_j ( \rho(z_1,z_3) - y_j )}{4h_j}-\frac{\ell_j}{4},z_3 ) \begin{pmatrix}
			1 & 0 & 0 \\
			-\frac{\ell_j}{4h_j} & 1 & 0 \\
			0 & 0 & 1\\
		\end{pmatrix} 
		& \text{if } \rho(z_1,z_3) = \rho(z_1,z_1) ,\\
			\nabla \beta_{1jk}(z_1,z_2-\frac{\ell_j ( \rho(z_1,z_3) - y_j )}{4h_j}-\frac{\ell_j}{4},z_3 ) \begin{pmatrix}
			1 & 0 & 0 \\
			0& 1 & -\frac{\ell_j}{4h_j}  \\
			0 & 0 & 1\\
		\end{pmatrix}  & \text{if } \rho(z_1,z_3) = \rho(z_3,z_3). \\
	\end{cases}
\end{align*}
 The estimate \eqref{ineq:uniformboundpartial2} allows us to control $\Vert\partial_2 \beta_{1jk}\Vert_{L^p}$ and we obtain due to the previous choices of $L$, $H_1$, $H_2$ and $N$ 
 and the transformation theorem 
\begin{align}\label{eq:affinetrans2}
 \int_{\omega_{j,k}^{(3)}}\dist^p( \sym (\nabla \beta_{3jk} ), \hat K ) \di x & \lesssim\int_{\omega_{j,k}^{(1)}}\dist^p\Big( \sym (\nabla \beta_{1jk}) , \hat K \Big) \di x \notag \\
& \qquad +  r_1 r_2^p \left(   2^{-j(p+1)} \theta^{j (1-p)} + 2^{-j (p+1)} \theta^j \varphi^{-jp}   \right). 
\end{align}
Similarly, the total variation can be controlled.
Eventually, using \eqref{eq:elasticfirstorder}, \eqref{eq:surfacefirstorder}, \eqref{eq:helpforelasticenergy}, \eqref{eq:secondorder1}, \eqref{eq:secondorder2}, $\frac{1}{2^{p/(p-1)}} < 2^{-2q^*/(2q^*-1)} <\theta < \varphi< \frac{1}{2}$ and \eqref{eq:affinetrans2} we sum over $j$ and $k$ and derive that
\begin{align}\label{eq:collecteverything}
& \quad \	E_\epsilon^{(p)}(u;\Omega) \notag  \\
& \lesssim  r_1^p + \epsilon r_1^{-1} + \sum_{j=0}^{j_0+1} \sum_{k= 0}^{N_1 2^j -1}  \Big(   2^{-2j}  r_2^{p+1} r_1^{-p+1} +   \theta^j 2^{-j} r_2^{2p} r_1^{1-2p} +  \epsilon \varphi^j  r_2^{-1} r_1 + \epsilon 2^{-j} \varphi^{-j} \theta^j  r_2  \notag \\ 
&  \qquad  \qquad \qquad \qquad \qquad \qquad \qquad \qquad \qquad \qquad 
 +  r_1 r_2^p \big(   2^{-j(p+1)} \theta^{j (1-p)} + 2^{-j (p+1)} \theta^j \varphi^{-jp}   \big)   \Big) \notag \\
& \lesssim  r_1^p + \epsilon r_1^{-1} + \sum_{j=0}^{\infty}  \Big(   2^{-j}  r_2^{p+1} r_1^{-p} +   \theta^j   r_2^{2p} r_1^{-2p}  + \epsilon \varphi^j 2^j r_2^{-1} + \epsilon   \varphi^{-j} \theta^j  r_2 r_1^{-1}  \notag \\
& \qquad \qquad \qquad \qquad \qquad \qquad  \qquad  +  r_2^p \big(   2^{-jp} \theta^{j (1-p)} + 2^{-j p} \varphi^{j (1-p)}  \big)    \Big)  \notag \\
& \lesssim  r_1^p +  r_2^{p+1} r_1^{-p} +   r_2^{2p} r_1^{-2p}   + r_2^p     + \epsilon (r_1^{-1} + r_2^{-1} + r_2 r_1^{-1}  ) .
\end{align} 
Here, we additionally have used that the cut-off energy is sufficiently small as remarked below \eqref{eq:cutoffbuildingblock}.
As we choose $r_2$ as a smaller length scale than $r_1$, i.e., $r_2 \lesssim r_1 $, it suffices to optimize
\begin{align}\label{ineq:optimizationproblemend}
	 r_1^p +  r_2^{p+1} r_1^{-p} +   r_2^{2p} r_1^{-2p}         + \epsilon   r_2^{-1}
\end{align}  
in terms of $\epsilon$. Notice that the second term appears due to the cut-off in Lemma~\ref{lem:2well3d}. We first optimize by neglecting this term and
choose $r_1^p \sim r_2^{2p} r_1^{-2p}$, which leads to
  $ r_1 \sim  r_2^{2/3}$. A substitution of this scaling into the previous equation leads to 
\begin{align}\label{eq:finalopti}
	E_\epsilon^{(p)}(u;\Omega)  \lesssim r_2^{2p/3} +  r_2^{p/3+1}      + \epsilon   r_2^{-1} .
\end{align}
Finally, as $p\in [1,3]$ the choice $r_2 \sim \epsilon^{3/(3+2p)}$ concludes the proof.
 \end{proof}
 For $p>3$, term $ r_2^{p+1} r_1^{-p}$ in \eqref{ineq:optimizationproblemend} which results from the cut-off dominates.
Here, the asymptotic equivalence $r_1^p \sim  r_2^{p+1} r_1^{-p} \sim \epsilon r_2^{-1} $ leads to the scaling
$\sim \epsilon^{(p+1)/(p+3)}$ which coincides with the bound in \eqref{prop:upperboundcons} for $p = 3$, but is worse for $p>3$. If one finds a bound on a cube without the suboptimal cut-off part which for instance does not appear in Proposition~\ref{prop:optimalfirst-order}, one can also improve the scaling for $p>3$.

\subsection{Heuristic argument for the upper scaling behaviour of higher order data}\label{sec:heuristics}
We conclude our discussion of upper bounds by providing a heuristic argument for the sharpness of the deduced lower bounds in Theorems~\ref{mainth:md} and \ref{mainth:mdperiodic} for general $\ell, m$. We emphasize that our argument does \emph{not} provide a proof in the form of a matching upper bound construction, but only serves as a rough heuristic. As the above arguments in the previous subsections show, we expect that making this rigorous requires major additional arguments and difficulties. We plan to return to this in future work.

Let $F \in K^{(\ell)}_m \setminus K^{(\ell-1)}_m$ for $\ell \in \{1,\dots,m\}$ be a datum of $\ell$-th order.
Following the construction in Proposition~\ref{prop:upperboundcons},  in the case of a Dirichlet boundary datum we expect that upper bound constructions require $\ell$ iterated branching constructions on $\ell$ different length scales $r_\ell $.
 In our heuristic argument, as a crucial simplification we \emph{ignore} the nontrivial cut-off in Lemma~\ref{lem:2well3d} in the following
  and perform the computations by assuming that the energy scaling from the two-dimensional setting in Lemma~\ref{lem:firstorder2d} can be transferred also to higher dimensions.
  The magnitude of the elastic energy generated by a branching construction on the $i$-th length scale strongly depends
  on the amount of compatibility directions.
  
  Extrapolating from the previous section, we expect that an elastic contribution of order $ \sim r_i^p r_{i-1}^{-p}$ should be generated
  if two compatibility directions exist, while one would hope that this could be refined to $\sim r_i^{2p} r_{i-1}^{-2p}$ in the case of a single direction, see the third term in \eqref{ineq:optimizationproblemend}. 
  Moreover, it is expected that it suffices to control the surface energy by the finest length scale which results in the estimate 
  $ E_{\rm surf} \lesssim r_\ell^{-1} $.
  	This heuristically leads to the estimate
\begin{align}\label{eq:optim}
	  \Energyinf{m+1}{p}{\dir} \lesssim \eps r_\ell^{-1} + \sum_{\substack{i=1 , \dots , \ell : \\ f(i) = 1} }  \left(\frac{r_i}{r_{i-1}} \right)^{2p} + \sum_{\substack{i=1 , \dots , \ell : \\ f(i) = 2} }  \left(\frac{r_i}{r_{i-1}}\right)^{p},
\end{align}
where the length scales $r_i = r_i(\epsilon)$, $i \in\{ 0, \dots, \ell\}$, are functions such that
	\begin{align*}
	\lim\limits_{\epsilon \to 0 }\frac{r_{i+1}}{r_i} = 0 \quad \text{ and } \quad r_0 \sim 1.
	\end{align*}
However, we stress once again that these are heuristic computations.

The case of periodic boundary conditions can be tackled similarly. We expect that the essential difference is that one does not need to resort to a branching construction for the lowest order. In fact, as in the periodic setting no fixed Dirichlet boundary conditions have to be attained, a (single) simple laminate (without fine-scale oscillation) always suffices for the coarsest, outermost construction. Hence, on this scale only surface energy is used. The Lipschitz regularity of the deformation, however, requires matching boundary conditions with the constructed first order laminate by turning to the next order laminates and by replacing the not yet compatible strain values with higher order laminates. 
To this end, one expects to use nested branching constructions on $\ell -1$ length scales, leading to an overall bound of the form 
\begin{align}\label{eq:optimper}
	  \Energyinf{m+1}{p}{\per} \lesssim \eps r_{\ell-1}^{-1} + \sum_{\substack{i=1 , \dots , \ell -1 : \\ f(i) = 1} }  \left(\frac{r_i}{r_{i-1}} \right)^{2p} + \sum_{\substack{i=1 , \dots , \ell -1 : \\ f(i) = 2} }  \left(\frac{r_i}{r_{i-1}}\right)^{p},
\end{align}
where the sums on the right-hand side vanish for $\ell = 1$.
\begin{lemma}[Heuristic scaling of the upper bound]\label{lem:heuristic}
	For an explicit choice of sequences $r_i$, $i \in \{ 0, \dots, \ell\}$, in \eqref{eq:optim} and \eqref{eq:optimper}, respectively, it holds that
	\begin{align}
	 \eps r_\ell^{-1} + \sum_{\substack{i=1 , \dots , \ell : \\ f(i) = 1} }  \left(\frac{r_i}{r_{i-1}} \right)^{2p} + \sum_{\substack{i=1 , \dots , \ell : \\ f(i) = 2} }  \left(\frac{r_i}{r_{i-1}}\right)^{p} & \lesssim \epsilon^{2p / (2p+   2(\ell-k_{\ell}) + k_{\ell} )} \quad \text{and}    \label{eq:finalstatement}  \\
	\eps r_{\ell-1}^{-1} + \sum_{\substack{i=1 , \dots , \ell -1 : \\ f(i) = 1} }  \left(\frac{r_i}{r_{i-1}} \right)^{2p} + \sum_{\substack{i=1 , \dots , \ell -1 : \\ f(i) = 2} }  \left(\frac{r_i}{r_{i-1}}\right)^{p} & \lesssim \epsilon^{2p / (2p+   2((\ell-1)- k_{\ell-1}) + k_{\ell-1} )} \label{eq:finalstatementper} ,
\end{align}
where $k_0 = 0$ and  $k_\ell$ is defined in \eqref{eq:kr}.
\end{lemma}

\begin{proof} We first show \eqref{eq:finalstatement}.
	To prove the desired upper scaling bound, we consider functions $r_i = r_i(\epsilon)$, $i \in \{ 0, \dots, \ell \}$,
	such that every term on the left-hand side of
	\eqref{eq:finalstatement} is asymptotically equivalent. We first analyze the relation between $r_\ell$ and $r_1$ and show the asymptotic formula
	\begin{align}\label{eq:inductionstatement}
		r_\ell            \sim r_1^{ (\delta_{2f(1)} + 2\delta_{1f(1)} )   (  \ell - \frac{k_\ell}{2})} 
	\end{align}
  by means of an induction over $\ell$.
The base case $\ell = 1$ follows from the fact that $k_\ell = 0$ if $f(1) = 2$ and  $k_\ell = 1$ if $f(1) = 1$.
For convenience, we prove the formula for $\ell = 2$ explicitly.
The asymptotic equivalence of the sums on the left-hand side of
	\eqref{eq:finalstatement} and $r_0 \sim 1$ yields
\begin{align*}
	\begin{cases}
		r_2 \sim  r_1^2  & \text{if } f(2) = f(1), \\
			r_2 \sim  r_1^{3/2}     & \text{if } f(2) = 1, \  f(1) = 2 , \\
			 r_2  \sim   r_1^3   & \text{if } f(2) = 2, \  f(1) = 1 , \\
	\end{cases}
\end{align*}
which is the desired formula for $\ell = 2$.
In order to show the induction step, we assume that the formula in \eqref{eq:inductionstatement} holds for $\ell - 1$ with
$\ell \geq 2$. 
The asymptotic equivalence of the sums particularly yields
\begin{align*}
	\left(\frac{r_\ell}{r_{\ell-1}} \right)^{4 \delta_{1f(\ell)}+ 2 \delta_{2 f(\ell)}} \sim r_1^{4 \delta_{1f(1)}+ 2 \delta_{2 f(1)}}  .
\end{align*}
Hence, by using the formula in \eqref{eq:inductionstatement}  for $\ell - 1$, we find that
\begin{align*}
	r_\ell^{4 \delta_{1f(\ell)}+ 2 \delta_{2 f(\ell)}} &\sim r_{\ell-1}^{4 \delta_{1f(\ell)}+ 2 \delta_{2 f(\ell)}} r_1^{4 \delta_{1f(1)}+ 2 \delta_{2 f(1)}}  \\
													&\sim r_1^{  \left( (4 \delta_{1f(\ell)}+ 2 \delta_{2 f(\ell)}) (\delta_{2f(1)} + 2\delta_{1f(1)} )   (  (\ell-1) - \frac{k_{\ell-1}}{2}) \right) + 2 (2 \delta_{1f(1)}+  \delta_{2 f(1)} ) }  \\
													&\sim r_1^{ (\delta_{2f(1)} + 2\delta_{1f(1)} ) \left( \left( (4 \delta_{1f(\ell)}+ 2 \delta_{2 f(\ell)})    (  (\ell-1) - \frac{k_{\ell-1}}{2}) \right) + 2 \right) }  .
\end{align*}
By taking the $(4 \delta_{1f(\ell)}+ 2 \delta_{2 f(\ell)})$-th root and the relation $k_\ell = k_{\ell -1} + \delta_{1f(\ell)}$ we obtain
\begin{align*}
	r_\ell 	&\sim r_1^{ (\delta_{2f(1)} + 2\delta_{1f(1)} ) \left(  (  (\ell-1) - \frac{k_{\ell-1}}{2})  + \frac{1}{2} \delta_{1f(\ell)} +  \delta_{2 f(\ell)} \right) }  \sim  r_1^{ (\delta_{2f(1)} + 2\delta_{1f(1)} )   (  \ell - \frac{k_\ell}{2})}   ,
\end{align*}
which is precisely \eqref{eq:inductionstatement}. This formula and the asymptotic equivalence on the left-hand side of
	\eqref{eq:finalstatement} imply that
	\begin{align*}
		 \eps r_\ell^{-1} + \sum_{\substack{i=1 , \dots , \ell : \\ f(i) = 1} }  \left(\frac{r_i}{r_{i-1}} \right)^{2p} + \sum_{\substack{i=1 , \dots , \ell : \\ f(i) = 2} }  \left(\frac{r_i}{r_{i-1}}\right)^{p} &\lesssim \eps  r_1^{ - (\delta_{2f(1)} + 2\delta_{1f(1)} )   (  \ell - \frac{k_\ell}{2})}  + r_1^{p (\delta_{2f(1)} + 2\delta_{1f(1)} )} \\
	 &=  \eps  z^{ -  (  \ell - \frac{k_\ell}{2})}  + z^p,
\end{align*} 
where the last equality follows by a substitution of $r_1^{ (\delta_{2f(1)} + 2\delta_{1f(1)} )}$ by $z$.
Choosing
$z \sim \eps^{   2 / (2p+ ( 2\ell - k_\ell))   }  $, we conclude the proof of \eqref{eq:finalstatement}. The estimate 
\eqref{eq:finalstatementper} can be deduced by the same optimization as the left-hand side of \eqref{eq:finalstatement}. Similarly, \eqref{eq:finalstatementper} is shifted by one index, where we additionally note that the case $\ell = 1$ is trivial.  
\end{proof}

\appendix
\section{Proofs of auxiliary results}\label{sec:Appendix} 

\subsection{Proof of Lemma~\ref{lem:compatibility}}
We provide a short proof of Lemma~\ref{lem:compatibility} if $A,B \in \R^{d \times d}_\diag$.
\begin{proof}[Proof of Lemma~\ref{lem:compatibility} for diagonal strains]
 We first show that compatibility implies the claimed eigenvalue condition.
	By the subadditivity of the rank, we have ${\rm rank} (A-B) \leq 2$.
	Thus, as $A$ and $B$ both are diagonal, at most two entries of $A-B$ are non-zero. Due to \eqref{eq:def:compatible} we can write these entries as $a_1 b_1 = \lambda_1 $ and $a_2 b_2 = \lambda_2$ with $a_2 b_1 + b_2 a_1 = 0$ for suitable $a_1, a_2, b_1, b_2 \in \R$. If $a_1 \neq 0 \neq a_2$, this implies that $\lambda_1 a_2 /a_1 + \lambda_2 a_1 / a_2 = 0$. This can only be guaranteed if $\lambda_1$ and $\lambda_2$ have the opposite sign.
	If $a_1=0$ or $a_2=0$, then $\lambda_1 =0$ or $\lambda_2 = 0$, respectively. 

	To show the converse statement, let $\lambda_1 > 0$ be the $i$-th eigenvalue of $A-B$ which is the only one not equal to zero. Then, we can write $a = (0,\dots, 0, \sqrt{\lambda_1}, 0, \dots, 0 ) ^T$ and $A-B = a \odot a$. In a similar fashion, we write $a \odot -a$ in the case $\lambda_1<0$ for a suitable $a$. If $\lambda_1<0$ and $\lambda_2 >0$, we have
	$a = (0, \dots, 1,0, \dots, 0 ,\sqrt{- \lambda_1 \lambda_2} / \lambda_1, \dots ,0)^T$ and $b = (0, \dots, \lambda_1, 0, \dots, 0, -  \sqrt{ - \lambda_1 \lambda_2}, \dots, 0)^T$, giving   $A-B = a \odot b$. By interchanging $\lambda_1$ and $\lambda_2$, we conclude the proof.
\end{proof}

\subsection{Proof of Lemma~\ref{lem:Fourierchar}}
The proof of Lemma~\ref{lem:Fourierchar} is strongly inspired by \cite[Lemma~3.1]{CO1} and \cite[Lemma~4.1]{KKO13}.
\begin{proof}[Proof of Lemma~\ref{lem:Fourierchar}]
  Our goal is to reformulate the elastic energy as a Fourier series.   
 Let us first consider the case $u \in \mathcal{A}_F^{\dir} $. 
  To this end, we subtract the boundary conditions of $u$, see \eqref{def:boundarycond}, such that we can view the function as a periodic function on $\mathbb{T}^d$.
  More precisely, we set $
	  v(x) = u(x) - Fx - b$ and we have $v \in W^{1,p}_0(\Omega;\R^d) \subset H_0^1(\Omega;\R^d)  \subset H^1(\mathbb{T}^d;\R^d)$. 
 We denote the Fourier transforms of $v$ and $\tilde \chi$ by  $\hat v $ and $\F \tilde\chi $, respectively. We again recall that $\F \tilde{\chi}(\xi) = \hat{\chi}(\xi)$ for $\xi \neq 0$ and  $\F \tilde{\chi}(0) = \hat{\chi}(0)-F$. 
In view of \eqref{eq:Dirichletsetandperiodic}, Hölder's inequality,
  Plancherel's identity (see e.g.\ \cite[Proposition~3.2.7]{Grafakos}), we obtain
\begin{align}\label{eq:writeinFour}
  \Eelinf{m+1}{p}^{2/p}    &\geq \inf\limits_{{v \in H^1(\mathbb{T}^d;\R^d) } } \int_{\mathbb{T}^d} \vert \sym (\nabla v) - \tilde \chi \vert^2 \, {\rm d}x \notag\\
  &= \inf\limits_{v \in H^1(\mathbb{T}^d;\R^d)    } \sum\limits_{ \xi \in \mathbb{Z}^d \setminus \{0 \}}  \vert 2 \pi \hat v \odot  i \xi - \hat \chi \vert^2 +  \vert  \hat \chi(0) -F \vert^2 .
\end{align}
The next goal is to minimize every term on the right-hand side pointwisely in $\hat v(\xi)$.   
In view of the Pythagorean theorem, we find that
\begin{align}\label{eq:complexexpansion}
	\vert 2 \pi  \hat v \odot  i \xi - \hat \chi \vert^2 &= \vert i \left( 2 \pi  {\rm Re} [\hat v ] \odot   \xi - {\rm Im} [ \hat \chi ] \right) - 2 \pi  {\rm Im} [\hat v ]  \odot  \xi - {\rm Re} [ \hat \chi ] \vert^2 \notag \\ 
	&= \vert  2 \pi  {\rm Re} [\hat v ] \odot   \xi - {\rm Im} [ \hat \chi ] \vert^2 + \vert 2 \pi  {\rm Im} [\hat v ]  \odot  \xi + {\rm Re} [ \hat \chi ] \vert^2,
\end{align}
where ${\rm Re} [ \cdot ] $ and $ {\rm Im} [\cdot] $ denote the real and imaginary part of a complex number, respectively.
Hence, we can view the real and imaginary part of $\hat v$ as independent variables in the minimization problem above. 
 To minimize the first term on the right-hand side of \eqref{eq:complexexpansion},
the strict convexity of the Euclidean norm implies  that it is equivalent to compute the unique value ${\rm Re}[\hat v_* ] $ with vanishing gradient, i.e.,
\begin{align}\label{eq:initialeulerlag}
	0 =    \left( 2 \pi {\rm Re}[\hat v_* ] \odot  \xi - {\rm Im} [ \hat \chi ] \right):   e_j \odot  \xi     =  \left( 2 \pi {\rm Re}[\hat v_* ] \odot  \xi - {\rm Im} [ \hat \chi ]  \right) \xi \cdot e_j    
\end{align}
 for all   $j = 1, \dots, d$ , where $e_j$ denotes the $j$-th standard unit vector in $\R^d$.
Given vectors $a,b \in \R^d$, we have $ 2 (a \odot b ) b = \vert b \vert^2 a + (b \cdot a ) b $. Thus,
using \eqref{eq:initialeulerlag},
we
solve for ${\rm Re}[\hat v_* ] $ in the equation
\begin{align*}
	\pi \left(\vert \xi \vert^2 {\rm Re}[\hat v_* ]  + (\xi \cdot {\rm Re}[\hat v_* ]  ) \xi \right) =  {\rm Im} [ \hat \chi ]    \xi  .
\end{align*}
We take the scalar product of the previous equation with $\vert \xi \vert^{-2} \xi$ and find that
\begin{align*}
	  \pi \xi \cdot {\rm Re}[\hat v_* ]   =  \frac{1}{2} \vert \xi \vert^{-2} {\rm Im} [ \hat \chi ]   \xi  \cdot \xi  .
\end{align*}
A combination of the two previous identities then provides an explicit expression for $ {\rm Re}[\hat v_* ]$. We can argue analogously for the minimization of the
second term in \eqref{eq:complexexpansion} which leads to 
\begin{align}\label{eq:identityminimizer}
	 \hat v_*   =  i  \left(     (2\pi)^{-1}    \vert \xi \vert^{-4}  \xi \cdot \hat \chi   \xi   - \pi^{-1} \vert \xi \vert^{-2} \hat \chi   \right) \xi.
\end{align}
Next, we compute the value of this minimizer. As a preparation
we recall the identities 
$ 2 \vert (a \odot b ) \vert^2 = \vert a \vert^2 \vert b \vert^2  + (a \cdot b)^2 $,
$ \vert a \odot a \vert^2 = \vert a \vert^4$,
$( a \odot b ) : ( b \otimes b ) = (a \cdot b) \vert b \vert^2$,
$(Ca) \odot a : C = \vert Ca \vert^2$,
$a \otimes a : C = a \cdot Ca$
for vectors $a,b \in \R^d$, $C \in \R^{d \times d}_{\sym}$.
 The identity \eqref{eq:identityminimizer} and these identities (employed in this respective order) imply that
\begin{align}\label{eq:tediousexp1}
	\vert 2\pi \hat v_* \odot   i \xi - \hat \chi  \vert^2 & = \big\vert 2 \vert \xi \vert^{-2} (\hat \chi  \xi)   \odot \xi - \vert \xi \vert^{-4} ( \xi \cdot \hat \chi \xi )  (\xi \odot \xi ) - \hat \chi \big\vert^2  \notag \\
	& = \big\vert 2 \vert \xi \vert^{-2} (\hat \chi  \xi)  \odot \xi  \big\vert^2 + \big\vert \vert \xi \vert^{-4} ( \xi \cdot \hat \chi \xi )  (\xi \odot \xi )   \big\vert^2  +  \vert \hat \chi  \vert^2  \notag \\
	& \quad \qquad - 4   \vert \xi \vert^{-2} (\hat \chi  \xi)  \odot \xi    :   \vert \xi \vert^{-4} ( \xi \cdot \hat \chi \xi )  (\xi \odot \xi )   -  4  \vert \xi \vert^{-2} (\hat \chi  \xi)  \odot \xi    :  \hat \chi  \notag \\
	&\quad \qquad + 2\vert \xi \vert^{-4} ( \xi \cdot \hat \chi \xi )  (\xi \odot \xi ) :\hat \chi \notag \\
	& =   2 \vert \xi \vert^{-4} \left(  \vert\hat \chi  \xi \vert^2  \vert \xi\vert^2 + (\xi \cdot \hat \chi  \xi )^2 \right) +   \vert \xi \vert^{-4}  \vert   ( \xi \cdot \hat \chi \xi )     \vert^2  + \vert \hat \chi \vert^2 \notag \\
	& \quad \qquad - 4   \vert \xi \vert^{-4}       ( \xi \cdot \hat \chi \xi )^2     -  4  \vert \xi \vert^{-2}   \vert \hat \chi \xi \vert^2 \notag \\
	&\quad \qquad + 2\vert \xi \vert^{-4} ( \xi \cdot \hat \chi \xi )^2  \notag \\
	 & =    \vert \xi \vert^{-4} \vert   ( \xi \cdot \hat \chi \xi )     \vert^2  + \vert \hat \chi  \vert^2   -  2  \vert \xi \vert^{-2}   \vert \hat \chi  \xi \vert^2 .
\end{align}
By the diagonal structure of $\hat \chi$, a further elementary computation yields  
\begin{align}\label{eq:tediousexp2}
	&\quad  \ \vert \xi \vert^{-4} \vert   ( \xi \cdot \hat \chi \xi )     \vert^2  + \vert \hat \chi (\xi) \vert^2   -  2   \vert \xi \vert^{-2} \vert \hat \chi (\xi) \xi \vert^2 \notag  \\
	& =   \vert \xi \vert^{-4} \left\vert \sum_{i = 1}^d \hat \chi_{ii} \xi_i \xi_i \right\vert^2 +    \sum_{i=1}^d \vert \hat \chi_{ii} \vert^2 - 2 \vert \xi \vert^{-2} \sum_{i =1}^d \vert \hat \chi_{ii} \xi_i \vert^2  \notag \\
	& =   \vert \xi \vert^{-4} \Bigg( \sum_{i = 1}^d  \xi_i^2 \vert \hat \chi_{ii}   \vert^2 +  \sum_{i = 1}^d \sum_{j \neq i}^d \xi_i^2 \xi_j^2 {\rm Re} [\hat \chi_{ii} \overline{\hat \chi_{jj}} ]  \notag \\
	& \qquad \qquad  +  \sum_{i=1}^d \sum_{m=1}^d \sum_{l \neq m}^d    \xi_m^2  \xi_l^2    \vert \hat \chi_{ii} \vert^2  
	      +  \sum_{i=1}^d \sum_{m=1}^d   \xi_m^4    \vert \hat \chi_{ii} \vert^2 - 2 \sum_{i =1}^d \sum_{m=1}^d   \xi_m^2   \vert \hat \chi_{ii} \vert^2   \xi_i  ^2 \Bigg) \notag \\
	& =  \vert \xi \vert^{-4} \left( \sum_{i=1}^d \sum_{j \neq i}^d    \xi_i^2 \xi_j^2   {\rm Re} [\hat \chi_{ii} \overline{\hat \chi_{jj}} ]   +\sum_{i = 1}^d \Bigg( \sum_{j \neq i}^d    \xi_j^2 \Bigg)^2 \vert \hat \chi_{ii} \vert^2   \right)   \notag  \\
	& = \vert \xi \vert^{-4} \left( \sum_{i = 1}^d \sum_{j = i+1}^d \vert  \xi_i^2 \hat \chi_{jj} +  \xi_j^2 \hat \chi_{ii} \vert^2 +  \sum_{i=1}^d  \sum_{j \neq i }^d \sum_{k \neq j , k \neq i}^d  \xi_k^2 \xi_j^2 \vert \hat \chi_{ii} \vert^2 \right).
\end{align}
Eventually, \eqref{eq:writeinFour}, \eqref{eq:tediousexp1}, and \eqref{eq:tediousexp2} yield \eqref{ineq:foursym}.
Let us comment on the case $u \in \mathcal{A}_F^{\per}$. Setting
$v(x) = u(x) - Fx$, \eqref{def:periodicboundarycond} implies that $\langle \nabla v \rangle = 0$.
Then, the fundamental theorem of calculus implies that $v \in H^1(\mathbb{T}^d, \R^d)$. Hence, we can proceed analogously to the case $u \in \mathcal{A}_F^{\dir}$, concluding the proof.
\end{proof}

\subsection*{Acknowledgements}
L.M.\ and A.R.\ gratefully acknowledge support by the Deutsche Forschungsgemeinschaft (DFG, German Research Foundation)
under Germany's Excellence Strategy -- EXC-2047/1 -- 3906858. A.R.~further also gratefully acknowledges support by the Deutsche Forschungsgemeinschaft (DFG, German Research Foundation) through SPP2256, project ID 44106824.

The authors gratefully thank Sergio Conti and Antonio Tribuzio for helpful discussions.
 Some of the figures were created with the assistance of ChatGPT and subsequently revised and suitably adapted. 
  
\printbibliography

\end{document}